\documentclass[11pt]{amsart}
\usepackage[letterpaper,margin=1.05in]{geometry}
\usepackage[T1]{fontenc}
\usepackage{lmodern,microtype}
\usepackage{amsmath,amssymb,amscd,mathtools,enumitem,booktabs,array}
\usepackage{xcolor}
\usepackage[colorlinks=true,linkcolor=blue!45!black,citecolor=blue!45!black,urlcolor=blue!45!black,bookmarksnumbered=true,bookmarksdepth=2]{hyperref}
\hypersetup{pdftitle={Solid categorical completion and Serre duality},pdfauthor={Yuan Gao}}
\numberwithin{equation}{section}
\newtheorem{theorem}{Theorem}[section]
\newtheorem{proposition}[theorem]{Proposition}
\newtheorem{lemma}[theorem]{Lemma}
\newtheorem{corollary}[theorem]{Corollary}
\theoremstyle{definition}
\newtheorem{definition}[theorem]{Definition}

\newtheorem{example}[theorem]{Example}
\theoremstyle{remark}
\newtheorem{remark}[theorem]{Remark}
\newcommand{\Z}{\mathbb Z}
\newcommand{\F}{\mathbb F}
\newcommand{\Q}{\mathbb Q}
\newcommand{\N}{\mathbb N}
\newcommand{\DZ}{\mathcal D(\mathbb Z)}
\newcommand{\DS}{\mathcal D(\mathbb Z_{\square})}
\newcommand{\Hom}{\operatorname{Hom}}
\newcommand{\sHom}{\underline{\operatorname{Hom}}}
\newcommand{\Mod}{\operatorname{Mod}}
\newcommand{\LMod}{\operatorname{LMod}}
\newcommand{\RMod}{\operatorname{RMod}}
\newcommand{\Fun}{\operatorname{Fun}}
\newcommand{\Perf}{\operatorname{Perf}}
\newcommand{\Ind}{\operatorname{Ind}}
\newcommand{\Cauchy}{\operatorname{Cauchy}}
\newcommand{\Kar}{\operatorname{Kar}}
\newcommand{\sqdual}{\mathbb D_{\square}}
\newcommand{\ordinarydual}{\mathbb D_k}
\newcommand{\HH}{\operatorname{HH}}
\newcommand{\Map}{\operatorname{Map}}

\newcommand{\PsPerf}{\operatorname{PsPerf}}

\newcommand{\cofib}{\operatorname{cofib}}
\newcommand{\fib}{\operatorname{fib}}
\newcommand{\id}{\mathrm{id}}
\newcommand{\op}{\mathrm{op}}
\newcommand{\qc}{\mathrm{qc}}
\newcommand{\thick}{\operatorname{thick}}
\DeclareMathOperator*{\colim}{colim}
\newcommand{\doi}[1]{\href{https://doi.org/#1}{\nolinkurl{#1}}}
\allowdisplaybreaks[1]
\setlist[enumerate]{leftmargin=*,itemsep=2pt,topsep=5pt}
\setlist[enumerate]{label=\textup{(\roman*)}}
\title[Solid completion and Serre duality]{Solid categorical completion and Serre duality}
\author{Yuan Gao}
\address{School of Mathematics, Nanjing University, Nanjing 210093, China}
\email{yuangao@nju.edu.cn}
\urladdr{}
\date{29 November 2026}
\begin{document}

\begin{abstract}
We construct a solid refinement of categorical formal completion that records profinite families of morphisms in that completion.
For a smooth proper scheme over $\Z$ or $\F_p$, we prove that enriched Cauchy completion recovers all perfect complexes on the formal neighborhood of any closed subscheme, together with their inverse-limit mapping objects.
We also establish Serre duality on the categorical puncture.
If a proper linear category over either ring has a coherent Serre autoequivalence preserving the full thick subcategory along which we complete, then the induced autoequivalence, shifted by $-1$, is an internal solid Serre autoequivalence of the puncture and its enriched Cauchy completion.
Our examples include unbounded mapping objects and show that the corresponding statement with derived algebraic duals can fail.
For a smooth proper $\F_p$-scheme, the underlying category of the enriched Cauchy completion of the puncture is equivalent to the category of all perfect complexes on the geometric puncture when the support is an effective Cartier divisor or a smooth closed subscheme.
\end{abstract}

\maketitle
\tableofcontents

\section{Introduction}
\label{sec:introduction}

Formal completion studies a space through the infinitesimal neighborhoods of a closed locus.
For a commutative noetherian ring $R$ and an ideal $I$, its ring of formal functions is $\widehat R_I=\varprojlim_nR/I^n$.
This operation also has a categorical description: completing $\Perf(R)$ along the perfect complexes supported on $V(I)$ changes the endomorphism algebra of the object $R$ to $\widehat R_I$.
More generally, categorical formal completion records an object through its morphisms from a chosen full subcategory and then takes morphisms between the resulting restricted Yoneda modules.
Efimov introduced this construction through derived double centralizers \cite{EfimovCompletion}; its restricted-Yoneda and module formulations are developed in \cite[Appendix~B]{BGW} and \cite{Gao}.

The completed ring carries an inverse-limit topology as well as an algebraic structure.
For example, the discrete and $p$-adic topologies on $\Z_p$ have the same underlying abelian group but determine different continuous families parametrized by profinite sets.
Solid modules provide a setting in which these families enter the algebraic formalism \cite{Scholze,ScholzeAnalytic}.
We apply the discrete embedding to the evaluation modules and compute their internal Homs over the solid coefficient category.
Evaluation at a point recovers categorical completion over the original coefficient ring, while sections over profinite sets record families of morphisms in that completion.
This additional structure also affects duality: mapping objects in the categorical puncture can satisfy internal solid Serre duality even when the corresponding statement for their derived algebraic duals fails.

The categorical constructions have close geometric counterparts.
Efimov relates categorical completion to perfect complexes on formal neighborhoods; for smooth proper schemes over a field, he describes all formal perfect complexes in terms of modules \cite[Theorem~5.4]{EfimovCompletion}, \cite[Theorem~2.14]{EfimovInfinity}.
For geometric boundaries, Pantev--To\"en establish Tate duality on boundary cohomology \cite[Proposition~5.5 and Corollary~5.6]{PT}.
The solid geometric duality of Clausen--Scholze identifies the dualizing complex of a smooth morphism of relative dimension $d$ with $\det\Omega^1[d]$ \cite[Theorem~7.7.21 and Corollary~7.7.23]{Camargo}.
These results motivate two questions about the solid categorical construction: which objects it contains after Cauchy completion, and how Serre duality behaves after quotienting by the subcategory along which completion is formed.

We answer these questions over $\Z$ and $\F_p$.
For a smooth proper scheme, Cauchy completion recovers every perfect complex on the formal neighborhood of any closed subscheme, together with its inverse-limit mapping objects.
For a proper linear category over either ring with a coherent Serre autoequivalence, we consider the invariant full stable subcategories of pointwise perfect modules specified below.
We prove internal solid Serre duality shifted by $-1$ on their quotients by the representable modules and on the Cauchy completions of these quotients.
Over $\F_p$, when the support is an effective Cartier divisor or a smooth closed subscheme, we identify the underlying category of the enriched Cauchy completion of the categorical puncture with the category of all perfect complexes on the geometric puncture.

\subsection{Solid completion and formal perfect complexes}

For a discrete commutative ring $k$, the relative coefficient category is
\[
 \mathcal D((k,\Z)_\square)=\Mod_{\delta k}(\DS),
 \qquad
 \delta:\mathcal D(k)\rightleftarrows
 \mathcal D((k,\Z)_\square):\Gamma.
\]
Here $\delta$ is the fully faithful symmetric monoidal discrete embedding and $\Gamma$ is evaluation at a point.
We call the image under $\delta$ of any complex in $\mathcal D(k)$ discrete, in every cohomological amplitude.
The notation $\DS$ denotes $\mathcal D(\operatorname{Solid}_{\Z})$.

Let $\mathcal C$ be small, stable, idempotent complete, and $\Perf(k)$-linear, and let $\mathcal T\subset\mathcal C$ be full and thick.
For $X\in\mathcal C$, write $r_X=\mathcal C(-,X)|_{\mathcal T}$.
The solid categorical completion has the objects of $\mathcal C$ and morphism objects
\begin{equation}\label{eq:intro-support-completion}
 \widehat{\mathcal C}_{\mathcal T}^{\square}(X,Y)
 =\sHom_{\delta\mathcal T}(\delta r_X,\delta r_Y).
\end{equation}
Composition and units are inherited from the solid module category.
Propositions~\ref{thm:solid-enrichment} and~\ref{thm:family-formula} identify its point evaluation and profinite sections with
\begin{align}
 \Gamma\widehat{\mathcal C}_{\mathcal T}^{\square}(X,Y)
 &\simeq\Hom_{\mathcal T}(r_X,r_Y),
 \label{eq:intro-underlying}\\
 \widehat{\mathcal C}_{\mathcal T}^{\square}(X,Y)(S)
 &\simeq\Hom_{\mathcal T}(r_X,C(S,k)\otimes_kr_Y),
 \label{eq:intro-families}
\end{align}
where $C(S,k)$ is the ring of locally constant functions on a profinite set $S$.
In \S~\ref{sec:completion}, we construct completion over a general locally rigid coefficient category.
The profinite-family formula~\eqref{eq:intro-families} follows by specializing to solid coefficients.

Let $\mathcal B$ be a small category enriched over $\mathcal D((k,\Z)_\square)$.
A right $\mathcal B$-module is an enriched functor $\mathcal B^{\op}\to\mathcal D((k,\Z)_\square)$.
\begin{samepage}
Such a module $P$ is Cauchy if the coefficient-valued internal Hom functor
\[
 \sHom_{\mathcal B}(P,-):
 \RMod_{\mathcal B}\bigl(\mathcal D((k,\Z)_\square)\bigr)
 \longrightarrow\mathcal D((k,\Z)_\square)
\]
preserves small colimits and is linear over $\mathcal D((k,\Z)_\square)$, meaning that it commutes with tensoring by coefficient objects.
\end{samepage}
The Cauchy right modules form a full enriched subcategory, the Cauchy completion $\Cauchy(\mathcal B)$.
Yoneda embeds $\mathcal B$ in $\Cauchy(\mathcal B)$ by $b\mapsto\mathcal B(-,b)$.

Equivalently, the representables generate $\Cauchy(\mathcal B)$ under absolute weighted colimits, those preserved by every functor enriched over the given base \cite[\S~5.5]{Kelly}, \cite[Proposition~5.49 and Corollary~5.54]{HeineWeighted}.
This completion contains the closure of the representables under finite cofibers and retracts, and over solid coefficients the inclusion can be strict.
In \S~\ref{subsec:cauchy}, we identify the condition in Definition~\ref{def:cauchy-module} with preservation of the Yoneda presentation~\eqref{eq:absolute-yoneda} by every enriched functor out of the module category.

Let $k=\Z$ or $\F_p$, let $X$ be smooth and proper over $k$, and let $Z\subseteq X$ be a closed subscheme defined by an ideal sheaf $\mathcal I$.
Write $i_n:Z_n=V(\mathcal I^n)\hookrightarrow X$.
The category $\Perf(\widehat X_Z)=\lim_n\Perf(Z_n)$ consists of compatible formal perfect complexes, and $\Perf_Z(X)$ consists of perfect complexes whose restriction to $X\setminus Z$ vanishes.
We equip $\Perf(\widehat X_Z)$ with the solid mapping objects
\begin{equation}\label{eq:intro-formal-mappings}
 \sHom_{\widehat X_Z}(L,M)
 =\varprojlim_n\delta\Hom_{Z_n}(L_n,M_n).
\end{equation}

\begin{theorem}
\label{thm:intro-cauchy}
Let $k=\Z$ or $\F_p$, let $X/k$ be smooth and proper, and let $Z\subseteq X$ be a closed subscheme.
There is a canonical equivalence of $\mathcal D((k,\Z)_\square)$-enriched categories
\begin{equation}\label{eq:intro-cauchy}
 \Cauchy\!\left(\widehat{\Perf(X)}^{\square}_{\Perf_Z(X)}\right)
 \simeq\Perf(\widehat X_Z).
\end{equation}
The formal mapping objects are those of~\eqref{eq:intro-formal-mappings}.
For $P\in\Perf(X)$, the equivalence sends the image of $P$ under the completion and Yoneda functors to the compatible system $(i_n^*P)_n$.
Every compatible formal perfect complex occurs, including those outside the thick closure of restrictions when such complexes exist.
\end{theorem}

We begin the proof by constructing the solid realization $\lim_n\delta(i_{n*}L_n)$ of a formal perfect complex.
For formal perfect complexes $L$ and $M$, the map~\eqref{eq:geometric-tensor-map} of Lemma~\ref{lem:geometric-tensor} compares the solid tensor product of their realizations with the inverse limit of $i_{n*}L_n\otimes_{\mathcal O_X}i_{m*}M_m$.
Lemma~\ref{lem:formal-pro-multiplication} shows that the map~\eqref{eq:formal-pro-product} identifies the corresponding pro-object with $\{i_{n*}(L_n\otimes_{\mathcal O_{Z_n}}M_n)\}_n$.
These identifications, including the case of the formal structure sheaf, prove that the realization is symmetric monoidal.
Conversely, using the module equivalence~\eqref{eq:geometric-free-envelope}, the proof of Theorem~\ref{thm:geometric-cauchy} identifies each Cauchy module over the solid categorical completion with a derived $\mathcal I$-complete solid object on $X$.
Its reductions $L_n$ are compatible perfect complexes on the thickenings $Z_n$, and the inverse limit $\lim_n\delta(i_{n*}L_n)$ recovers that realization.
Theorem~\ref{thm:geometric-cauchy} proves that these constructions are mutually inverse.

The solid categorical completion $\widehat{\Perf(X)}^{\square}_{\Perf_Z(X)}$ has objects indexed by $\Perf(X)$; under Yoneda followed by~\eqref{eq:intro-cauchy}, these objects map to restrictions of perfect complexes from $X$.
The closure of their images under finite cofibers and retracts is the thick closure of these restrictions.
The enriched Cauchy completion on the left-hand side of~\eqref{eq:intro-cauchy} identifies with the entire category $\Perf(\widehat X_Z)$.
This distinction between the thick closure of restrictions and the full formal category also appears in the earlier geometric descriptions.
Efimov's double centralizer recovers the thick closure of restrictions from $X$ \cite[Theorem~5.4]{EfimovCompletion}, whereas his description by pseudo-perfect $\Perf_Z(X)$-modules recovers all formal perfect complexes over a field \cite[Theorem~2.14]{EfimovInfinity}.
Here pseudo-perfect means that the module takes values in $\Perf(k)$ on every object of $\Perf_Z(X)$.
His lecture announces the latter description for smooth separated $X$ of finite type over a noetherian base with support proper over that base \cite[slide~7/35]{EfimovNuclear}.
For $k=\Z$ or $\F_p$, Theorem~\ref{thm:intro-cauchy} characterizes the formal objects by the solid Cauchy condition and recovers their mapping objects~\eqref{eq:intro-formal-mappings}.

\subsection{Internal Serre duality on the puncture}

Let $\mathcal H$ be the stable idempotent-complete envelope of $\widehat{\mathcal C}_{\mathcal T}^{\square}$ inside the solid module category.
For $t\in\mathcal T$, Yoneda over $\mathcal D((k,\Z)_\square)$ identifies the map of the completion functor~\eqref{eq:support-completion-functor} with
\[
 \delta\mathcal C(t,X)\xrightarrow{\ \sim\ }\mathcal H(t,X).
\]
The completion functor therefore restricts to a fully faithful exact inclusion $\delta\mathcal T\hookrightarrow\mathcal H$.
Define
\begin{equation}\label{eq:intro-puncture-orders}
 \mathcal Q_0=\mathcal H/\delta\mathcal T,
 \qquad
 \mathcal Q=\Cauchy(\mathcal Q_0),
 \qquad
 \mathcal Q_{\mathrm{form}}
 =\Kar\bigl(\Cauchy(\mathcal H)/\delta\mathcal T\bigr),
\end{equation}
where $\Kar$ denotes stable idempotent completion.
The quotient is formed with the internal solid mapping objects already in place.
In Proposition~\ref{prop:cauchy-quotient}, we construct a canonical fully faithful functor $\iota_{\mathrm{form}}:\mathcal Q_{\mathrm{form}}\to\mathcal Q$.
Theorem~\ref{thm:intro-puncture-geometry} proves that this functor is an equivalence in the stated geometric cases.

For $M\in\mathcal D((k,\Z)_\square)$, write
\[
 \sqdual M=\sHom_{\delta k}(M,\delta k),
 \qquad
 \ordinarydual E=\Hom_k(E,k)\quad(E\in\mathcal D(k)).
\]
A $k$-linear category is proper when all its mapping complexes are perfect over $k$.
A coherent Serre structure on an autoequivalence $S$ is an equivalence $\ordinarydual\mathcal C(X,Y)\simeq\mathcal C(Y,SX)$ of $\mathcal D(k)$-enriched bimodules.
In particular, it includes the compatibility with composition and its coherent homotopies.

\begin{theorem}\label{thm:intro-punctured-serre}
Let $k=\Z$ or $\F_p$.
Suppose that $\mathcal C$ is a small stable idempotent-complete proper $k$-linear category with a coherent $\mathcal D(k)$-enriched Serre autoequivalence $S$, and that $\mathcal T\subset\mathcal C$ is full and thick with $S(\mathcal T)=\mathcal T$.
Then $S$ induces an autoequivalence $\overline S$ of $\mathcal Q$, and there is an equivalence of $\mathcal D((k,\Z)_\square)$-valued bimodules
\begin{equation}\label{eq:intro-punctured-serre}
 \sqdual\mathcal Q(P,R)
 \simeq\mathcal Q(R,\overline S P)[-1].
\end{equation}
The equivalence holds on $\mathcal Q_0$, extends to every Cauchy object over $\mathcal D((k,\Z)_\square)$, and restricts to $\mathcal Q_{\mathrm{form}}$.
Thus the internal solid Serre autoequivalence is $\overline S[-1]$.
\end{theorem}

In the proof, we use the Serre structure on $\mathcal T$ and pointwise perfection of the restricted modules.
Theorem~\ref{thm:punctured-solid-serre} therefore treats a full stable subcategory of pointwise perfect $\mathcal T$-modules containing the representables and invariant under precomposition by $S_{\mathcal T}^{-1}$; it requires properness only of $\mathcal T$.
For Theorem~\ref{thm:intro-punctured-serre}, the full stable idempotent-complete subcategory generated by the restricted modules $r_X$ has these properties by properness of $\mathcal C$ and the condition $S(\mathcal T)=\mathcal T$.
In the Calabi--Yau case the induced action is a shift, so invariance follows from stability of the subcategory.

The proof uses the bimodule $I$ defined in~\eqref{eq:punctured-factorization-coend} by the coend
\[
 I(M,N)=\int^{t\in\delta\mathcal T}
 \mathcal H(M,h_t)\otimes_{\delta k}^{\square}\mathcal H(h_t,N),
 \qquad h_t=\delta\mathcal T(-,t).
\]
Lemma~\ref{lem:serre-mixed-maps} identifies $I(M,N)$ with the discrete realization in~\eqref{eq:serre-discrete-coend}, and Lemma~\ref{lem:unbounded-solid-biduality} proves canonical biduality for these complexes, including the unbounded case.
For $k=\Z$, the proof uses length-one free resolutions and the decomposition~\eqref{eq:hereditary-complex-splitting}; the case $k=\F_p$ follows by restriction and coinduction from the integral case.
Equation~\eqref{eq:serre-diagonal-dual-coend} identifies the mapping objects of $\mathcal H$ with internal solid duals of values of $I$.
The trace identity~\eqref{eq:testing-serre-cyclicity} makes the resulting duality equivalences compatible with composition, as expressed by the commuting square of Lemma~\ref{lem:serre-transpose-square}.
Applying this square to the cofiber sequence~\eqref{eq:punctured-quotient-triangle}, we obtain~\eqref{eq:serre-quotient-transpose}, and hence duality on $\mathcal Q_0$ with the shift $[-1]$ in~\eqref{eq:intro-punctured-serre}.
For a Cauchy module, every enriched functor preserves the Yoneda presentation~\eqref{eq:absolute-yoneda}.
Using this property, the calculation in~\eqref{eq:cauchy-serre-tensor-hom}--\eqref{eq:cauchy-solid-serre} extends the duality to $\mathcal Q$.

Related categorical boundary pairings arise from Calabi--Yau structures.
The chain-level residue in \cite[Definition~1.1 and Theorem~1.7]{GaoResidue} induces pairings on morphism cohomology from a smooth Calabi--Yau structure.
Hanihara--Liu construct the quotient bimodule map from a cofiber and identify it with the Hochschild connecting trace.
Their criterion for nondegeneracy requires biduality of the ambient diagonal and its recovery from Homs between the restricted modules \cite[\S~4.2, Proposition~4.2.2 and Theorem~4.3.1(a),(b)]{HL}.
For a category with a weak proper Calabi--Yau structure over a field, Gao--Gao obtain algebraic-dual nondegeneracy when every cohomology group of every Hom between pseudo-perfect modules is finite-dimensional, and exhibit its failure without finiteness \cite[Theorem~5.6 and Example~5.7]{GaoGao}.
Here we construct the quotient pairing on the full solid mapping objects through~\eqref{eq:serre-quotient-transpose}, using Lemmas~\ref{lem:unbounded-solid-biduality} and~\ref{lem:serre-mixed-maps} and the compatibility with composition in Lemma~\ref{lem:serre-transpose-square}.
The distinction between bounded and unbounded objects also appears in the realization of derived Tate objects: for rings of finite type over $\Z$, Melani--Pourcelot--Vezzosi prove full faithfulness for bounded countable derived Tate objects and give a counterexample to fullness in the unbounded case \cite[Theorem~4.3 and Proposition~5.1]{MPV}.

The examples illustrate why the solid mapping objects are essential to this duality.
In Example~\ref{ex:unbounded-frobenius}, we construct a quotient endomorphism complex with infinitely generated cohomology in arbitrarily positive and negative degrees.
In Example~\ref{ex:product-quotient}, the internal solid dual of the quotient endomorphism algebra is its shift by $-1$.
The point evaluation of the algebra itself is a nonzero vector space concentrated in degree zero, so its algebraic dual cannot satisfy the corresponding shifted self-duality.

\subsection{The full geometric puncture}

Return to the smooth proper scheme $X/k$ and closed subscheme $Z$ above, and specialize $\mathcal C=\Perf(X)$ and $\mathcal T=\Perf_Z(X)$.
Suppose that $X$ has pure relative dimension $d$.
Its Serre functor $(-)\otimes_{\mathcal O_X}\omega_{X/k}[d]$ preserves $\Perf_Z(X)$ for every closed $Z$.
Tensoring by $\omega_{X/k}$ therefore induces autoequivalences of $\mathcal Q_0=\mathcal H/\delta\Perf_Z(X)$ and $\mathcal Q=\Cauchy(\mathcal Q_0)$, denoted by the same tensor symbol.
By Theorem~\ref{thm:intro-punctured-serre}, we have
\begin{equation}\label{eq:intro-geometric-serre}
 \sqdual\mathcal Q(P,R)
 \simeq\mathcal Q(R,P\otimes_{\mathcal O_X}\omega_{X/k})[d-1]
\end{equation}
on the categorical solid puncture.
Theorem~\ref{thm:intro-cauchy} identifies $\Cauchy(\mathcal H)$ with all formal perfect complexes.
The next theorem shows that $\iota_{\mathrm{form}}:\mathcal Q_{\mathrm{form}}\to\mathcal Q$ is an equivalence when $Z$ is an effective Cartier divisor or a smooth closed subscheme over $\F_p$.

\begin{theorem}\label{thm:intro-puncture-geometry}
Let $X$ be smooth and proper over $\F_p$, and let $Z\subset X$ be an effective Cartier divisor or a smooth closed subscheme.
The canonical $\mathcal D((\F_p,\Z)_\square)$-enriched functor
\begin{equation}\label{eq:intro-geometric-puncture}
 \iota_{\mathrm{form}}:
 \Kar\!\left(\Perf(\widehat X_Z)/\delta\Perf_Z(X)\right)
 \longrightarrow
 \Cauchy\!\left(\mathcal H/\delta\Perf_Z(X)\right)
\end{equation}
is an equivalence, where $\mathcal H$ is the stable idempotent-complete envelope of $\widehat{\Perf(X)}^{\square}_{\Perf_Z(X)}$ and the formal category has the enrichment~\eqref{eq:intro-formal-mappings}.
After point evaluation, we obtain an equivalence of $\F_p$-linear stable categories
\[
 \Gamma\mathcal Q\simeq\Perf(\widehat X_Z\setminus Z),
\]
where the right-hand side is the category of all perfect complexes on the punctured formal neighborhood of Hennion--Porta--Vezzosi.
Via $\iota_{\mathrm{form}}$ in~\eqref{eq:intro-geometric-puncture}, the same identification holds for $\Gamma\mathcal Q_{\mathrm{form}}$.
\end{theorem}

The classical quotient description of the geometric puncture comes from coherent localization \cite[Theorem~7.3 and proof of Corollary~7.5]{HPV}, together with regularity.
To prove that $\iota_{\mathrm{form}}$ in~\eqref{eq:intro-geometric-puncture} is an equivalence, we identify the Cauchy objects over $\mathcal D((\F_p,\Z)_\square)$ with dualizable modules over the solid algebra defined in~\eqref{eq:punctured-structure-algebra}.
For a principal-adic affine chart, Andreychev's analytic completeness and perfect-complex theorems identify these modules with perfect complexes over the ring of point sections \cite[Lemma~5.27 and Corollary~5.51.1]{Andreychev}.
For an effective Cartier divisor, we prove the equivalence by a finite \v Cech calculation on the punctured completions of affine charts.
For a smooth center, we reduce to the Cartier case by blowing up and verify the pushforwards of all exceptional ideal powers, together with compatibility with formal perfect complexes that are not algebraizable.
We carry out these arguments in the proof of Theorem~\ref{thm:geometric-puncture}.
For smooth proper $X$ of pure relative dimension $d$ over $\Z$ or $\F_p$, formula~\eqref{eq:intro-geometric-serre} expresses internal duality on the categorical puncture for every closed $Z\subseteq X$.

We also compare Hochschild cochains.
Restriction along the Yoneda inclusion $\delta\mathcal T\hookrightarrow\mathcal H$ induces an $E_2$-equivalence
\begin{equation}\label{eq:intro-center-restriction}
 \HH^\bullet_\square(\mathcal H)
 \longrightarrow\HH^\bullet_\square(\delta\mathcal T).
\end{equation}
This follows from evaluation on representables for any full module subcategory containing them (Proposition~\ref{prop:restriction-centers}).
For a commutative ring $R$ essentially of finite type over $\Z$, an ideal $I\subseteq R$ and $A=\lim_n\delta(R/I^n)$, the scalar-extension map $A\otimes_{\delta R}^{\square}\delta\HH^\bullet(R/\Z)\to\HH^\bullet_\square(A/\Z)$ is an equivalence of $E_1$ $A$-algebras (Proposition~\ref{thm:hochschild-comparison}).

\subsection*{Conventions}
All derived and module categories are stable $\infty$-categories.
For a ring or algebra object $R$, $\mathcal D(R)$ denotes the unbounded derived category of left modules; the restricted Yoneda modules over $\mathcal T$ in the support-completion and duality constructions are right modules.
We use cohomological grading, so connective algebras have no positive cohomology.
All Hom objects, tensor products, ends, coends, limits, Kan extensions, and quotients are derived.
We omit $\mathrm R$ and $\mathbf L$ from their notation.
We write $\Hom$ for derived mapping complexes over the stated discrete base ring, $\sHom$ for internal mapping objects over the declared coefficient category, and $\Map$ for mapping spaces.
Subscripts specify the category or module argument.
In particular, global $\Hom_X$ takes values in $\mathcal D(k)$ for a $k$-scheme $X$.
On a category enriched over $\mathcal D((k,\Z)_\square)$, $\Gamma$ evaluates the mapping objects at a point and induces the $\mathcal D(k)$-enriched composition through its lax monoidal structure.
We write $\otimes_\Z$ for complexes in $\DZ$, $\otimes_{\Z_\square}$ in $\DS$, and $\otimes_A^\square$ for a relative solid tensor product.
The quotient and Cauchy completions use compatible universes, enlarged together when necessary.

\subsection*{Structure of the paper}

In \S~\ref{sec:completion}, we construct completion by evaluation and realize it by enriched Kan extension.
In \S~\ref{sec:solid}, we describe profinite families and Cauchy completion, define the punctured category, and compare quotienting before and after Cauchy completion.
In \S~\ref{sec:punctured-duality}, we prove internal solid Serre duality by computing the coend~\eqref{eq:punctured-factorization-coend}, transposing its composition map, and extending the resulting pairing to Cauchy objects.
In \S\S~\ref{sec:adic} and~\ref{sec:scalar-extension}, we compute morphism objects of adic categorical completions and compare derived completion with solid scalar extension.
In \S~\ref{sec:formal-geometry}, we identify perfect complexes on formal completions and punctures using enriched Cauchy completion.
In \S~\ref{sec:hochschild}, we prove the cochain comparisons, and in \S~\ref{sec:examples}, we give examples with infinite and unbounded mapping objects.
In Appendix~\ref{app:centralizers}, we examine coefficient change for internal Homs between modules.

\subsection*{Acknowledgments and AI Declaration}

The original ideas and proofs underlying this work were developed while the author was writing the paper \cite{Gao}, whose first arXiv version appeared in 2024.
The initial construction expressed the morphism objects of categorical formal completion through an explicit homotopy-end formula.
The author thanks Christopher Brav, Yingdi Qin, and Umut Varolgunes for helpful discussions.

The author used ChatGPT 6 and 6.2 Astra to assist with revising and polishing the exposition, identifying and checking relevant references, and developing some arguments---in particular, to recast the explicit homotopy-end argument in terms of enriched internal Homs between evaluation modules over solid coefficients.
The author reviewed the suggestions provided by the LLM, revised the final presentation, and takes full responsibility for the mathematical content of the manuscript.

\section{Categorical completion}\label{sec:completion}

A full subcategory of modules determines a completion through internal Homs between its evaluation modules.
We construct this completion over a general coefficient category and then study its behavior under change of coefficients.
For completion along a full subcategory $\mathcal T\subset\mathcal C$, we use enriched Kan extension to realize objects of the completion as modules on $\mathcal C$; this realization also determines a localization of the free module category of the completion.

\subsection{Source categories and modules}\label{subsec:source-categories}

Throughout this subsection and \S~\ref{subsec:completion-evaluation}, $\mathcal E$ is a compactly generated presentable stable symmetric monoidal $\infty$-category whose tensor product preserves colimits separately, whose unit is compact, and whose compact objects are dualizable.
Thus $\mathcal E\simeq\Ind(\mathcal E^\omega)$, with $\mathcal E^\omega$ a small stable idempotent-complete rigid symmetric monoidal category; the unit need not generate $\mathcal E$.
These are the locally rigid hypotheses of \cite[Definition~D.7.4.1]{LurieSAG}.
They hold for $\mathcal E=\mathcal D(k)$ for every commutative ring $k$, with $\mathcal E^\omega=\Perf(k)$.
In these two subsections enrichment is over $\mathcal E$.

\begin{definition}\label{def:linear-category}
A small $\mathcal E^\omega$-linear category is a small stable idempotent-complete $\infty$-category $\mathcal C$ with a coherent unital action
\[
 \mathcal E^\omega\times\mathcal C\longrightarrow\mathcal C,
 \qquad(a,X)\longmapsto a\otimes_{\mathcal E}X,
\]
exact in each variable.
Linear functors are exact functors with coherent equivalences $F(a\otimes_{\mathcal E}X)\simeq a\otimes_{\mathcal E}F(X)$; their natural transformations respect these structures.
\end{definition}

We use the convention of \cite[Definition~4.3]{HSS}, applied to $\mathcal E^\omega$.
A presentable stable $\mathcal E$-module category has a unital $\mathcal E$-action preserving colimits separately.
Extending the action to $\Ind\mathcal C$ determines the $\mathcal E$-valued mapping object characterized by
\begin{equation}\label{eq:linear-enrichment}
 \Map_{\mathcal E}\bigl(a,\sHom_{\mathcal C}(X,Y)\bigr)
 \simeq\Map_{\Ind\mathcal C}(a\otimes_{\mathcal E}X,Y),
 \qquad a\in\mathcal E.
\end{equation}
We also write $\mathcal C(X,Y)$ for this object of $\mathcal E$, which need not be compact.
For $a\in\mathcal E^\omega$, the opposite category has action $a\otimes_{\mathcal E}X^{\op}=(a^\vee\otimes_{\mathcal E}X)^{\op}$ and mapping objects $\mathcal C^{\op}(X,Y)=\mathcal C(Y,X)$.

Write $\Fun_{\mathcal E}$ for the category of $\mathcal E$-enriched functors and enriched natural transformations.
Let $\mathcal C$ be small $\mathcal E^\omega$-linear and $\mathcal D$ either small $\mathcal E^\omega$-linear or a presentable stable $\mathcal E$-module category.
By \cite[Proposition~4.5]{HSS}, an exact linear functor extends to a colimit-preserving $\mathcal E$-linear functor on $\Ind\mathcal C$, with target $\Ind\mathcal D$ in the small case.
Applying this extension to evaluation and taking tensor--Hom adjoints defines
\begin{equation}\label{eq:linear-to-enriched-functors}
 \Fun^{\mathrm{ex}}_{\mathcal E^\omega}(\mathcal C,\mathcal D)
   \longrightarrow\Fun_{\mathcal E}(\mathcal C,\mathcal D).
\end{equation}
Conversely, for $a\in\mathcal E^\omega$, an enriched functor $F$ determines a comparison $s_{a,X}:a\otimes_{\mathcal E}F(X)\to F(a\otimes_{\mathcal E}X)$, adjoint to
$a\to\mathcal C(X,a\otimes_{\mathcal E}X)\to\mathcal D(FX,F(a\otimes_{\mathcal E}X))$.

\begin{lemma}\label{lem:enriched-linear-functors}
The functor \eqref{eq:linear-to-enriched-functors} is an equivalence.
For every enriched $F$, the maps $s_{a,X}$ are equivalences, define its coherent linear structure, and make its underlying functor exact.
\end{lemma}

\begin{proof}
An inverse to $s_{a,X}$ is the composite of coevaluation, $a\otimes_{\mathcal E}s_{a^\vee,a\otimes_{\mathcal E}X}$ and $F$ of evaluation.
The triangle identities prove the inverse identities, and the action coherences follow from compatibility with enriched composition.
Taking $a=0$ and the suspension cocone of the unit shows that $F$ is reduced and its suspension comparison invertible, hence exact \cite[Corollary~1.4.2.14]{LurieHA}.
Applying $\Map_{\mathcal E}(a,-)$ for compact generators $a$ recovers the original maps on Hom objects, and the same identification on natural transformations follows from tensor--Hom adjunction.
\end{proof}

For presentable $\mathcal E$-module categories, $\Fun^L_{\mathcal E}$ denotes colimit-preserving $\mathcal E$-linear functors.
For a colimit-preserving $\mathcal E$-enriched functor between presentable $\mathcal E$-modules, the comparisons for compact scalars extend to all of $\mathcal E$ by colimits.

Let $\mathcal B$ be a small $\mathcal E$-enriched category.
A left $\mathcal B$-module is an $\mathcal E$-enriched functor $\mathcal B\to\mathcal E$, and a right $\mathcal B$-module is an $\mathcal E$-enriched functor $\mathcal B^{\op}\to\mathcal E$.
Morphisms of modules are enriched natural transformations.
We denote these categories by
\begin{align}
 \LMod_{\mathcal B}(\mathcal E)
 &:=\Fun_{\mathcal E}(\mathcal B,\mathcal E),
 \label{eq:left-enriched-modules}\\
 \RMod_{\mathcal B}(\mathcal E)
 &:=\Fun_{\mathcal E}(\mathcal B^{\op},\mathcal E).
 \label{eq:right-enriched-modules}
\end{align}
The representable left and right modules at $b\in\mathcal B$ are $\mathcal B(b,-)$ and $\mathcal B(-,b)$, respectively.
Both module categories are presentable and stable, with colimits and the $\mathcal E$-action formed objectwise.
We write $\sHom_{\mathcal B}(M,N)$ for the $\mathcal E$-valued enriched mapping object between two left modules or two right modules, as specified at each use.
When $\mathcal E=\mathcal D(k)$, we write $\Hom_{\mathcal B}(M,N)$ for the corresponding derived $k$-valued mapping complex.

Write $\thick$ for closure under finite cofibers, shifts, and retracts, and $\thick_{\mathcal E^\omega}$ for closure also under compact scalar tensors.
For $a\in\mathcal E^\omega$, the enriched Yoneda lemma implies
\[
 \Map_{\LMod_{\mathcal B}(\mathcal E)}(a\otimes_{\mathcal E}\mathcal B(b,-),N)\simeq\Map_{\mathcal E}(a,N(b)).
\]
These modules are compact generators, so the perfect left modules are
\begin{equation}\label{eq:compact-enriched-modules}
 \LMod_{\mathcal B}(\mathcal E)^\omega
 =\thick_{\mathcal E^\omega}\{\mathcal B(b,-)\}
 =\thick\{a\otimes_{\mathcal E}\mathcal B(b,-)\mid a\in\mathcal E^\omega,\ b\in\mathcal B\}.
\end{equation}
The same statements hold for right modules.
For $\mathcal E=\mathcal D(k)$, the scalar closure is the thick closure because $\Perf(k)=\thick\{k\}$.
If $\mathcal C$ is small $\mathcal E^\omega$-linear, duality identifies $a\otimes_{\mathcal E}\mathcal C(X,-)$ with $\mathcal C(a^\vee\otimes_{\mathcal E}X,-)$.
The colimit extensions of the two Yoneda embeddings are the equivalences
\[
 \Ind\mathcal C\xrightarrow{\ \sim\ }\RMod_{\mathcal C}(\mathcal E),
 \qquad
 \Ind(\mathcal C^{\op})\xrightarrow{\ \sim\ }\LMod_{\mathcal C}(\mathcal E)
\]
of \cite[Corollary~4.11]{HSS}, using Lemma~\ref{lem:enriched-linear-functors} to identify the enriched and linear module conventions.

\subsection{Completion by evaluation}\label{subsec:completion-evaluation}

Let $\mathcal C$ be small $\mathcal E$-enriched and $\mathcal M\subset\LMod_{\mathcal C}(\mathcal E)$ a small full $\mathcal E$-enriched subcategory.
The evaluation pairing
\begin{equation}\label{eq:evaluation-pairing}
 \mathcal C\otimes_{\mathcal E}\mathcal M\longrightarrow\mathcal E,
 \qquad(X,M)\longmapsto M(X)
\end{equation}
defines $e:\mathcal C\to\LMod_{\mathcal M}(\mathcal E)$ by $e_X(M)=M(X)$.

\begin{definition}\label{def:categorical-completion}
The completed category for the inclusion $\mathcal M\subset\LMod_{\mathcal C}(\mathcal E)$ has objects indexed by $X\in\mathcal C$, with mapping objects
\begin{equation}\label{eq:ordinary-completion}
 \widehat{\mathcal C}_{\mathcal M}(X,Y)=\sHom_{\mathcal M}(e_X,e_Y),
\end{equation}
with units and composition inherited from the module category.
Evaluation defines the completion functor $\mathcal C\to\widehat{\mathcal C}_{\mathcal M}$.
Its idempotent-complete stable envelope and its $\mathcal E^\omega$-linear envelope are respectively $\thick\{e_X\}$ and $\thick_{\mathcal E^\omega}\{e_X\}$ inside $\LMod_{\mathcal M}(\mathcal E)$.
\end{definition}

The definition by evaluation and internal Hom works over any presentably symmetric monoidal coefficient category with tensor product preserving colimits in each variable.
The envelope assertions here use the locally rigid coefficient category $\mathcal E$ of \S~\ref{subsec:source-categories}.
The two envelopes agree for $\mathcal E=\mathcal D(k)$, or when $\mathcal C$ is $\mathcal E^\omega$-linear, since in the latter case $a\otimes_{\mathcal E}e_X\simeq e_{a\otimes_{\mathcal E}X}$ by Lemma~\ref{lem:enriched-linear-functors}.

For a full enriched inclusion $u:\mathcal T\hookrightarrow\mathcal C$, take the full subcategory of covariant representable modules $\mathcal C(t,-)$, $t\in\mathcal T$, identified with $\mathcal T^{\op}$ by Yoneda.
Under this identification, \eqref{eq:left-enriched-modules} and~\eqref{eq:right-enriched-modules} identify left $\mathcal T^{\op}$-modules with right $\mathcal T$-modules.
Evaluation therefore becomes restricted right Yoneda,
\begin{equation}\label{eq:restricted-yoneda}
 r_{\mathcal T}:\mathcal C\longrightarrow\RMod_{\mathcal T}(\mathcal E),
 \qquad r_{\mathcal T}X=\mathcal C(-,X)|_{\mathcal T},
\end{equation}
and we write $\widehat{\mathcal C}_{\mathcal T}(X,Y)=\sHom_{\mathcal T}(r_{\mathcal T}X,r_{\mathcal T}Y)$.
This is the module convention of \cite{Gao,BGW}; the small subcategories called localizing in \cite[Appendix~B.2.1]{BGW} are called thick here.
The analogous construction for a family of right modules is the opposite of the completion of $\mathcal C^{\op}$, and the two conventions need not agree \cite[Proposition~4.4]{EfimovCompletion}.

Let $\mathcal T$ be either $\thick(\mathcal M)$ or $\thick_{\mathcal E^\omega}(\mathcal M)$ in $\LMod_{\mathcal C}(\mathcal E)$, and let $j:\mathcal M\hookrightarrow\mathcal T$ be inclusion.
Restriction carries evaluation at $X$ on $\mathcal T$ to $e_X$, defining a functor $\widehat{\mathcal C}_{\mathcal T}\to\widehat{\mathcal C}_{\mathcal M}$ that is the identity on objects.

\begin{proposition}\label{prop:testing-envelope}
Restriction along $j$ induces enriched equivalences
\[
 \begin{aligned}
 j^*:\LMod_{\mathcal T}(\mathcal E)&\xrightarrow{\ \sim\ }\LMod_{\mathcal M}(\mathcal E),\\
 j^*:\RMod_{\mathcal T}(\mathcal E)&\xrightarrow{\ \sim\ }\RMod_{\mathcal M}(\mathcal E).
 \end{aligned}
\]
For stable envelopes, with or without retracts, these module equivalences remain valid when $\mathcal E$ is any presentable stable symmetric monoidal category whose tensor product preserves colimits separately.
Consequently, restriction induces an enriched equivalence $\widehat{\mathcal C}_{\mathcal T}\to\widehat{\mathcal C}_{\mathcal M}$ under $\mathcal C$.
\end{proposition}

\begin{proof}
For a small full enriched subcategory of a presentable stable module category, let $j$ be its inclusion into its stable envelope, with or without retracts.
The base may be any presentably symmetric monoidal stable category whose tensor product preserves colimits separately.
The unit of $j_!\dashv j^*$ is invertible on representables by fullness, hence on all modules since both functors preserve weighted colimits.
The image of $j_!$ is stable, closed under retracts and contains the original representables.
By exactness of enriched Hom in each variable, the image contains every representable of the envelope.
Every enriched module is a weighted colimit of representables, so $j_!$ is essentially surjective on both left and right modules.

For the locally rigid source $\mathcal E$, linearity of $j_!$ and the identity $\mathcal T(a\otimes_{\mathcal E}T,-)\simeq a^\vee\otimes_{\mathcal E}\mathcal T(T,-)$ for $a\in\mathcal E^\omega$ imply the same conclusion for $\thick_{\mathcal E^\omega}(\mathcal M)$.
Restriction identifies the evaluation modules in both cases and hence induces the equivalence under $\mathcal C$.
\end{proof}

Over a ring these are the usual double-centralizer invariance properties \cite[Propositions~3.2(2) and~3.4]{EfimovCompletion}.

Fix full small enriched subcategories $\mathcal M\subset\LMod_{\mathcal C}(\mathcal E)$ and $\mathcal N\subset\LMod_{\mathcal D}(\mathcal E)$.

\begin{proposition}\label{prop:completion-functoriality}
Let $F:\mathcal C\to\mathcal D$ be enriched and suppose $F^*\mathcal N\subset\overline{\mathcal M}=\thick_{\mathcal E^\omega}(\mathcal M)$.
Then $F$ induces an enriched functor $\widehat F:\widehat{\mathcal C}_{\mathcal M}\to\widehat{\mathcal D}_{\mathcal N}$ compatible with the completion functors.
It extends to the idempotent-complete stable envelopes and to the $\mathcal E^\omega$-linear envelopes, compatibly with composition in $F$.
\end{proposition}

\begin{proof}
Put $j:\mathcal M\hookrightarrow\overline{\mathcal M}$ and $f=F^*|_{\mathcal N}:\mathcal N\to\overline{\mathcal M}$.
The counit of $j_!\dashv j^*$ is invertible by Proposition~\ref{prop:testing-envelope}, so
\begin{equation}\label{eq:functorial-evaluation-modules}
 (f^*j_!e_X)(N)\simeq(F^*N)(X)=N(FX)=e_{FX}(N).
\end{equation}
Thus $f^*j_!$ induces $\widehat F$ on evaluation modules and preserves both envelopes by exactness and $\mathcal E$-linearity.
On the $\mathcal E^\omega$-linear stable envelopes of $\mathcal M$ and $\mathcal N$, this is restriction along $F^*:\thick_{\mathcal E^\omega}(\mathcal N)\to\overline{\mathcal M}$, so the induced functors compose coherently.
\end{proof}

For full inclusions $\mathcal T\subset\mathcal C$ and $\mathcal U\subset\mathcal D$, an enriched adjunction $G\dashv F$ with $G(\mathcal U)\subset\mathcal T$ supplies the hypothesis through $\mathcal D(u,FX)\simeq\mathcal C(Gu,X)$.
The condition $F(\mathcal T)\subset\mathcal U$ alone is insufficient.
For a field $k$, the sum functor $\Perf(k)\oplus\Perf(k)\to\Perf(k)$ takes $\Perf(k)\oplus0$ into $\Perf(k)$ and sends $(0,k)$ to $k$.
The image of $(0,k)$ under completion along $\Perf(k)\oplus0$ is zero, whereas the image of $k$ under completion along the target $\Perf(k)$ is nonzero.

\begin{proposition}\label{prop:completion-morita}
Suppose an $\mathcal E$-linear colimit-preserving equivalence $H:\LMod_{\mathcal C}(\mathcal E)\simeq\LMod_{\mathcal D}(\mathcal E)$ restricts to an equivalence $h:\mathcal M\simeq\mathcal N$.
Then $h_!$ identifies the $\mathcal E^\omega$-linear stable envelopes of $\widehat{\mathcal C}_{\mathcal M}$ and $\widehat{\mathcal D}_{\mathcal N}$.
These are the idempotent-complete stable envelopes when $\mathcal E=\mathcal D(k)$ or both original categories are $\mathcal E^\omega$-linear.
\end{proposition}

\begin{proof}
For $K_X=H(\mathcal C(X,-))$, Yoneda identifies $(h_!e_X)(N)$ with $\sHom_{\mathcal D}(K_X,N)$.
The compact module $K_X$ belongs to the $\mathcal E^\omega$-linear stable envelope of the representables by \eqref{eq:compact-enriched-modules}.
Mapping into $N$ turns finite colimits into finite limits and compact scalar tensors into tensors with their duals, so $h_!e_X$ lies in the $\mathcal E^\omega$-linear stable envelope of the evaluation modules for $\mathcal D$.
Applying the inverse $h^*$ proves the reverse inclusion.
When $\mathcal E=\mathcal D(k)$, the identity $\Perf(k)=\thick\{k\}$ identifies the two envelopes.
When both source categories are $\mathcal E^\omega$-linear, the equivalences $a\otimes_{\mathcal E}e_X\simeq e_{a\otimes_{\mathcal E}X}$ show that their stable envelopes are already closed under compact scalar tensors.
\end{proof}

If $H$ identifies the full subcategories of representable left modules, then $h_!$ restricts to an enriched equivalence $\widehat{\mathcal C}_{\mathcal M}\simeq\widehat{\mathcal D}_{\mathcal N}$ before taking their stable envelopes.

For $\mathcal C$ small $\mathcal E^\omega$-linear, replace $\mathcal M$ by its $\mathcal E^\omega$-linear stable envelope as in Proposition~\ref{prop:testing-envelope}.
The inclusion $\mathcal M\hookrightarrow\LMod_{\mathcal C}(\mathcal E)$ extends by colimits to $g:\Ind\mathcal M\to(\Ind\mathcal C)^\vee$ via the equivalences of \cite[Propositions~4.5, 4.10 and Corollary~4.11]{HSS},
\begin{equation}\label{eq:presentable-dual}
 (\Ind\mathcal C)^\vee
 \simeq\Fun^L_{\mathcal E}(\Ind\mathcal C,\mathcal E)
 \simeq\LMod_{\mathcal C}(\mathcal E).
\end{equation}
Its transpose $g^\vee:\Ind\mathcal C\to(\Ind\mathcal M)^\vee$ is the $\mathcal E$-linear colimit extension of $e$ because $(g^\vee X)(M)=(gM)(X)=e_X(M)$ for $X\in\mathcal C$ and $M\in\mathcal M$.
The thick closure of $g^\vee((\Ind\mathcal C)^\omega)$ is therefore the stable envelope of the completion, already closed under compact scalar tensors since $a\otimes_{\mathcal E}e_X\simeq e_{a\otimes_{\mathcal E}X}$.

Write $e$ also for its colimit extension, with right adjoint $e^R$.
By the enriched adjunction,
\[
 \widehat{\mathcal C}_{\mathcal M}(X,Y)
 \simeq\sHom_{\Ind\mathcal C}(X,e^ReY),
\]
with composition transported from the category of evaluation modules.
Let $\eta$ and $\epsilon$ be the unit and counit of $e\dashv e^R$, and let $\mu=e^R\epsilon e$ be the multiplication of the monad $e^Re$.
The transported composition is called \emph{Kleisli composition}:
\begin{equation}\label{eq:kleisli-composition}
 \begin{aligned}
 &\sHom_{\Ind\mathcal C}(Y,e^ReZ)\otimes_{\mathcal E}
   \sHom_{\Ind\mathcal C}(X,e^ReY)
   \longrightarrow\sHom_{\Ind\mathcal C}(X,e^ReZ),\\
 &(g,f)\longmapsto\mu_Z\circ(e^Re)(g)\circ f.
 \end{aligned}
\end{equation}
Here the map on enriched Hom objects applies $e^Re$ to the first factor, composes in $\Ind\mathcal C$, and postcomposes with $\mu_Z$; the identity at $X$ is $\eta_X$.
The enriched adjunction transports the unit and associativity coherences from the module category.
This monad need not be idempotent; for a full linear inclusion $u:\mathcal T\hookrightarrow\mathcal C$, it is $u_*u^*$, since $e$ is restriction and $e^R$ is enriched right Kan extension.

\subsection{Change of coefficients}\label{subsec:change-coefficients}

Let $\mathcal E$ satisfy the hypotheses of \S~\ref{subsec:source-categories}, and let $\mathcal V$ be presentable, stable, and symmetric monoidal, with tensor product preserving colimits separately.
Fix an adjunction
\begin{equation}\label{eq:coefficient-adjunction}
 \delta:\mathcal E\rightleftarrows\mathcal V:\Gamma
\end{equation}
with $\delta$ fully faithful and strong symmetric monoidal.
We write $\otimes_{\mathcal E}$ and $\otimes_{\mathcal V}$ for coefficient tensor products and their module actions.

For a small $\mathcal E$-enriched $\mathcal B$, let the $\mathcal V$-enriched category $\delta\mathcal B$ have the same objects and mapping objects $\delta(\mathcal B(b,c))$, with composition induced by monoidality.
Applying $\delta$ objectwise to modules and their actions defines $\delta:\RMod_{\mathcal B}(\mathcal E)\to\RMod_{\delta\mathcal B}(\mathcal V)$.
The lax monoidal structure of $\Gamma$ and the unit $\mathcal B\to\Gamma\delta\mathcal B$ similarly define an objectwise functor in the reverse direction.
These are the module forms of change of enrichment \cite[Lemma~4.3.19 and Proposition~5.7.17]{GepnerHaugseng}; the same definitions apply to left modules.
The image under $\delta$ of module evaluation $\sHom_{\mathcal B}(M,N)\otimes_{\mathcal E}M\to N$ has tensor--Hom adjoint
\begin{equation}\label{eq:module-coefficient-comparison}
 \kappa_{M,N}:\delta\sHom_{\mathcal B}(M,N)
 \longrightarrow\sHom_{\delta\mathcal B}(\delta M,\delta N).
\end{equation}

Let $\mathcal C$ be small $\mathcal E$-enriched and let $\mathcal M\subset\LMod_{\mathcal C}(\mathcal E)$ be a full small $\mathcal E$-enriched subcategory.
Define $\widehat{\mathcal C}_{\mathcal M}^{\mathcal V}$ to be the full $\mathcal V$-enriched category on the family of modules $(\delta e_X)_{X\in\mathcal C}$ in $\LMod_{\delta\mathcal M}(\mathcal V)$, indexed by the objects of $\mathcal C$, so
\begin{equation}\label{eq:solid-completion}
 \widehat{\mathcal C}_{\mathcal M}^{\mathcal V}(X,Y)
 =\sHom_{\delta\mathcal M}(\delta e_X,\delta e_Y).
\end{equation}
The functor $\delta\mathcal C\to\widehat{\mathcal C}_{\mathcal M}^{\mathcal V}$ is the identity on objects and has maps on mapping objects
\[
 \delta\mathcal C(X,Y)\longrightarrow\delta\sHom_{\mathcal M}(e_X,e_Y)
 \longrightarrow\sHom_{\delta\mathcal M}(\delta e_X,\delta e_Y),
\]
given by evaluation and $\kappa_{e_X,e_Y}$ in~\eqref{eq:module-coefficient-comparison}.
Taking the coefficient adjoint of the second arrow defines
\begin{equation}\label{eq:completion-coefficient-comparison}
 \widehat{\mathcal C}_{\mathcal M}
 \longrightarrow\Gamma\widehat{\mathcal C}_{\mathcal M}^{\mathcal V},
\end{equation}
where the lax monoidal structure of $\Gamma$ defines the $\mathcal E$-enrichment on the target.

After change of coefficients, we form the envelope generated by the $\delta e_X$ under finite cofibers, shifts, retracts, and tensoring with $\delta a$ for $a\in\mathcal E^\omega$.
These scalars are dualizable and need not be compact in $\mathcal V$.

\begin{proposition}\label{thm:solid-enrichment}\label{lem:coefficient-adjunction}
The objectwise module functors form an adjunction $\delta\dashv\Gamma$ with fully faithful left adjoint.
The coefficient adjoint of $\kappa_{M,N}$ in~\eqref{eq:module-coefficient-comparison} is an equivalence
\begin{equation}\label{eq:module-point-comparison}
 \sHom_{\mathcal B}(M,N)
 \xrightarrow{\ \sim\ }\Gamma\sHom_{\delta\mathcal B}(\delta M,\delta N).
\end{equation}
For perfect $M$, the map $\kappa_{M,N}$ in~\eqref{eq:module-coefficient-comparison} is itself an equivalence in $\mathcal V$.
In particular, \eqref{eq:completion-coefficient-comparison} is an $\mathcal E$-enriched equivalence under $\mathcal C$.

With the hypotheses and notation of Propositions~\ref{prop:testing-envelope}--\ref{prop:completion-morita}, restriction $(\delta j)^*$ identifies the completion for either $\thick(\mathcal M)$ or $\thick_{\mathcal E^\omega}(\mathcal M)$ with that for $\mathcal M$.
For the functoriality assertion, put
\[
 \overline{\mathcal M}=\thick_{\mathcal E^\omega}(\mathcal M),
 \qquad j:\mathcal M\hookrightarrow\overline{\mathcal M},
 \qquad f=F^*|_{\mathcal N}:\mathcal N\to\overline{\mathcal M}.
\]
The enriched functors $(\delta f)^*(\delta j)_!$ and $(\delta h)_!$ induce respectively the functor on completions and the Morita equivalence of the envelopes defined above, compatibly with composition.
These envelopes are idempotent-complete stable envelopes for $\mathcal E=\mathcal D(k)$, and in the Morita comparison also when both original categories are $\mathcal E^\omega$-linear.
\end{proposition}

\begin{proof}
The objectwise adjunction on free modules extends over their bar resolutions, and its unit is invertible objectwise.
For $a\in\mathcal E$, we apply $\Map_{\mathcal E}(a,-)$ and use $\delta a\otimes_{\mathcal V}\delta M\simeq\delta(a\otimes_{\mathcal E}M)$ together with tensor--Hom adjunction to identify
\[
 \Gamma\sHom_{\delta\mathcal B}(\delta M,P)\simeq\sHom_{\mathcal B}(M,\Gamma P).
\]
At $P=\delta N$ this is the adjoint of $\kappa_{M,N}$ in~\eqref{eq:module-coefficient-comparison}, compatible with units and composition by evaluation.
For $M=a\otimes_{\mathcal E}\mathcal B(-,b)$ and $a$ compact, $\kappa_{M,N}$ is the equivalence $\delta(a^\vee\otimes_{\mathcal E}N(b))\simeq(\delta a)^\vee\otimes_{\mathcal V}\delta N(b)$.
Exactness and retracts extend this calculation to perfect sources by \eqref{eq:compact-enriched-modules}.
Applying~\eqref{eq:module-point-comparison} to $e_X,e_Y$ proves the assertion for the completion.

For the envelope assertion, let $j:\mathcal M\hookrightarrow\mathcal T$ be either inclusion.
In the proof of Proposition~\ref{prop:testing-envelope}, each covariant representable $\mathcal T(t,-)$ is obtained from the modules $\mathcal T(m,-)$, $m\in\mathcal M$, using shifts, finite limits, and retracts, together with tensoring by $a^\vee$, $a\in\mathcal E^\omega$, when $\mathcal T=\thick_{\mathcal E^\omega}(\mathcal M)$.
Applying $\delta$ preserves these operations, including $\delta(a^\vee)\simeq(\delta a)^\vee$, and sends $\mathcal T(t,-)$ to $(\delta\mathcal T)(t,-)$.
Since $(\delta j)_!$ sends representables to representables and its essential image is closed under these operations, that image contains every representable of $\delta\mathcal T$.
The functor $\delta$ preserves tensor products and colimits, and hence the coends defining enriched Kan extension.
Applying it to~\eqref{eq:functorial-evaluation-modules} identifies the images of the evaluation modules, proving the assertion about the induced functor on completions.
The calculation in the proof of Proposition~\ref{prop:completion-morita} therefore shows that $(\delta h)_!$ and $(\delta h)^*$ identify the specified envelopes.
These Kan-extension identifications preserve the coherent composition laws.
Objectwise full faithfulness identifies the underlying $\mathcal E$-enriched category of this envelope with the $\mathcal E^\omega$-linear stable envelope before change of coefficients.
\end{proof}

\begin{remark}\label{rem:fixed-testing-representation}
The category $\delta\mathcal M$ obtained by change of coefficients has $\mathcal V$-valued mapping objects $\delta\sHom_{\mathcal C}(M,N)$.
The full $\mathcal V$-enriched category on the modules $\delta M$ instead has mapping objects $\sHom_{\delta\mathcal C}(\delta M,\delta N)$.
The maps $\kappa_{M,N}$ in~\eqref{eq:module-coefficient-comparison} identify these mapping objects when every $M\in\mathcal M$ is perfect as a $\mathcal C$-module; this applies in particular when the objects of $\mathcal M$ are representable.
Proposition~\ref{prop:centralizer-counterexample} shows that these two constructions can yield different completions.
\end{remark}

\subsection{Solid coefficients over a discrete ring}
\label{subsec:solid-coefficients}
\label{rem:relative-solid-coefficients}

Fix an arbitrary discrete commutative ring $k$ and set $\mathcal E=\mathcal D(k)$.
The Clausen--Scholze relative analytic ring $(k,\Z)_\square$ of \cite[Example~7.3(iii) and Proposition~7.8(ii)]{Scholze} has derived category
\[
 \mathcal V=\mathcal D((k,\Z)_\square)
 \simeq\Mod_{\delta k}(\DS)
\]
by \cite[Proposition~7.5(ii)]{Scholze}: its objects are condensed $k$-complexes with solid underlying complexes.
Write $\otimes_{\delta k}^{\square}$ and $\sHom_{\delta k}$ for its tensor product and internal Hom.
The discrete embedding and point evaluation define
\begin{equation}\label{eq:discrete-adjunction}
 \delta:\mathcal D(k)\rightleftarrows\mathcal D((k,\Z)_\square):\Gamma.
\end{equation}
The integral adjunction lifts to module objects: $\delta$ is fully faithful, exact, colimit preserving, and strong symmetric monoidal, since it preserves the bar constructions for relative tensor products.
Point evaluation $\Gamma$ is exact and preserves limits and colimits by compact projectivity of the solid unit \cite[Corollary~6.1]{Scholze}.

For $S=\varprojlim_iS_i$ profinite with finite $S_i$, define
\[
 \Z[S]^\square=\varprojlim_i\Z[S_i],\qquad
 k[S]^\square=\delta k\otimes_{\Z_\square}\Z[S]^\square.
\]
The first limit is formed in condensed abelian groups.
By free-module adjunction and \cite[Corollaries~5.5 and~6.1]{Scholze}, the $k[S]^\square$ are compact projective generators of the heart and compact generators of $\mathcal D((k,\Z)_\square)$.
After extending scalars in \cite[Theorem~6.2(ii) and Proposition~6.3]{Scholze}, we obtain the free-module tensor equivalence
\begin{equation}\label{eq:relative-profinite-tensor}
 k[S]^\square\otimes_{\delta k}^{\square}k[T]^\square
 \simeq k[S\times T]^\square.
\end{equation}
For $k=\Z$ or $\F_p$, dualizability of $\delta k$ over $\Z_\square$ also implies $k[S]^\square\simeq\varprojlim_i\delta(k[S_i])$.
For general $k$ the scalar-extension definition is used.

For $H\in\mathcal D((k,\Z)_\square)$, its sections on $S$ form the $\mathcal D(k)$-valued complex
\[
 H(S)=\Hom_{\mathcal D((k,\Z)_\square)}(k[S]^\square,H).
\]
For a category enriched over $\mathcal D((k,\Z)_\square)$, we obtain the $\mathcal D(k)$-valued mapping complexes between its modules by applying $\Gamma$ to their internal Homs.
Put $C(S,k)=\colim_i k^{S_i}=C(S,\Z)\otimes_{\Z}k$, the flat $k$-algebra of locally constant functions.
For a topological abelian group or ring $A$, write $\underline A$ for the associated condensed object; its topology will be specified when used.
We write $\widehat{\mathcal C}_{\mathcal M}^{\square}$ for the completion over these coefficients.
Proposition~\ref{thm:solid-enrichment} applies for every $k$.
Forgetting coefficient module structures recovers the integral construction, since the units $k\to\mathcal B(b,b)$ recover the $k$-action on modules over the underlying $\DZ$-enriched category, and the solid analogue follows by free-module adjunction.

The filtered-colimit construction of \cite[Proposition~13.16]{ScholzeAnalytic} also defines the absolute analytic ring $(k,k)_\square$.
Analytic scalar extension is a symmetric monoidal localization with fully faithful right adjoint \cite[Propositions~12.4 and~12.6, Remark~12.5]{ScholzeAnalytic}:
\[
 \begin{gathered}
 j^*:\mathcal D((k,\Z)_\square)\rightleftarrows\mathcal D((k,k)_\square):j_*,\\
 j^*H=H\otimes_{(k,\Z)_\square}(k,k)_\square.
 \end{gathered}
\]
Discrete complexes lie in $\mathcal D((k,k)_\square)$ \cite[Lemma~5.7]{Andreychev}, so $\delta_{\mathrm{abs}}=j^*\delta$ is fully faithful and strong symmetric monoidal, with right adjoint $\Gamma j_*$.

For modules $M,N$ over a small $\mathcal D(k)$-enriched $\mathcal B$, the lax monoidal inclusion and evaluation induce
\[
 j_*\sHom_{\delta_{\mathrm{abs}}\mathcal B}
 (\delta_{\mathrm{abs}}M,\delta_{\mathrm{abs}}N)
 \longrightarrow
 \sHom_{\delta\mathcal B}(\delta M,\delta N).
\]
The end formula for the internal Hom between modules and the adjunction
\[
 \sHom_{\mathcal D((k,\Z)_\square)}(X,j_*Y)
 \simeq j_*\sHom_{\mathcal D((k,k)_\square)}(j^*X,Y)
\]
show that this map is an equivalence.
Thus every mapping object of $\widehat{\mathcal C}_{\mathcal M}^{\square}$, and of its stable idempotent-complete envelope, lies in the essential image of
\[
 j_*:\mathcal D((k,k)_\square)\longrightarrow\mathcal D((k,\Z)_\square).
\]
For every profinite set $S$, we deduce from adjunction and Proposition~\ref{thm:family-formula} that
\[
 \begin{aligned}
 &\Hom_{\mathcal D((k,k)_\square)}
 \bigl(j^*k[S]^\square,
 \sHom_{\delta_{\mathrm{abs}}\mathcal B}
 (\delta_{\mathrm{abs}}M,\delta_{\mathrm{abs}}N)\bigr)\\
 &\qquad\simeq\Hom_{\mathcal B}(M,C(S,k)\otimes_kN).
 \end{aligned}
\]
Moreover, $M$ is perfect if and only if the canonical map
\[
 \delta_{\mathrm{abs}}\Hom_{\mathcal B}(M,N)
 \longrightarrow
 \sHom_{\delta_{\mathrm{abs}}\mathcal B}
 (\delta_{\mathrm{abs}}M,\delta_{\mathrm{abs}}N)
\]
is an equivalence for every $N$.
Indeed, the fully faithful functor $j_*$ identifies this map with $\kappa_{M,N}$ in~\eqref{eq:module-coefficient-comparison}, so Proposition~\ref{prop:perfection-detection} applies.
The two coefficient categories agree for $k=\Z$ or $\F_p$ \cite[Proposition~3.22]{Andreychev}.
For $k=\Z[T]$ the localization has a nonzero kernel \cite[Observation~8.10]{Scholze}.
Throughout the solid construction, we form the module categories and impose the Cauchy condition over $\mathcal D((k,\Z)_\square)$.

\subsection{Completion along a full subcategory}\label{sec:support}\label{subsec:full-subcategory}

Fix a discrete commutative ring $k$.
The source categories are enriched over $\mathcal D(k)$.
Kan extensions, internal Homs between modules, and coends after applying $\delta$ are taken over $\mathcal D((k,\Z)_\square)$, as in~\eqref{eq:discrete-adjunction}.
The corresponding constructions over $\mathcal D(k)$ are obtained by taking $\delta=\id$.
Let $\mathcal C$ be small $\Perf(k)$-linear and $u:\mathcal T\hookrightarrow\mathcal C$ a full thick inclusion; $\mathcal T$ is automatically $\Perf(k)$-linear.
Put $h_X=\delta\mathcal C(-,X)$ and denote restriction along $\delta u$ by $u^*:\RMod_{\delta\mathcal C}(\mathcal D((k,\Z)_\square))\to\RMod_{\delta\mathcal T}(\mathcal D((k,\Z)_\square))$, so $(u^*M)(t)=M(ut)$.
Left and right Kan extension over $\mathcal D((k,\Z)_\square)$ are given by
\begin{equation}\label{eq:support-kan-extensions}
 \begin{aligned}
 (u_!N)(X)&=\int^{t\in\delta\mathcal T}\delta\mathcal C(X,ut)\otimes_{\delta k}^{\square}N(t),\\
 (u_*N)(X)&=\sHom_{\delta\mathcal T}(u^*h_X,N),
 \end{aligned}
\end{equation}
with actions induced by composition in $\mathcal C$.
The enriched Kan-extension theorem establishes the adjunctions $u_!\dashv u^*\dashv u_*$; see \cite[Theorems~4.50 and~4.99]{Kelly} and the $\infty$-categorical adjoint theorem \cite[Theorem~3.73]{HeineWeighted}.
Fullness and Yoneda make the unit $a:\operatorname{id}\to u^*u_!$ and counit $b:u^*u_*\to\operatorname{id}$ invertible.

Set $\Gamma_{\mathcal T}=u_!u^*$, $\Lambda_{\mathcal T}=u_*u^*$, and $\mathcal N=\ker(u^*)$.
Write $\epsilon:\Gamma_{\mathcal T}\to\operatorname{id}$ for the counit of $u_!\dashv u^*$ and $\eta:\operatorname{id}\to\Lambda_{\mathcal T}$ for the unit of $u^*\dashv u_*$.
The monad multiplication of $\Lambda_{\mathcal T}$ is $\mu=u_*bu^*$.
Call $M$ $\mathcal T$-complete if $\eta_M$ is invertible, and put $\mathcal N^\perp=\{M\mid\sHom_{\delta\mathcal C}(K,M)=0\text{ for all }K\in\mathcal N\}$.
This also equals the right orthogonal for underlying mapping complexes, since $\mathcal N$ is closed under the $\mathcal D((k,\Z)_\square)$-action.

The category $\widehat{\mathcal C}_{\mathcal T}^{\square}$ is the full enriched category on the modules $u^*h_X=\delta r_{\mathcal T}X$, indexed by $X\in\mathcal C$.
Restricted enriched Yoneda factors as the completion functor $\delta\mathcal C\to\widehat{\mathcal C}_{\mathcal T}^{\square}$, followed by the fully faithful inclusion of $\widehat{\mathcal C}_{\mathcal T}^{\square}$ into the category of solid $\delta\mathcal T$-modules:
\begin{equation}\label{eq:support-completion-functor}
 \begin{aligned}
 \delta\mathcal C&\longrightarrow\widehat{\mathcal C}_{\mathcal T}^{\square}
   \hookrightarrow\RMod_{\delta\mathcal T}(\mathcal D((k,\Z)_\square)),\\
 X&\longmapsto X\longmapsto u^*h_X=\delta r_{\mathcal T}X.
 \end{aligned}
\end{equation}
Left and right Kan extension realize the second functor of~\eqref{eq:support-completion-functor} fully faithfully in the category of solid $\delta\mathcal C$-modules, carrying $X$ to $\Gamma_{\mathcal T}h_X$ and $\Lambda_{\mathcal T}h_X$, respectively.
For the second realization, precomposition with the completion unit is an equivalence:
\begin{equation}\label{eq:completed-unit-precomposition}
 \sHom_{\delta\mathcal C}(\Lambda_{\mathcal T}h_X,\Lambda_{\mathcal T}h_Y)
 \xrightarrow[\ \sim\ ]{\eta_{h_X}^{*}}
 \sHom_{\delta\mathcal C}(h_X,\Lambda_{\mathcal T}h_Y).
\end{equation}
Its inverse applies the enriched functor $\Lambda_{\mathcal T}$ to Hom objects and postcomposes with $\mu_{h_Y}$, sending $f:h_X\to\Lambda_{\mathcal T}h_Y$ to $\mu_{h_Y}\Lambda_{\mathcal T}(f)$.
Under~\eqref{eq:completed-unit-precomposition}, composition and units are those of~\eqref{eq:kleisli-composition} for the monad $\Lambda_{\mathcal T}=u_*u^*$.

\begin{proposition}\label{thm:support-completion}
The functors $\Gamma_{\mathcal T}$ and $\Lambda_{\mathcal T}$ are idempotent colocalization and localization, respectively.
From the adjunctions $u_!\dashv u^*\dashv u_*$, we deduce
\begin{equation}\label{eq:torsion-completion-identities}
 \Gamma_{\mathcal T}\dashv\Lambda_{\mathcal T},\qquad
 \Gamma_{\mathcal T}\Lambda_{\mathcal T}\xrightarrow[\ \sim\ ]{u_!bu^*}\Gamma_{\mathcal T},\qquad
 \Lambda_{\mathcal T}\xrightarrow[\ \sim\ ]{u_*au^*}\Lambda_{\mathcal T}\Gamma_{\mathcal T}.
\end{equation}
The essential images of $\Gamma_{\mathcal T}$ and $\Lambda_{\mathcal T}$ are enriched equivalent to $\RMod_{\delta\mathcal T}(\mathcal D((k,\Z)_\square))$, and $\operatorname{im}(\Lambda_{\mathcal T})=\mathcal N^\perp$.
The maps induced by $u_!$, $u_*$, and precomposition~\eqref{eq:completed-unit-precomposition} are equivalences
\begin{align}\label{eq:completed-kan-morphisms}
 \widehat{\mathcal C}_{\mathcal T}^{\square}(X,Y)
 &\simeq\sHom_{\delta\mathcal C}(h_X,\Lambda_{\mathcal T}h_Y)\notag\\
 &\simeq\sHom_{\delta\mathcal C}(\Lambda_{\mathcal T}h_X,\Lambda_{\mathcal T}h_Y)\notag\\
 &\simeq\sHom_{\delta\mathcal C}(\Gamma_{\mathcal T}h_X,\Gamma_{\mathcal T}h_Y),
\end{align}
compatible with units and composition, where $\sHom_{\delta\mathcal C}(h_X,\Lambda_{\mathcal T}h_Y)$ carries the Kleisli composition of~\eqref{eq:kleisli-composition}.
\end{proposition}

\begin{proof}
Full faithfulness of $u_!$ and $u_*$ implies idempotence and the two absorption equivalences in~\eqref{eq:torsion-completion-identities}, while the composite adjunction is $u_!u^*\dashv u_*u^*$.
Adjunction puts $\operatorname{im}(u_*)$ in $\mathcal N^\perp$; for $M\in\mathcal N^\perp$, the fiber of $\eta_M$ lies in both $\mathcal N$ and $\mathcal N^\perp$, so vanishes.
Applying full faithfulness and $u^*\dashv u_*$ to the restricted representables establishes~\eqref{eq:completed-kan-morphisms}, with composition as computed in~\eqref{eq:completed-unit-precomposition}.
\end{proof}

The formal local-duality identities in~\eqref{eq:torsion-completion-identities} also occur in \cite[Theorem~2.21]{BHV}; here they follow directly from the Kan-extension adjunctions.
The image of $\Gamma_{\mathcal T}$ is generated by the modules $h_t$, $t\in\mathcal T$, under colimits and the $\mathcal D((k,\Z)_\square)$-action; over $\mathcal D(k)$ it is the colimit closure of $\mathcal T$ in $\Ind\mathcal C$.
The idempotent-complete stable envelope of $\widehat{\mathcal C}_{\mathcal T}^{\square}$ is $\thick\{\Lambda_{\mathcal T}h_X\}$ inside $\operatorname{im}(\Lambda_{\mathcal T})$.

\subsection{Modules over the completed category}

For a full small enriched subcategory $\mathcal M\subset\LMod_{\mathcal C}(\mathcal D(k))$, the category of solid right modules over the completion of $\mathcal C$ along $\mathcal M$ is
\begin{equation}\label{eq:free-solid-envelope}
 \RMod_{\widehat{\mathcal C}_{\mathcal M}^{\square}}(\mathcal D((k,\Z)_\square)).
\end{equation}
Right modules over a small category $\mathcal H$ enriched over $\mathcal V=\mathcal D((k,\Z)_\square)$ form a presentable stable category with objectwise $\mathcal V$-action.
For every presentable stable $\mathcal V$-module category $\mathcal D$, restriction along the Yoneda embedding $yX=\mathcal H(-,X)$ is an equivalence:
\begin{equation}\label{eq:enriched-presheaf-universal}
 \Fun^L_{\mathcal D((k,\Z)_\square)}(\RMod_{\mathcal H}(\mathcal D((k,\Z)_\square)),\mathcal D)
 \simeq\Fun_{\mathcal D((k,\Z)_\square)}(\mathcal H,\mathcal D).
\end{equation}
The inverse is enriched left Kan extension, $F\mapsto(N\mapsto\int^{X\in\mathcal H}N(X)\otimes_{\mathcal V}F(X))$, by the enriched presheaf theorem \cite[Theorem~4.51]{Kelly} in its higher-categorical form \cite[Proposition~3.23(2) and Corollary~3.68]{HeineWeighted}.
Restriction from the stable idempotent-complete envelope of $\mathcal H$ is an equivalence on modules by Proposition~\ref{prop:testing-envelope}.

For the full inclusion $u$ above, transport the presentable $\mathcal D((k,\Z)_\square)$-module structure along $u_*:\RMod_{\delta\mathcal T}(\mathcal D((k,\Z)_\square))\simeq\operatorname{im}(u_*)$.
Its action on complete modules is
\begin{equation}\label{eq:complete-solid-action}
 (P,N)\longmapsto u_*(P\otimes_{\delta k}^{\square}u^*N),
\end{equation}
and colimits are obtained by applying $\Lambda_{\mathcal T}$ to colimits in the category of solid $\delta\mathcal C$-modules.
Let $j:\delta\mathcal T\to\widehat{\mathcal C}_{\mathcal T}^{\square}$ be the restriction of the completion functor~\eqref{eq:support-completion-functor}; it is fully faithful by Yoneda.
Define the composite
\begin{equation}\label{eq:free-to-complete-functor}
 \RMod_{\widehat{\mathcal C}_{\mathcal T}^{\square}}(\mathcal D((k,\Z)_\square))
 \xrightarrow{j^*}\RMod_{\delta\mathcal T}(\mathcal D((k,\Z)_\square))
 \xrightarrow{u_*}\operatorname{im}(u_*).
\end{equation}
\begin{proposition}\label{prop:free-to-complete}
The functor $u_*j^*$ in~\eqref{eq:free-to-complete-functor} is a $\mathcal D((k,\Z)_\square)$-linear localization for the action~\eqref{eq:complete-solid-action} on $\operatorname{im}(u_*)$.
It carries the representable $\widehat{\mathcal C}_{\mathcal T}^{\square}(-,X)$ to $\Lambda_{\mathcal T}h_X$.
\end{proposition}

\begin{proof}
Restriction $j^*$ preserves colimits and the action, and has fully faithful right Kan extension because $j$ is full.
By Yoneda, $j^*(\widehat{\mathcal C}_{\mathcal T}^{\square}(-,X))\simeq\delta r_{\mathcal T}X$.
Composing with the equivalence $u_*$ for the transported action establishes the stated image and localization.
\end{proof}

For $\mathcal C=\mathcal T=\Perf(k)$ the free solid module category is $\mathcal D((k,\Z)_\square)$, whereas the Ind-completion of the underlying stable $\infty$-category of the envelope is equivalent to $\mathcal D(k)$.

\begin{remark}\label{prop:solid-swindle}
The compact profinite generators and \eqref{eq:relative-profinite-tensor} make the compact solid coefficients symmetric monoidal.
The compact scalar $P_k=\delta k\otimes_{\Z_\square}\prod_{n\geq0}\Z$ splits as $\delta k\oplus P_k$ \cite[Corollary~6.1]{Scholze}.
If $\Perf(k)$ admitted an exact unital action of all compact solid scalars, tensoring with $P_k$ would define an endofunctor $F\simeq\id\oplus F$, forcing every additive invariant of $\Perf(k)$ to vanish by the Eilenberg swindle \cite[Lemma~9.4]{BGT}.
For $k\ne0$, scalar extension to a residue field detects $[k]$ in $K_0(\Perf(k))$, excluding such an action even though $\delta\Perf(k)$ has solid mapping objects.
\end{remark}

\section{Cauchy completion and punctures}\label{sec:cauchy-punctures}\label{sec:solid}

We now study families of morphisms in the completed category and form its puncture.
Throughout this section, $k$ is a discrete commutative ring.
We form solid module categories, enriched coends, enriched quotients, and Cauchy completions over $\mathcal D((k,\Z)_\square)$.
The categories and modules to which we apply discrete realization are enriched over $\mathcal D(k)$.
Point evaluation $\Gamma$ takes solid mapping objects to complexes in $\mathcal D(k)$.

\subsection{Families of morphisms parametrized by profinite sets}

The projections $k[S]^\square\to\delta(k[S_i])$ and finite evaluation pairings define
\begin{equation}\label{eq:profinite-evaluation}
 k[S]^\square\otimes_{\delta k}^{\square}\delta C(S,k)
 \longrightarrow\delta k.
\end{equation}
Tensoring with $\delta N$ and taking its adjoint defines, for $N\in\mathcal D(k)$,
\begin{equation}\label{eq:solid-profinite-hom}
 \delta(C(S,k)\otimes_kN)
 \longrightarrow\sHom_{\delta k}(k[S]^\square,\delta N).
\end{equation}

\begin{lemma}\label{lem:solid-profinite-hom}
The map \eqref{eq:solid-profinite-hom} is an equivalence for every profinite $S$ and every $N\in\mathcal D(k)$.
\end{lemma}

\begin{proof}
For $N=k$, applying $\Hom_{\mathcal D((k,\Z)_\square)}(k[T]^\square,-)$ reduces the comparison to $C(T,C(S,k))\simeq C(T\times S,k)$, which follows from finite clopen partitions; compact projectivity excludes higher cohomology.
Both sides preserve colimits in $N$: on the right, composition with $\Hom_{\mathcal D((k,\Z)_\square)}(k[T]^\square,-)$ is the exact mapping functor from the compact object $k[T\times S]^\square$.
Since $k$ generates $\mathcal D(k)$, the assertion follows for unbounded $N$, compatibly with coefficient actions and hence objectwise for modules.
\end{proof}

For right modules $M,N$ over a small $\mathcal D(k)$-enriched $\mathcal B$, apply $\delta$ to universal module evaluation into $C(S,k)\otimes_kN$, tensor with $k[S]^\square$, and compose with \eqref{eq:profinite-evaluation} tensored with $\delta N$.
Taking adjoints of this composite defines
\begin{equation}\label{eq:family-formula}
 \theta_{S;M,N}:\Hom_{\mathcal B}(M,C(S,k)\otimes_kN)
 \longrightarrow\sHom_{\delta\mathcal B}(\delta M,\delta N)(S).
\end{equation}
Module evaluation over $\mathcal D(k)$ similarly defines
\begin{equation}\label{eq:discreteness-comparison}
 \chi_{S;M,N}:C(S,k)\otimes_k\Hom_{\mathcal B}(M,N)
 \longrightarrow\Hom_{\mathcal B}(M,C(S,k)\otimes_kN).
\end{equation}

For a continuous map of profinite sets $u:S'\to S$, pullback of functions is
\[
 u^*:C(S,k)\longrightarrow C(S',k),\qquad a\longmapsto a\circ u.
\]
Postcomposition with $u^*\otimes\id_N$ induces the pullback map
\[
 \Hom_{\mathcal B}(M,C(S,k)\otimes_kN)
 \longrightarrow\Hom_{\mathcal B}(M,C(S',k)\otimes_kN).
\]
On the target of $\theta_{S;M,N}$, pullback of solid sections along $u$ is induced by precomposition with $k[S']^\square\to k[S]^\square$.
For a fixed profinite set $S$, extension--restriction adjunction identifies
\[
 \begin{aligned}
 &\Hom_{C(S,k)\otimes_k\mathcal B}
 \bigl(C(S,k)\otimes_kM,C(S,k)\otimes_kN\bigr)\\
 &\hspace{2cm}\simeq\Hom_{\mathcal B}(M,C(S,k)\otimes_kN),
 \end{aligned}
\]
where the underlying $k$-complex is taken on the left.
For a third right $\mathcal B$-module $Q$, composition of the scalar-extended modules therefore defines a pairing of derived mapping complexes
\[
 \begin{aligned}
 &\Hom_{\mathcal B}(N,C(S,k)\otimes_kQ)
 \otimes_k\Hom_{\mathcal B}(M,C(S,k)\otimes_kN)\\
 &\hspace{2cm}\longrightarrow\Hom_{\mathcal B}(M,C(S,k)\otimes_kQ).
 \end{aligned}
\]
On module morphisms $f:M\to C(S,k)\otimes_kN$ and $g:N\to C(S,k)\otimes_kQ$, this pairing is the composite
\[
 \begin{aligned}
 M&\xrightarrow{f}C(S,k)\otimes_kN
 \xrightarrow{\id\otimes g}C(S,k)\otimes_kC(S,k)\otimes_kQ\\
 &\xrightarrow{\operatorname{mult}\otimes\id_Q}C(S,k)\otimes_kQ,
 \end{aligned}
\]
where $\operatorname{mult}:C(S,k)\otimes_kC(S,k)\to C(S,k)$ is multiplication of functions.
The identity at $M$ is the map $M\to C(S,k)\otimes_kM$ induced by the constant-function unit $k\to C(S,k)$.
Composition of solid sections over $S$ is induced by enriched composition and the diagonal-induced map
\[
 k[S]^\square\longrightarrow k[S\times S]^\square
 \simeq k[S]^\square\otimes_{\delta k}^{\square}k[S]^\square.
\]
The identity section is induced by $S\to *$ and the enriched identity.

\begin{proposition}\label{thm:family-formula}
For every profinite set $S$ and every pair of right $\mathcal B$-modules $M,N$, the map $\theta_{S;M,N}$ in~\eqref{eq:family-formula} is an equivalence, natural in $M,N$.
Under the canonical identification
\[
 \bigl(\delta\Hom_{\mathcal B}(M,N)\bigr)(S)
 \simeq C(S,k)\otimes_k\Hom_{\mathcal B}(M,N),
\]
the section map $\kappa_{M,N}(S)$ of~\eqref{eq:module-coefficient-comparison} satisfies
\[
 \kappa_{M,N}(S)=\theta_{S;M,N}\circ\chi_{S;M,N},
\]
where $\chi_{S;M,N}$ is the map in~\eqref{eq:discreteness-comparison}.
For every continuous map $u:S'\to S$, the maps $\theta_{S;M,N}$ and $\theta_{S';M,N}$ intertwine postcomposition with $u^*\otimes\id_N$ and pullback of solid sections along $u$.
For every $S$ and triple of right modules $M,N,Q$, the maps $\theta$ carry composition of the scalar-extended $C(S,k)\otimes_k\mathcal B$-modules to composition of solid $S$-sections induced by the diagonal $S\to S\times S$.
They also preserve the identity maps for every $M$.
\end{proposition}

\begin{proof}
Tensor--Hom adjunction identifies the target of $\theta_{S;M,N}$ with the $\mathcal D(k)$-valued mapping complex
\[
 \Hom_{\RMod_{\delta\mathcal B}(\mathcal D((k,\Z)_\square))}
 \bigl(\delta M,\sHom_{\delta k}(k[S]^\square,\delta N)\bigr).
\]
Lemma~\ref{lem:solid-profinite-hom} and Proposition~\ref{thm:solid-enrichment} identify this with the source of $\theta_{S;M,N}$, through its defining adjunctions.
Applying the same adjunctions to module evaluation establishes the identity $\kappa_{M,N}(S)=\theta_{S;M,N}\circ\chi_{S;M,N}$.
Naturality of the evaluation pairing~\eqref{eq:profinite-evaluation} under $u:S'\to S$ identifies the two pullback maps.
The same pairing intertwines multiplication on $C(S,k)$ with the map on $k[S]^\square$ induced by the diagonal $S\to S\times S$, and the constant-function unit with $S\to *$.
Applying the adjunctions defining $\theta_{S;M,N}$ therefore identifies the composition and identity maps on the two sides.
\end{proof}

If $f:M\to C(S,k)\otimes_kN$ and $g:N\to C(T,k)\otimes_kQ$, pull them back to $S\times T$ along its two projections and compose there.
Equivalently, their composite is
\[
 \begin{aligned}
 M&\xrightarrow{f}C(S,k)\otimes_kN
 \xrightarrow{\id\otimes g}C(S,k)\otimes_kC(T,k)\otimes_kQ\\
 &\simeq C(S\times T,k)\otimes_kQ.
 \end{aligned}
\]
For $T=S$, pullback along the diagonal $S\to S\times S$ recovers the composition with common parameter $S$.
Combining Proposition~\ref{thm:family-formula} with extension--restriction adjunction for $k\to C(S,k)$ identifies the $S$-sections of $\widehat{\mathcal C}_{\mathcal T}^{\square}(X,Y)$ with the derived mapping complex between $C(S,k)\otimes_k r_{\mathcal T}X$ and $C(S,k)\otimes_k r_{\mathcal T}Y$ as modules over $C(S,k)\otimes_k\mathcal T$.
Taking stable envelopes after scalar extension preserves these mapping objects.

Proposition~\ref{prop:perfection-detection} characterizes perfect right modules by requiring that $\kappa_{M,N}$ in~\eqref{eq:module-coefficient-comparison} be an equivalence for every target module $N$.

\subsection{Cauchy modules over the solid base}\label{subsec:cauchy}

Let $\mathcal B$ be a small category enriched in $\mathcal D((k,\Z)_\square)$.
We use right modules and write $h_b=\mathcal B(-,b)$.

\begin{definition}\label{def:cauchy-module}
A right $\mathcal B$-module $P$ is \emph{Cauchy} if
\[
 \sHom_{\mathcal B}(P,-):
 \RMod_{\mathcal B}\bigl(\mathcal D((k,\Z)_\square)\bigr)
 \longrightarrow \mathcal D((k,\Z)_\square)
\]
preserves colimits and is $\mathcal D((k,\Z)_\square)$-linear.
The \emph{Cauchy completion} of $\mathcal B$, denoted $\Cauchy(\mathcal B)$, is the full enriched subcategory of such modules.
\end{definition}

By \cite[Proposition~5.49]{HeineWeighted}, the Cauchy condition is equivalent to absoluteness of the enriched Yoneda presentation
\begin{equation}\label{eq:absolute-yoneda}
 P\simeq\int^{b\in\mathcal B}P(b)\otimes_{\delta k}^{\square} h_b:
\end{equation}
every $\mathcal D((k,\Z)_\square)$-enriched functor out of the module category preserves this weighted colimit.
This identifies our definition with the enriched Cauchy completion of \cite[Notation~5.46, Proposition~5.52 and Corollary~5.54]{HeineWeighted}.

Representable modules are Cauchy, and Cauchy modules are closed under finite colimits and retracts.
They are also compact in the underlying presentable category, since the unit $\delta k$ is compact in $\mathcal D((k,\Z)_\square)$.
If a linear left adjoint has a linear, colimit-preserving right adjoint, it preserves Cauchy objects by the enriched adjunction formula.
In particular, extension of scalars between solid algebras preserves Cauchy objects.

Let $\mathcal B$ now be a small $\mathcal D(k)$-enriched category, and let $\Perf(\mathcal B)$ be the thick subcategory generated by its representable right modules.
Applying $\delta$ objectwise defines a functor
\begin{equation}\label{eq:discrete-cauchy-functor}
 \delta:\delta\Perf(\mathcal B)\longrightarrow
 \RMod_{\delta\mathcal B}\bigl(\mathcal D((k,\Z)_\square)\bigr),
 \qquad P\longmapsto\delta P.
\end{equation}
Its enriched maps are $\kappa_{M,N}$ in~\eqref{eq:module-coefficient-comparison}, which respect composition and units.

\begin{lemma}\label{lem:discrete-cauchy}
The functor $\delta$ in~\eqref{eq:discrete-cauchy-functor} takes values in Cauchy modules and induces a $\mathcal D((k,\Z)_\square)$-enriched equivalence
\[
 \delta\Perf(\mathcal B)\xrightarrow{\ \sim\ }\Cauchy(\delta\mathcal B).
\]
\end{lemma}

\begin{proof}
Let $P$ be a Cauchy right $\delta\mathcal B$-module, and write $h_b=\delta\mathcal B(-,b)$.
For the compact generators $k[S]^{\square}$, Lemma~\ref{lem:solid-profinite-hom} identifies
\begin{equation}\label{eq:discrete-profinite-test}
 \sHom_{\delta k}(k[S]^{\square},\delta V)
 \simeq \delta\bigl(C(S,k)\otimes_k V\bigr)
 \qquad (V\in \mathcal D(k)).
\end{equation}
Apply the enriched functor
$\sHom_{\delta\mathcal B}(k[S]^{\square}\otimes_{\delta k}^{\square} h_b,-)$
to the absolute weighted colimit~\eqref{eq:absolute-yoneda}.
By~\eqref{eq:discrete-profinite-test} and co-Yoneda,
\begin{align*}
 \sHom_{\delta k}(k[S]^{\square},P(b))
 &\simeq
 \int^{c\in\delta\mathcal B}P(c)\otimes_{\delta k}^{\square}
   \sHom_{\delta k}\bigl(k[S]^{\square},\delta\mathcal B(b,c)\bigr)\\
 &\simeq P(b)\otimes_{\delta k}^{\square}\delta C(S,k).
\end{align*}
Evaluating at the point identifies the right-hand side with
$C(S,k)\otimes_k\Gamma P(b)$.
Indeed, $C(S,k)$ is the filtered colimit of the finite free modules of functions on finite quotients of $S$, and point evaluation preserves colimits.
These identifications carry evaluation of locally constant functions to the map on $S$-sections of the canonical counit
\[
 \delta\Gamma P(b)\longrightarrow P(b).
\]
Thus this counit is an equivalence after evaluation on every profinite set, and hence an equivalence of solid complexes.
These equivalences respect the module structure, so $P\simeq\delta\Gamma P$, where $\Gamma$ is applied objectwise.

The compactness of $P$ implies the compactness of $\Gamma P$: the objectwise discrete embedding is fully faithful and preserves colimits on module categories.
The compact objects of $\RMod_{\mathcal B}(\mathcal D(k))$ are exactly $\Perf(\mathcal B)$.
Conversely, the images under $\delta$ of representable modules are Cauchy, and this property is preserved under shifts, finite colimits, and retracts.
Enriched full faithfulness is Proposition~\ref{thm:solid-enrichment} for a perfect source.
\end{proof}

\begin{lemma}\label{lem:cauchy-discrete-evaluation}
Let $\mathcal B$ be a small $\mathcal D((k,\Z)_\square)$-enriched category, and let $b\in\mathcal B$.
Suppose that $\mathcal B(b,c)$ is discrete for every $c\in\mathcal B$.
Then $P(b)$ is discrete for every Cauchy right $\mathcal B$-module $P$.
\end{lemma}

\begin{proof}
For every profinite set $S$ and every $c\in\mathcal B$, the comparison induced by~\eqref{eq:profinite-evaluation} identifies
\[
 \mathcal B(b,c)\otimes_{\delta k}^{\square}\delta C(S,k)
 \simeq\sHom_{\delta k}(k[S]^\square,\mathcal B(b,c))
\]
by~\eqref{eq:discrete-profinite-test}.
As in the proof of Lemma~\ref{lem:discrete-cauchy}, apply $\sHom_{\mathcal B}(k[S]^\square\otimes_{\delta k}^{\square}h_b,-)$ to the absolute presentation~\eqref{eq:absolute-yoneda}.
Co-Yoneda and point evaluation show that the counit $\delta\Gamma P(b)\to P(b)$ induces an equivalence on $S$-sections.
Since this holds for every profinite $S$, the counit is an equivalence, so $P(b)$ is discrete.
\end{proof}

For a quasi-compact separated $k$-scheme $Y$, restricted Yoneda $N\mapsto(P\mapsto\Hom_Y(P,N))$ presents $\mathcal D_{\mathrm{qc}}(Y)$, the derived category of quasi-coherent complexes, as the module category on $\Perf(Y)$.
Lemma~\ref{lem:discrete-cauchy} therefore identifies the Cauchy objects of $\mathcal D((k,\Z)_\square)\otimes_{\mathcal D(k)}\mathcal D_{\mathrm{qc}}(Y)$ with $\delta\Perf(Y)$.

\subsection{The punctured completion}
\label{subsec:punctured-completion}

Use the full thick inclusion of \S~\ref{sec:support}, with source over $\mathcal D(k)$.
The following quotient and its enriched Cauchy completion are enriched over $\mathcal D((k,\Z)_\square)$, and the quotient bimodules take values in this coefficient category.
Let $\mathcal H$ be the stable idempotent-complete envelope of $\widehat{\mathcal C}_{\mathcal T}^{\square}$ inside $\RMod_{\delta\mathcal T}(\mathcal D((k,\Z)_\square))$.
Thus $\mathcal H$ is the closure of the discrete restricted modules $\delta r_{\mathcal T}X$ under finite cofibers, shifts, and retracts, with the full internal mapping objects of this module category.
For $t\in\mathcal T$, put $h_t=\delta\mathcal T(-,t)$.
The completion functor of~\eqref{eq:support-completion-functor} restricts to a fully faithful exact functor
\[
 j:\delta\mathcal T\longrightarrow\mathcal H,
 \qquad t\longmapsto h_t,
\]
since $\mathcal H(h_t,\delta r_{\mathcal T}X)\simeq\delta\mathcal C(t,X)$.
We use the same notation for this full inclusion when $\mathcal H$ is any small full stable subcategory of right $\delta\mathcal T$-modules containing the representables.

\begin{definition}\label{def:punctured-completion}
Define the enriched stable quotient $\mathcal Q_0$ and its enriched Cauchy completion $\mathcal Q$ by
\begin{equation}\label{eq:punctured-completion}
 \mathcal Q_0=\mathcal H/\delta\mathcal T,
 \qquad
 \mathcal Q=\Cauchy(\mathcal Q_0).
\end{equation}
Write $q:\mathcal H\to\mathcal Q_0$ for the quotient functor.
More precisely, $\mathcal Q_0$ is the full subcategory of $\ker(j^*)$ in~\eqref{eq:punctured-module-kernel} generated, under finite cofibers and shifts, by the images of the representable $\mathcal H$-modules under the localization~\eqref{eq:punctured-module-triangle}.
The coend~\eqref{eq:punctured-factorization-coend} and the cofiber sequence~\eqref{eq:punctured-quotient-triangle} compute its mapping objects between the images of objects of $\mathcal H$.
These closure operations change neither the enriched module category nor the enriched Cauchy completion of the localized representables.
For the order in which Cauchy completion precedes quotienting, put
\begin{equation}\label{eq:formal-quotient-order}
 \mathcal K=\Cauchy(\mathcal H),
 \qquad
 \mathcal Q_{\mathrm{form}}=\Kar(\mathcal K/\delta\mathcal T).
\end{equation}
All these categories are enriched over $\mathcal D((k,\Z)_\square)$.
\end{definition}

The enriched quotient has the following description.
Define the $\mathcal H$-bimodule $I$ by
\begin{equation}\label{eq:punctured-factorization-coend}
 I(M,N)=\int^{t\in\delta\mathcal T}
 \mathcal H(M,h_t)\otimes_{\delta k}^{\square}\mathcal H(h_t,N),
 \qquad M,N\in\mathcal H.
\end{equation}
This is a derived coend over $\delta\mathcal T$.
For example, its two-sided bar construction in degree $n$ sums the objects
\[
 \mathcal H(M,h_{t_0})\otimes_{\delta k}^{\square}
 \delta\mathcal T(t_0,t_1)\otimes_{\delta k}^{\square}\cdots
 \otimes_{\delta k}^{\square}\delta\mathcal T(t_{n-1},t_n)
 \otimes_{\delta k}^{\square}\mathcal H(h_{t_n},N).
\]
The face maps use composition and the degeneracies use units.
Composition defines $\alpha:I\to\Delta_{\mathcal H}$, and the quotient mapping objects fit into the bimodule cofiber sequence
\begin{equation}\label{eq:punctured-quotient-triangle}
 I\xrightarrow{\alpha}\Delta_{\mathcal H}
 \longrightarrow(q^e)^*\Delta_{\mathcal Q_0}.
\end{equation}
Here $\Delta_{\mathcal H}(M,N)=\mathcal H(M,N)$ and $q^e$ is the induced functor on enveloping categories, with $\mathcal H^e=\mathcal H^{\op}\otimes_{\mathcal D((k,\Z)_\square)}\mathcal H$.
In particular, the quotient mapping object $\mathcal Q_0(qM,qN)$ is the cofiber of the composition map $\alpha_{M,N}:I(M,N)\to\mathcal H(M,N)$.
This is the enriched form of the quotient calculation in \cite[\S~3]{Drinfeld}; the completion-then-quotient construction for $\mathcal D(k)$-enriched categories is developed in \cite{Gao}.

To derive~\eqref{eq:punctured-quotient-triangle}, consider the localization of right $\mathcal H$-modules away from the image of $j_!$.
Restriction $j^*$ has a fully faithful enriched left adjoint $j_!$, and this localization has triangle
\begin{equation}\label{eq:punctured-module-triangle}
 j_!j^*R\longrightarrow R\longrightarrow q^*q_!R.
\end{equation}
Evaluating on representables recovers~\eqref{eq:punctured-quotient-triangle}.
Its local objects are exactly the modules restricting to zero on $\delta\mathcal T$, so
\begin{equation}\label{eq:punctured-module-kernel}
 \RMod_{\mathcal Q_0}\bigl(\mathcal D((k,\Z)_\square)\bigr)
 \simeq
 \ker\!\left(
 j^*:\RMod_{\mathcal H}\bigl(\mathcal D((k,\Z)_\square)\bigr)
 \longrightarrow\RMod_{\delta\mathcal T}\bigl(\mathcal D((k,\Z)_\square)\bigr)
 \right).
\end{equation}
Let $L$ denote the localization in~\eqref{eq:punctured-module-triangle} and put $h_x=\mathcal H(-,x)$.
For $R\in\ker j^*$, enriched adjunction identifies $\sHom_{\ker j^*}(Lh_x,R)$ with $R(x)$, so $Lh_x$ is Cauchy there.
Localized co-Yoneda expresses $R$ as the weighted colimit $\int^{x\in\mathcal H}Lh_x\otimes_{\delta k}^{\square}R(x)$.
Thus the localized representables present $\ker j^*$, and their finite stable envelope has the same module category by Proposition~\ref{prop:testing-envelope}.
When $\mathcal H$ is the stable idempotent-complete envelope of $\widehat{\mathcal C}_{\mathcal T}^{\square}$, restriction identifies its enriched right-module category with that of $\widehat{\mathcal C}_{\mathcal T}^{\square}$.
Under this equivalence, \eqref{eq:punctured-module-kernel} identifies $\RMod_{\mathcal Q_0}(\mathcal D((k,\Z)_\square))$ with the kernel of $u_*j^*$ from Proposition~\ref{prop:free-to-complete}.
The quotient functor is an enriched homological epimorphism:
\begin{equation}\label{eq:quotient-homological-epimorphism}
 \mathcal Q_0\otimes_{\mathcal H}^{\square}\mathcal Q_0
 \simeq\Delta_{\mathcal Q_0}.
\end{equation}
Indeed, its extension--restriction counit is invertible on every module, and evaluation on the representables identifies the counit with~\eqref{eq:quotient-homological-epimorphism}.
The same statement for left and right module actions makes restriction along $q^e$ fully faithful on bimodules.

\subsection{Quotienting before and after Cauchy completion}

\begin{proposition}\label{prop:cauchy-quotient}
Extension of scalars along $q$ induces a canonical fully faithful enriched functor
\begin{equation}\label{eq:cauchy-quotient-comparison}
 \iota_{\mathrm{form}}:\mathcal Q_{\mathrm{form}}
 =\Kar\bigl(\Cauchy(\mathcal H)/\delta\mathcal T\bigr)
 \longrightarrow\Cauchy(\mathcal H/\delta\mathcal T)=\mathcal Q.
\end{equation}
If $\mathcal H(h_t,x)$ is discrete for every $t\in\mathcal T$ and $x\in\mathcal H$, then point evaluation induces
\begin{equation}\label{eq:ordinary-formal-quotient}
 \begin{aligned}
 \Gamma\mathcal Q_0&\simeq(\Gamma\mathcal H)/\mathcal T,\\
 \Gamma\mathcal Q_{\mathrm{form}}
 &\simeq\Kar\bigl((\Gamma\mathcal K)/\mathcal T\bigr),
 \end{aligned}
\end{equation}
where the quotients on the right use the $\mathcal D(k)$-enrichment.
\end{proposition}

\begin{proof}
The functor $q_!$ preserves Cauchy modules because its right adjoint $q^*$ is linear and preserves colimits.
Let $P,R$ be Cauchy right $\mathcal H$-modules.
Apply $\sHom_{\mathcal H}(P,-)$ to~\eqref{eq:punctured-module-triangle} for $R$.
Its linearity and preservation of colimits identify the first term with
\[
 \int^{t\in\delta\mathcal T}
 \sHom_{\mathcal H}(P,h_{jt})\otimes_{\delta k}^{\square}R(jt),
\]
where $h_{jt}=\mathcal H(-,jt)$ and $R(jt)=\sHom_{\mathcal H}(h_{jt},R)$.
This is the coend~\eqref{eq:punctured-factorization-coend} with $\mathcal H$ replaced by $\mathcal K$, evaluated at $(P,R)$.
Adjunction identifies the last term with $\sHom_{\mathcal Q_0}(q_!P,q_!R)$.
The resulting cofiber comparison proves full faithfulness between the images of the objects of $\mathcal K$ in $\mathcal K/\delta\mathcal T$.
Taking finite stable envelopes and then retracts extends this to the fully faithful functor $\iota_{\mathrm{form}}$ in~\eqref{eq:cauchy-quotient-comparison}.

Assume now that $\mathcal H(h_t,x)$ is discrete for every $t\in\mathcal T$ and $x\in\mathcal H$, as for the stable idempotent-complete envelope of $\widehat{\mathcal C}_{\mathcal T}^{\square}$ and the full stable subcategories of discrete pointwise perfect modules in \S~\ref{sec:punctured-duality}.
Lemma~\ref{lem:cauchy-discrete-evaluation}, applied to $\mathcal H$ and $b=jt$, shows that $\mathcal K(h_{jt},P)\simeq P(jt)$ is discrete for every $P\in\mathcal K$.
For $M\in\mathcal D((k,\Z)_\square)$ and $N\in\mathcal D(k)$, the restricted projection formula is
\begin{equation}\label{eq:point-projection}
 \Gamma(M\otimes_{\delta k}^{\square}\delta N)
 \simeq\Gamma(M)\otimes_kN.
\end{equation}
Both sides preserve colimits in $N$, and the comparison is the identity on the generator $N=k$.
In the quotient bar constructions for $\mathcal H$ and $\mathcal K$, the middle factors are Homs in $\delta\mathcal T$.
The final factor is $\mathcal H(h_t,N)$ or $\mathcal K(h_{jt},R)$, respectively; the first is discrete by hypothesis and the second by Lemma~\ref{lem:cauchy-discrete-evaluation}.
Only the first mixed Hom, $\mathcal H(M,h_t)$ or $\mathcal K(P,h_{jt})$, can be nondiscrete.
Applying $\Gamma$ to the two solid bar constructions and using~\eqref{eq:point-projection} recovers the bar constructions in $\mathcal D(k)$ for $(\Gamma\mathcal H)/\mathcal T$ and $(\Gamma\mathcal K)/\mathcal T$, respectively.
Since $\Gamma$ preserves colimits, these identifications also identify the realizations, composition maps, and their cofibers.
Passing to finite stable envelopes, and to retracts for the quotient $\mathcal K/\delta\mathcal T$, proves~\eqref{eq:ordinary-formal-quotient}.
\end{proof}

\section{Internal solid Serre duality}
\label{sec:punctured-duality}

We prove internal Serre duality by applying biduality to the values $I(M,N)$ of the bimodule defined in~\eqref{eq:punctured-factorization-coend}.
Under the hypotheses of Theorem~\ref{thm:punctured-solid-serre}, Lemma~\ref{lem:serre-mixed-maps} identifies these values with the discrete realizations in~\eqref{eq:serre-discrete-coend}.
The underlying complexes may be unbounded, so we begin with internal biduality for arbitrary discrete complexes.\footnote{P\"utzst\"uck announced reflexivity of discrete complexes for condensed $\Z$-linear duality in joint work with Nikolaus \cite{PutzstuckReflexivity}.
Their paper \cite[reference~{[NP]}]{NikolausPutzstuck} lists \emph{Condensed Anderson duality} as in preparation.}
Example~\ref{ex:unbounded-frobenius} illustrates the resulting duality for a quotient endomorphism complex with cohomology in arbitrarily high and low degrees.

Throughout this section, $k=\Z$ or $\F_p$, and duality is taken in $\mathcal D((k,\Z)_\square)$:
\[
 \sqdual F
 :=\sHom_{\mathcal D((k,\Z)_\square)}(F,\delta k).
\]
We write $\ordinarydual E=\Hom_k(E,k)$ for derived $k$-linear duality.
The category $\mathcal T$, its Serre structure, and its right modules before discrete realization are defined over $\mathcal D(k)$.
The full subcategory $\mathcal H$, the quotient $\mathcal Q_0=\mathcal H/\delta\mathcal T$, and the enriched Cauchy completion $\mathcal Q=\Cauchy(\mathcal Q_0)$ are enriched over $\mathcal D((k,\Z)_\square)$.
Their mapping objects, the coends defining the quotient, and the internal duals $\sqdual$ are formed in this solid coefficient category.

\subsection{Canonical biduality for unbounded discrete complexes}

\begin{lemma}
\label{lem:unbounded-solid-biduality}
For every $E\in\mathcal D(k)$, the canonical evaluation map
\begin{equation}
 \label{eq:canonical-solid-biduality}
 \eta_{\delta E}:\delta E\longrightarrow\sqdual^{\,2}(\delta E)
\end{equation}
is an equivalence.
\end{lemma}

\begin{proof}
We first take $k=\Z$.
If $F=\bigoplus_J\Z$ is free, then
\[
 \sqdual(\delta F)\simeq\prod_J\delta\Z.
\]
For a profinite set $P$, put $G_P=\Z[P]^\square$.
Each $G_P$ is a compact projective solid abelian group, and the $G_P$ form a compact generating family of $\DS$ \cite[Theorem~5.8 and Corollary~6.1(i),(ii)]{Scholze}.
The group $C(P,\Z)$ is free, and a choice of basis $I$ identifies $G_P$ with $\prod_I\delta\Z$ \cite[Theorem~5.4 and Corollary~5.5]{Scholze}.
Computing the tensor product and derived Hom, we obtain
\begin{align*}
 \Hom_{\DS}\bigl(G_P,\sqdual^{\,2}(\delta F)\bigr)
 &\simeq
 \Hom_{\DS}\left(
 G_P\otimes_{\Z_\square}\prod_J\delta\Z,\delta\Z\right)\\
 &\simeq C(P,\Z)\otimes_\Z F
 \simeq\Hom_{\DS}(G_P,\delta F).
\end{align*}
Indeed, Proposition~6.3 of \cite{Scholze} identifies the tensor product with $\prod_{I\times J}\delta\Z$, and the proof of Proposition~5.7 there computes its derived mapping dual as $\bigoplus_{I\times J}\Z$.
The displayed Homs are derived mapping complexes in $\DZ$.
Compact projectivity of the product sources excludes higher cohomology in the mapping complexes to $\delta\Z$ used in this calculation.
Under these identifications, $\Hom_{\DS}(G_P,\eta_{\delta F})$ is the identity on $C(P,\Z)\otimes_\Z F$, since the identifications use the coordinate evaluation pairings.
Thus $\eta_{\delta F}$ is an equivalence, because the $G_P$ generate $\DS$.

Every abelian group $A$ has a free resolution of length one.
Apply naturality of biduality to the cofiber sequence obtained from this resolution.
The free-group case and exactness of $\sqdual^{\,2}$ show that $\eta_{\delta A}$ is an equivalence.
Dualizing the resolution once also shows that $\sqdual(\delta A)$ has cohomological amplitude $[0,1]$.
Finite truncation triangles then prove biduality for every bounded complex.

For an arbitrary complex $E$, we construct the noncanonical decomposition
\begin{equation}
 \label{eq:hereditary-complex-splitting}
 E\simeq\bigoplus_{n\in\Z}H^n(E)[-n].
\end{equation}
Choose a free resolution $0\to F_{1,n}\to F_{0,n}\to H^n(E)\to0$ for each $n$.
Let $P_n$ be the two-term complex with $F_{1,n}$ and $F_{0,n}$ in degrees $n-1$ and $n$.
Projectivity lifts $F_{0,n}\to H^n(E)$ to the cycles of $E^n$.
The composite from $F_{1,n}$ then lands in the boundaries, and projectivity lifts it to $E^{n-1}$.
These lifts define a chain map $\phi_n:P_n\to E$ that induces the identity on $H^n(E)$.
The maps $\phi_n$ jointly induce
\[
 \phi:\bigoplus_{n\in\Z}P_n\longrightarrow E.
\]
The degree-$q$ term of the source complex is $F_{0,q}\oplus F_{1,q+1}$.
Only $P_q$ contributes to its degree-$q$ cohomology, and $H^q(\phi)$ is the identity on $H^q(E)$.
Thus $\phi$ is a quasi-isomorphism.
Since $P_n\simeq H^n(E)[-n]$, the quasi-isomorphism $\phi$ proves \eqref{eq:hereditary-complex-splitting}.

The object $\sqdual(\delta H^n(E))[n]$ has cohomological amplitude $[-n,1-n]$.
For each fixed $q$, only $n=-q$ and $n=1-q$ can contribute to its degree-$q$ cohomology.
Products and direct sums are exact in solid abelian groups \cite[Theorems~1.10 and~5.8(i)]{Scholze}, so the canonical map
\[
 \bigoplus_n\sqdual(\delta H^n(E))[n]
 \longrightarrow
 \prod_n\sqdual(\delta H^n(E))[n]
\]
is an equivalence.
For the family $\delta H^n(E)[-n]$, only $n=q$ contributes to degree-$q$ cohomology, so the canonical map
\[
 \bigoplus_n\delta H^n(E)[-n]
 \longrightarrow\prod_n\delta H^n(E)[-n]
\]
is also an equivalence.
Dualizing \eqref{eq:hereditary-complex-splitting} twice and using these two comparisons, we obtain
\[
 \sqdual^{\,2}(\delta E)
 \simeq\prod_n\sqdual^{\,2}\bigl(\delta H^n(E)[-n]\bigr)
 \simeq\prod_n\delta H^n(E)[-n]
 \simeq\delta E.
\]
To identify the canonical evaluation map, first apply this calculation to the direct sum in~\eqref{eq:hereditary-complex-splitting}.
Naturality with respect to its projections shows that the $n$th component of its biduality map is the projection to $\delta H^n(E)[-n]$ followed by the bounded evaluation equivalence for that summand.
The two sum-to-product comparisons therefore imply that this canonical biduality map is an equivalence.
Naturality along the chosen decomposition proves the assertion for $\eta_{\delta E}$.

For $k=\F_p$, let
\[
 U:\mathcal D((\F_p,\Z)_\square)
   =\Mod_{\delta\F_p}(\DS)
   \longrightarrow\DS
\]
be restriction of scalars along $\delta\Z\to\delta\F_p$.
Write
\[
 \mathbb D_{\Z,\square}(M)=\sHom_{\DS}(M,\delta\Z)
\]
for integral internal duality, and use $\sqdual$ for internal duality relative to $\delta\F_p$.
Put $\omega=\Hom_\Z(\F_p,\Z)\simeq\F_p[-1]\in\mathcal D(\F_p)$, with $\Hom_\Z$ denoting derived Hom.
Since $\F_p$ is perfect over $\Z$, coinduction sends $\delta\Z$ to the solid $\F_p$-module $\delta\omega$.
By restriction--coinduction adjunction,
\[
 \mathbb D_{\Z,\square}(UN)
 \simeq U\sHom_{\delta\F_p}(N,\delta\omega)
\]
for every $N\in\mathcal D((\F_p,\Z)_\square)$.
The object $\delta\omega$ is invertible, so
\[
 \sHom_{\delta\F_p}(N,\delta\omega)
 \simeq\sqdual N\otimes^{\square}_{\delta\F_p}\delta\omega.
\]
Applying restriction--coinduction twice, then using tensor--Hom adjunction and $\sHom_{\delta\F_p}(\delta\omega,\delta\omega)\simeq\delta\F_p$, we obtain
\[
 \mathbb D_{\Z,\square}^{\,2}(UN)
 \simeq
 U\sHom_{\delta\F_p}
   \bigl(\sHom_{\delta\F_p}(N,\delta\omega),\delta\omega\bigr)
 \simeq U\sqdual^{\,2}N.
\]
Under these adjunctions, evaluation with values in $\delta\omega$ is $\delta\F_p$-valued evaluation tensored with $\delta\omega$.
Contracting the invertible object $\delta\omega$ by its canonical evaluation therefore identifies the integral biduality map with $U\eta_N$.
For $N=\delta E$, the solid $\Z$-complex $UN$ is the discrete realization of the underlying $\Z$-complex of $E$.
The integral case and conservativity of $U$ therefore prove that $\eta_N$ is an equivalence.
\end{proof}

\subsection{The internal dual of the coend \texorpdfstring{$I$}{I}}

Let $\mathcal T$ be a small stable idempotent-complete proper $k$-linear category.
Suppose that it has a coherent Serre autoequivalence $S_{\mathcal T}$, that is, a natural equivalence
\begin{equation}
 \label{eq:testing-serre-structure}
 \ordinarydual\mathcal T(t,u)
 \simeq\mathcal T(u,S_{\mathcal T}t)
\end{equation}
of $\mathcal D(k)$-enriched bimodules.
In particular, its traces $\operatorname{tr}_t:\mathcal T(t,S_{\mathcal T}t)\to k$ satisfy
\begin{equation}
 \label{eq:testing-serre-cyclicity}
 \operatorname{tr}_t(gf)
 =\operatorname{tr}_u\bigl((S_{\mathcal T}f)g\bigr),
 \qquad f:t\to u,\quad g:u\to S_{\mathcal T}t.
\end{equation}
This is an identity of pairings, with the symmetric braiding included.
After realization by $\delta$, the Serre traces take values in $\delta k$ and satisfy the realized identity of pairings.

A right $\mathcal T$-module $M$ is \emph{pointwise perfect} if $M(t)\in\Perf(k)$ for every $t$.
Let $\PsPerf^{\square}(\mathcal T)$ be the full $\mathcal D((k,\Z)_\square)$-enriched subcategory of $\RMod_{\delta\mathcal T}(\mathcal D((k,\Z)_\square))$ on the discrete realizations of these modules.
Its mapping objects are the internal mapping objects in that solid module category.
In $\mathcal H$ and its solid quotients, we write $M$ for the realization $\delta M$ of a $\mathcal D(k)$-valued right $\mathcal T$-module.
Values $M(t)$ and tensors over $k$ inside $\delta(\cdot)$ refer to the module before discrete realization.
Put
\[
 h_t=\delta\mathcal T(-,t),\qquad
 \sigma M=M\circ S_{\mathcal T}^{-1}.
\]
Thus $\sigma h_t\simeq h_{S_{\mathcal T}t}$, and $\sigma$ is an enriched autoequivalence of $\PsPerf^{\square}(\mathcal T)$.

\begin{theorem}
\label{thm:punctured-solid-serre}
\label{thm:pseudoperfect}
Let $k=\Z$ or $\F_p$ and let $\mathcal T$ be small, stable, idempotent complete and proper over $k$, with a coherent $\mathcal D(k)$-enriched Serre autoequivalence $S_{\mathcal T}$.
Let $\mathcal H$ be a small full $\mathcal D((k,\Z)_\square)$-enriched stable subcategory of $\PsPerf^{\square}(\mathcal T)$ containing the representables and invariant under $\sigma=(-)\circ S_{\mathcal T}^{-1}$.
Over the same relative solid base form
\[
 \mathcal Q_0=\mathcal H/\delta\mathcal T,
 \qquad
 \mathcal Q=\Cauchy(\mathcal Q_0).
\]
There is an equivalence of $\mathcal D((k,\Z)_\square)$-valued $\mathcal Q_0$-bimodules
\begin{equation}
 \label{eq:punctured-solid-serre}
 \sqdual\mathcal Q_0(M,N)
 \simeq\mathcal Q_0(N,\sigma M)[-1].
\end{equation}
It extends to every pair of objects of $\mathcal Q$ and restricts along $\iota_{\mathrm{form}}:\mathcal Q_{\mathrm{form}}\hookrightarrow\mathcal Q$ in~\eqref{eq:cauchy-quotient-comparison}.
The extension of $\sigma[-1]$ is the internal solid Serre autoequivalence of $\mathcal Q$.
\end{theorem}

Let $I$ be the bimodule defined by the coend~\eqref{eq:punctured-factorization-coend} in $\mathcal D((k,\Z)_\square)$, with composition map $\alpha:I\to\Delta_{\mathcal H}$ as in~\eqref{eq:punctured-quotient-triangle}.
The coends indexed by $\mathcal T$ inside $\delta(\cdot)$ below are formed in $\mathcal D(k)$.

\begin{lemma}
\label{lem:serre-mixed-maps}
For $M,N\in\mathcal H$, there are natural equivalences
\begin{align}
 \mathcal H(h_t,N)&\simeq\delta(N(t)),
 \label{eq:serre-forward-mixed}\\
 \mathcal H(M,h_t)&\simeq
 \delta\bigl(\ordinarydual M(S_{\mathcal T}^{-1}t)\bigr).
 \label{eq:serre-reverse-mixed}
\end{align}
Consequently
\begin{equation}
 \label{eq:serre-discrete-coend}
 I(M,N)\simeq\delta\left(\int^{t\in\mathcal T}
 \ordinarydual M(S_{\mathcal T}^{-1}t)\otimes_k N(t)\right).
\end{equation}
In particular, the canonical biduality map of $I$ is an equivalence.
\end{lemma}

\begin{proof}
The first formula is Yoneda over $\mathcal D((k,\Z)_\square)$.
For the second, let $\ell_v=\delta\mathcal T(v,-)$ be the left representable.
The Serre structure identifies $h_t$ with the pointwise internal dual of $\ell_{S_{\mathcal T}^{-1}t}$.
By relative tensor--Hom adjunction and co-Yoneda,
\begin{align*}
 \mathcal H(M,h_t)
 &\simeq\sqdual\bigl(
 \ell_{S_{\mathcal T}^{-1}t}\otimes_{\delta\mathcal T}^\square M
 \bigr)\\
 &\simeq\sqdual\delta\bigl(M(S_{\mathcal T}^{-1}t)\bigr)\\
 &\simeq\delta\bigl(\ordinarydual M(S_{\mathcal T}^{-1}t)\bigr).
\end{align*}
The last equivalence uses pointwise perfection.
Substituting these equivalences in~\eqref{eq:punctured-factorization-coend}, we obtain~\eqref{eq:serre-discrete-coend}, since the strong symmetric monoidal left adjoint $\delta$ preserves the tensors and realization of the two-sided bar construction.
The coend in $\mathcal D(k)$ can be unbounded even though its factors are pointwise perfect.
The right-hand side of~\eqref{eq:serre-discrete-coend} is discrete, so Lemma~\ref{lem:unbounded-solid-biduality} proves that the canonical biduality map of $I(M,N)$ is an equivalence.
\end{proof}

Changing variables by $S_{\mathcal T}$ rewrites \eqref{eq:serre-discrete-coend} as
\begin{equation}
 \label{eq:serre-reindexed-coend}
 I(N,\sigma M)
 \simeq\delta\left(\int^{t\in\mathcal T}
 \ordinarydual N(t)\otimes_k M(t)\right).
\end{equation}
Let $N^D$ denote the left solid module obtained by pointwise internal duality from $N$.
Its values are $\delta(\ordinarydual N(t))$.
The right side of \eqref{eq:serre-reindexed-coend} is $N^D\otimes_{\delta\mathcal T}^\square M$.
Pointwise perfection identifies $\sqdual N^D$ with $N$.
Tensor--Hom adjunction now induces an equivalence
\begin{equation}
 \label{eq:serre-diagonal-dual-coend}
 \beta_{M,N}:\mathcal H(M,N)
 \xrightarrow{\ \sim\ }\sqdual I(N,\sigma M).
\end{equation}
These are equivalences of enriched bimodules: pointwise duality, module evaluation, and relative tensor--Hom adjunction are enriched constructions.
Their naturality is expressed on the full internal mapping objects of $\mathcal H$, which may be nondiscrete.

\subsection{Duality for the quotient bimodule}

We describe the pairing adjoint to $\beta_{M,N}$ in~\eqref{eq:serre-diagonal-dual-coend} on the tensor products defining the coend~\eqref{eq:punctured-factorization-coend}.
In the notation $h:M\to N$, $c:N\to h_u$, and $d:h_u\to\sigma M$, it is given by
\begin{equation}
 \label{eq:serre-balanced-pairing}
 \langle h,[c,d]\rangle
 =\operatorname{tr}_{S_{\mathcal T}^{-1}u}
       \bigl(c h\,\sigma^{-1}d\bigr).
\end{equation}
The displayed composite is a map $h_{S_{\mathcal T}^{-1}u}\to h_u$.
Its trace is the realization of $\operatorname{tr}_{S_{\mathcal T}^{-1}u}$, with target $\delta k$.
Formula \eqref{eq:serre-balanced-pairing} denotes a morphism from the corresponding solid tensor product to $\delta k$.
Naturality of the Serre bimodule structure makes the pairing~\eqref{eq:serre-balanced-pairing} balanced in $u$.
It therefore induces a pairing
\[
 \mathcal H(M,N)\otimes_{\delta k}^{\square}I(N,\sigma M)
 \longrightarrow\delta k,
\]
whose adjoint is the equivalence $\beta_{M,N}$ in~\eqref{eq:serre-diagonal-dual-coend}.

The bimodules $I$ and $\Delta_{\mathcal H}$ carry their natural diagonal equivariances under $\sigma$: on $I$ one applies $\sigma$ to the two mixed mapping objects and reindexes $t$ by $S_{\mathcal T}t$.
Composition is equivariant for these structures.
The same holds for $\beta$.
To see the trace compatibility explicitly, specialize \eqref{eq:testing-serre-cyclicity} to obtain
\[
 \operatorname{tr}_{S_{\mathcal T}t}(S_{\mathcal T}f)
 =\operatorname{tr}_t(f),\qquad f:t\to S_{\mathcal T}t.
\]
Insert the $\mathcal D(k)$-enriched unit $k\to\mathcal T(S_{\mathcal T}t,S_{\mathcal T}t)$ in~\eqref{eq:testing-serre-cyclicity} with $u=S_{\mathcal T}t$.
After applying $\delta$, the same identity holds over the solid coefficient category.
Substituting it in~\eqref{eq:serre-balanced-pairing} identifies that pairing with its simultaneous $\sigma$-translate, with the coherences supplied by~\eqref{eq:testing-serre-structure}.

For such equivariant bimodules, use the contravariant twisted dual
\[
 d_\sigma K(M,N)=\sqdual K(N,\sigma M).
\]
Reversing the variance twice, we compute
\[
 d_\sigma^2K(M,N)
 =\sqdual^{\,2}K(\sigma M,\sigma N)
 \simeq\sqdual^{\,2}K(M,N),
\]
where the last equivalence is the diagonal equivariance.
Via this identification, internal evaluation defines the canonical map $\eta_K:K\to d_\sigma^2K$.
Equation~\eqref{eq:serre-diagonal-dual-coend} is $\beta:\Delta_{\mathcal H}\xrightarrow{\sim}d_\sigma I$.
Define its transpose by
\begin{equation}
 \label{eq:serre-transpose-rho}
 \rho=(d_\sigma\beta)\eta_I:
 I\longrightarrow d_\sigma\Delta_{\mathcal H}.
\end{equation}
Both factors are equivalences: $\eta_I$ by Lemma~\ref{lem:serre-mixed-maps}, and $d_\sigma\beta$ because $\beta$ is an equivalence.

\begin{lemma}
\label{lem:serre-transpose-square}
There is a commuting square of equivariant enriched bimodules
\begin{equation}
 \label{eq:serre-transpose-square}
 \begin{CD}
 I @>{\alpha}>> \Delta_{\mathcal H}\\
 @V{\rho}V{\simeq}V @V{\beta}V{\simeq}V\\
 d_\sigma\Delta_{\mathcal H}
       @>{d_\sigma\alpha}>> d_\sigma I.
 \end{CD}
\end{equation}
\end{lemma}

\begin{proof}
Consider two factorizations
\[
 a:M\to h_t,\quad b:h_t\to N,\qquad
 c:N\to h_u,\quad d:h_u\to\sigma M.
\]
The two sides of the asserted square, paired with the second factorization, are
\begin{equation}
 \label{eq:serre-two-transposes}
 \operatorname{tr}_{S_{\mathcal T}^{-1}u}
       (cba\,\sigma^{-1}d),
 \qquad
 \operatorname{tr}_t\bigl((\sigma a)dcb\bigr).
\end{equation}
For the second expression, the definition of $\rho$ first transports $[a,b]$ to $[\sigma a,\sigma b]$ by the diagonal equivariance and then evaluates $\beta_{N,\sigma M}$.
This is the pairing obtained from $\rho$ as defined in~\eqref{eq:serre-transpose-rho}.

Apply the trace identity~\eqref{eq:testing-serre-cyclicity}, after realization by $\delta$, to the two maps between representables
\[
 a\,\sigma^{-1}d:h_{S_{\mathcal T}^{-1}u}\longrightarrow h_t,
 \qquad cb:h_t\longrightarrow h_u.
\]
It identifies the expressions in \eqref{eq:serre-two-transposes}.
More precisely, both expressions are maps from
\[
 \mathcal H(M,h_t)\otimes_{\delta k}^{\square}\mathcal H(h_t,N)
 \otimes_{\delta k}^{\square}\mathcal H(N,h_u)
 \otimes_{\delta k}^{\square}\mathcal H(h_u,\sigma M)
 \longrightarrow\delta k.
\]
They are the two sides of the naturality identity for the Serre pairing, after composing the indicated pairs of factors.
The symmetric braiding in this identity supplies the Koszul signs.
The identity is balanced in both variables $t,u\in\mathcal T$.
Applying this identity to the bar constructions, with their iterated compositions and unit maps, therefore induces a compatible homotopy in every simplicial degree.
The coherence of \eqref{eq:testing-serre-structure} supplies the compatibilities among these homotopies.
Its geometric realization is the asserted square of bimodule maps.
\end{proof}

Put $B=(q^e)^*\Delta_{\mathcal Q_0}$.
By~\eqref{eq:punctured-quotient-triangle} this is the cofiber of $\alpha:I\to\Delta_{\mathcal H}$.
The vertical equivalences in~\eqref{eq:serre-transpose-square} identify $B=\cofib(\alpha)$ with $\cofib(d_\sigma\alpha)$.
Since $d_\sigma$ is contravariant and exact, we have
\[
 \cofib(d_\sigma\alpha)\simeq(d_\sigma B)[1],
\]
and hence
\begin{equation}
 \label{eq:serre-quotient-transpose}
 B(M,N)\simeq\sqdual B(N,\sigma M)[1].
\end{equation}
Replace $(M,N)$ here by $(N,\sigma M)$ and use diagonal $\sigma$-equivariance.
The result is
\[
 B(N,\sigma M)\simeq\sqdual B(M,N)[1],
\]
which is precisely \eqref{eq:punctured-solid-serre} after shifting.
In particular, the shift $[-1]$ comes from the dual of the quotient cofiber, with no boundedness assumption on its terms.

Equation~\eqref{eq:quotient-homological-epimorphism}, applied in both variables, implies that $(q^e)^*$ is fully faithful on $\mathcal D((k,\Z)_\square)$-valued bimodules.
Since $\sigma$ preserves $\delta\mathcal T$, it acts on the images of the objects of $\mathcal H$ in the quotient and on their finite stable envelope $\mathcal Q_0$.
Both sides of~\eqref{eq:punctured-solid-serre} therefore have quotient bimodule structures, and the equivalence just constructed descends uniquely to them.

\subsection{Extension of duality to Cauchy modules}

We extend~\eqref{eq:punctured-solid-serre} to Cauchy modules by using the Yoneda presentation~\eqref{eq:absolute-yoneda}.
For $x\in\mathcal Q_0$, write $h_x=\mathcal Q_0(-,x)$ for its right representable module.
Define, on the whole category of right solid $\mathcal Q_0$-modules,
\[
 P^\vee(x)=\sHom_{\mathcal Q_0}(P,h_x),
 \qquad \nu(P)=\sqdual(P^\vee).
\]
Here $P\mapsto P^\vee$ is contravariant from right modules to left modules, and pointwise $\sqdual$ is contravariant from left modules to right modules.
Both are enriched over $\mathcal D((k,\Z)_\square)$ on the corresponding opposite categories, so $\nu$ is a covariant functor enriched over this base.

If $P$ is Cauchy, its internal Hom functor is linear over the solid coefficient category and preserves colimits.
Applying this functor to the weighted Yoneda presentation of an arbitrary right module $R$, we obtain
\begin{equation}
 \label{eq:cauchy-serre-tensor-hom}
 \sHom_{\mathcal Q_0}(P,R)
 \simeq P^\vee\otimes_{\mathcal Q_0}^\square R.
\end{equation}
On representables, \eqref{eq:punctured-solid-serre} identifies
\[
 \nu(h_x)(y)=\sqdual\mathcal Q_0(x,y)
 \simeq\mathcal Q_0(y,\sigma x)[-1].
\]
Thus $\nu(h_x)\simeq h_{\sigma x}[-1]$, naturally as an enriched functor of $x$.

Since $\nu$ is covariant and enriched, it preserves the absolute Yoneda presentation~\eqref{eq:absolute-yoneda} of the Cauchy module $P$.
Substituting the values of $\nu$ on representables, we obtain
\begin{equation}
 \label{eq:cauchy-serre-nakayama}
 \nu(P)\simeq\int^{x\in\mathcal Q_0}
 P(x)\otimes_{\delta k}^{\square} h_{\sigma x}[-1]
 \simeq(\sigma_!P)[-1],
\end{equation}
where $\sigma_!$ is extension of scalars along $\sigma$.
Dualizing \eqref{eq:cauchy-serre-tensor-hom} and applying relative tensor--Hom adjunction, we deduce
\begin{align}
 \sqdual\sHom_{\mathcal Q_0}(P,R)
 &\simeq\sHom_{\mathcal Q_0}(R,\nu(P))\notag\\
 &\simeq\sHom_{\mathcal Q_0}(R,\sigma_!P)[-1].
 \label{eq:cauchy-solid-serre}
\end{align}
All the identifications are enriched and natural in both variables.
The autoequivalence $\sigma_!$ preserves Cauchy modules, as does its inverse.
This proves~\eqref{eq:cauchy-solid-serre} on $\mathcal Q$.
Since $\sigma$ preserves the image of $\delta\mathcal T\hookrightarrow\mathcal H$, extension along $q$ commutes with $\sigma$ and with shifts.
The essential image of the fully faithful functor $\iota_{\mathrm{form}}$ in~\eqref{eq:cauchy-quotient-comparison} is therefore preserved by $\sigma_![-1]$, so the duality restricts as asserted in Theorem~\ref{thm:punctured-solid-serre}.

\subsection{Proper source categories and the Calabi--Yau case}

\begin{proof}[Proof of Theorem~\ref{thm:intro-punctured-serre}]
Let $S$ be the coherent Serre autoequivalence of the proper category $\mathcal C$, and suppose $S(\mathcal T)=\mathcal T$.
Its restriction is a coherent Serre autoequivalence of $\mathcal T$.
Properness makes every restricted right Yoneda module $r_X=\mathcal C(-,X)|_{\mathcal T}$ pointwise perfect, and full faithfulness of $S$ implies
\[
 r_{SX}\simeq r_X\circ S_{\mathcal T}^{-1}.
\]
The stable idempotent-complete envelope of the discrete modules $\delta r_X$ still consists of discrete pointwise perfect modules.
Indeed, the discrete realization on module categories is fully faithful on the underlying $\infty$-categories and exact, while evaluation preserves finite cofibers and retracts.
This envelope is $\sigma$-invariant, so Theorem~\ref{thm:punctured-solid-serre} applies.

The restriction to $\mathcal Q_{\mathrm{form}}$ supplied by Theorem~\ref{thm:punctured-solid-serre} proves the remaining assertion of Theorem~\ref{thm:intro-punctured-serre}.
\end{proof}

Write $\HH_\bullet(\mathcal T/k)$ for the Hochschild chain complex of $\mathcal T$ over $k$.

\begin{corollary}
\label{cor:punctured-cy}
Suppose that $\mathcal T$ is proper over $k=\Z$ or $\F_p$ and carries a weak proper $d$-Calabi--Yau trace $\HH_\bullet(\mathcal T/k)\to k[-d]$ whose adjoint is a $\mathcal D(k)$-enriched bimodule equivalence
\[
 \Delta_{\mathcal T}[d]\simeq\ordinarydual\Delta_{\mathcal T},
 \qquad
 (\ordinarydual\Delta_{\mathcal T})(t,u)
       =\ordinarydual\mathcal T(u,t).
\]
Let $\mathcal H$ be any small full stable subcategory of $\PsPerf^{\square}(\mathcal T)$ containing the representables, and put
\[
 \mathcal Q_0=\mathcal H/\delta\mathcal T,
 \qquad\mathcal Q=\Cauchy(\mathcal Q_0),
\]
over $\mathcal D((k,\Z)_\square)$.
For $P,R\in\mathcal Q$, there is a natural internal duality
\begin{equation}
 \label{eq:cy-bimodule}
 \sqdual\mathcal Q(P,R)\simeq\mathcal Q(R,P)[d-1].
\end{equation}
The same formula holds for the quotient mapping objects after restriction to $\mathcal Q_0$.
\end{corollary}

\begin{proof}
The chosen trace makes $[d]$ a coherent Serre functor on $\mathcal T$.
Its induced action on right modules is $M\mapsto M\circ[-d]\simeq M[d]$.
Since $\mathcal H$ is stable, it is closed under the shift $[d]$ and hence invariant under this action.
The result follows from Theorem~\ref{thm:punctured-solid-serre}.
\end{proof}

\section{Adic completion}
\label{sec:adic}

We identify the solid categorical completion of $\Perf(R)$ along the perfect complexes generated by a Koszul complex.
We first compare derived completion in solid modules with derived completion in $\mathcal D(R)$.
Under weak proregularity, the quotient tower expresses the morphism objects of this categorical completion as scalar extensions.
In this section and \S~\ref{sec:scalar-extension}, we use $\DZ$ as the coefficient category before discrete realization and $\DS$ for the solid construction; modules carry additional $R$- or $\delta R$-actions, which we forget when specifying coefficient-valued mapping objects.

\subsection{Koszul description}

Let $R$ be a commutative ring and let $I=(f_1,\ldots,f_r)$ be a finitely generated ideal.
For $n\geq1$, set
\begin{equation}
  K_n=\bigotimes_{i=1}^r\operatorname{cofib}(f_i^n:R\longrightarrow R),
  \qquad
  K_\infty=\colim_n K_n^\vee,
  \label{eq:koszul-objects}
\end{equation}
where the tensor products are over $R$, $K_n$ lies in degrees $[-r,0]$, and $K_n^\vee=\Hom_R(K_n,R)$.
On each two-term factor, the transition $K_{n+1}\to K_n$ is multiplication by $f_i$ in degree $-1$ and the identity in degree zero.
In particular,
\[
  K_\infty\simeq
  \bigotimes_{i=1}^r\operatorname{fib}\bigl(R\longrightarrow R[1/f_i]\bigr).
\]
The maps from these fibers to $R$ induce an augmentation $\epsilon:K_\infty\to R$.

For a solid $\delta R$-module $M$, define
\begin{equation}
  \Gamma_I M=\delta K_\infty\otimes_{\delta R}^{\square}M,
  \qquad
  \Lambda_I M=\sHom_{\delta R}(\delta K_\infty,M)
  \simeq\lim_n\bigl(\delta K_n\otimes_{\delta R}^{\square}M\bigr).
  \label{eq:koszul-completion}
\end{equation}
These are endofunctors of solid $\delta R$-modules.
The augmentation induces the counit $\Gamma_I M\to M$ and the completion unit
\begin{equation}
  \eta_M:M\simeq\sHom_{\delta R}(\delta R,M)
       \xrightarrow{\epsilon^*}\Lambda_I M.
  \label{eq:adic-completion-unit}
\end{equation}

For a profinite set $S$, write $M(S)$ for the complex of sections from \S~\ref{subsec:solid-coefficients}, with its induced $R$-action.
Using evaluation and tensor--Hom adjunction, we define the projection-formula map, which is an equivalence for every $V\in\mathcal D(R)$:
\begin{equation}\label{eq:adic-section-projection}
 V\otimes_R M(S)
 \xrightarrow{\ \sim\ }
 \bigl(\delta V\otimes_{\delta R}^{\square}M\bigr)(S).
\end{equation}
Indeed, both sides preserve colimits in $V$ and agree for $V=R$, which generates $\mathcal D(R)$; on the right, use compactness of $\Z[S]^\square$.
Evaluation $M\mapsto M(S)$ is right adjoint to
$V\mapsto\delta V\otimes_{\Z_\square}\Z[S]^\square$,
and therefore preserves limits.
Combining preservation of limits with~\eqref{eq:adic-section-projection}, we obtain
\[
 (\Lambda_I M)(S)\simeq\lim_n\bigl(K_n\otimes_R M(S)\bigr).
\]
The right-hand side is Koszul completion in $\mathcal D(R)$, the reflector onto derived $I$-complete complexes by \cite[Tag~0920]{Stacks}, and the identification carries $\eta_M(S)$ to its completion unit.

The augmented localization fibers identify $\Gamma_I$ with the colocalization onto the subcategory generated by $\delta K_1$ under colimits and the $\DS$-action.
Apply Proposition~\ref{thm:support-completion} to $\thick(K_1)\subseteq\Perf(R)$, identifying modules over $\delta\Perf(R)$ with solid $\delta R$-modules by Morita theory.
The right adjoint $\Lambda_I$ of $\Gamma_I$ is therefore the reflector onto the derived $I$-complete solid modules, characterized by
\begin{equation}
  \sHom_{\delta R}(\delta(R[1/f_i]),M)=0
  \qquad (1\leq i\leq r).
  \label{eq:adic-completeness-test}
\end{equation}
By tensor--Hom adjunction, this condition is equivalent to $\Hom_R(R[1/f_i],M(S))=0$ for every $i$ and every profinite $S$.
Thus a solid module is derived $I$-complete exactly when all its profinite section complexes are derived $I$-complete.
For $E\in\mathcal D(R)$, write $\Gamma_I E=K_\infty\otimes_R E$ and $\Lambda_I E=\Hom_R(K_\infty,E)$; here $\Gamma_I$ is the colocalization onto the localizing subcategory generated by $K_1$.
These functors, both in $\mathcal D(R)$ and in solid $\delta R$-modules, depend only on $\sqrt I$.

\begin{lemma}
\label{lem:adic-koszul-detection}
For every solid $\delta R$-module $M$, the map
\[
 \delta K_1\otimes_{\delta R}^{\square}\eta_M:
 \delta K_1\otimes_{\delta R}^{\square}M
 \longrightarrow
 \delta K_1\otimes_{\delta R}^{\square}\Lambda_I M
\]
is an equivalence.
A morphism between derived $I$-complete solid modules is an equivalence if and only if it becomes an equivalence after tensoring with $\delta K_1$ over $\delta R$.
The same assertions hold in $\mathcal D(R)$, with $K_1\otimes_R-$.
\end{lemma}

\begin{proof}
In $\mathcal D(R)$, Koszul completion is localization at the maps that become equivalences after tensoring with $K_1$ by \cite[Proposition~2.34]{BHV}; both assertions follow from this characterization.
Apply these assertions on every profinite set and use~\eqref{eq:adic-section-projection}, the sectionwise completion formula, and detection of equivalences on profinite sections.
\end{proof}

Lemma~\ref{lem:adic-koszul-detection} identifies $\ker\Lambda_I$ with the tensor ideal $\ker(\delta K_1\otimes_{\delta R}^{\square}-)$.
Hence completion is a symmetric monoidal localization \cite[Example~2.2.1.7 and Proposition~2.2.1.9]{LurieHA}.
Its target has tensor product $(M,N)\mapsto\Lambda_I(M\otimes_{\delta R}^{\square}N)$ and unit $\Lambda_I\delta R$.
The lax monoidal structure of the inclusion defines the pairings
\begin{equation}
  \Lambda_I M\otimes_{\delta R}^{\square}\Lambda_I N
       \longrightarrow\Lambda_I(M\otimes_{\delta R}^{\square}N),
  \label{eq:adic-completion-pairing}
\end{equation}
compatible with the completion unit, associativity, and symmetry.

The quotient maps define $M\to\lim_n(\delta(R/I^n)\otimes_{\delta R}^{\square}M)$.
Each term is derived $I$-complete because the $f_i$ act nilpotently, so the limit is complete.
The factorization through $\eta_M$ defines
\begin{equation}
  q_M:\Lambda_I M\longrightarrow
  \lim_n\bigl(\delta(R/I^n)\otimes_{\delta R}^{\square}M\bigr).
  \label{eq:quotient-completion}
\end{equation}
The quotient-limit pairing is defined by projection to each level and
the multiplication maps
\[
\begin{aligned}
 &\bigl(\delta(R/I^n)\otimes_{\delta R}^{\square}M\bigr)
  \otimes_{\delta R}^{\square}
  \bigl(\delta(R/I^n)\otimes_{\delta R}^{\square}N\bigr)\\
 &\qquad\longrightarrow
 \delta(R/I^n)\otimes_{\delta R}^{\square}
       \bigl(M\otimes_{\delta R}^{\square}N\bigr).
\end{aligned}
\]
The composite of the quotient-limit pairing with $(q_M\circ\eta_M)\otimes_{\delta R}^{\square}(q_N\circ\eta_N)$ is
$q_{M\otimes_{\delta R}^{\square}N}\circ
 \eta_{M\otimes_{\delta R}^{\square}N}$,
since both composites at level $n$ are induced by $R\to R/I^n$.
The same construction for any finite number of factors, with unit
$\delta R\to\lim_n\delta(R/I^n)$, is coherent because the quotient
transitions are morphisms of commutative algebras.
Precomposition with finite tensor products of completion units induces
equivalences of mapping spaces into complete objects, so this
compatibility extends uniquely and coherently to $q$.

We recall weak proregularity \cite[Definition~4.21]{PSYTorsion}, which suffices for~\eqref{eq:quotient-completion} to be an equivalence.

\begin{definition}[Weak proregularity]\label{def:weak-proregularity}
An inverse system $(M_n)_{n\geq1}$ of $R$-modules is \emph{pro-zero} if for each $n\geq1$ there exists $m\geq n$ such that the transition $M_m\to M_n$ is zero.
With $K_n$ as in~\eqref{eq:koszul-objects}, the sequence $(f_1,\ldots,f_r)$ is \emph{weakly proregular} if $(H^{-q}(K_n))_{n\geq1}$ is pro-zero for every $q>0$.
An ideal is \emph{weakly proregular} if it admits a finite generating sequence with this property.
\end{definition}

Every finite generating sequence of a weakly proregular ideal is weakly proregular \cite[Corollary~6.3]{PSYTorsion}.

\begin{proposition}
\label{prop:quotient-towers}
Assume that $I=(f_1,\ldots,f_r)$ is weakly proregular and put $I_n=(f_1^n,\ldots,f_r^n)$.
\begin{enumerate}
\item\label{it:quotient-pro-objects}
The augmentations $K_n\to R/I_n$ induce an equivalence of pro-objects in $\mathcal D(R)$.
\item\label{it:quotient-solid-limit}
For every $M\in\Mod_{\delta R}(\DS)$, the map $q_M$ in~\eqref{eq:quotient-completion} is an equivalence.
\end{enumerate}
In particular, both assertions hold for every ideal of a noetherian ring.
\end{proposition}

\begin{proof}
The fibers of $K_n\to R/I_n$ have cohomology in $[-r,-1]$
and pro-zero cohomology systems.
By \cite[Tag~0G3B]{Stacks}, for each fiber some later transition to it
is zero in the homotopy category of $\mathcal D(R)$.
Their mapping spectra into any fixed target therefore have zero
filtered colimit, so the fiber pro-object is zero,
proving~\ref{it:quotient-pro-objects}.
Applying $\delta(-)\otimes_{\delta R}^{\square}M$ preserves zero transition morphisms, and the resulting fiber tower has zero limit, by the Milnor sequence on profinite sections.
The inclusions $I^{r(n-1)+1}\subseteq I_n\subseteq I^n$ identify the resulting quotient limit with the $I$-adic quotient limit.
The comparison is $q_M$ because its composite with $\eta_M$ is induced by the quotient maps, proving~\ref{it:quotient-solid-limit}.
For noetherian $R$, the required pro-zero Koszul cohomology follows from Artin--Rees as in \cite[Tag~0921]{Stacks}.
\end{proof}

Under weak proregularity, Koszul completion in $\mathcal D(R)$ also computes the left derived functor of $I$-adic completion \cite[Corollary~5.25]{PSYTorsion}.

\subsection{Adic morphism objects}\label{subsec:adic-morphisms}

For a small $\Perf(R)$-linear category $\mathcal C$, define the $\DS$-enriched category $\Lambda_I\delta\mathcal C$ to have the objects of $\mathcal C$ and mapping objects $\Lambda_I\delta\mathcal C(X,Y)$.
Its units come from the completion unit, and its composition is the pairing~\eqref{eq:adic-completion-pairing} followed by completion of the original composition.

Take $\mathcal C=\Perf(R)$ and $\mathcal T=\thick(K_1)$.
Composing~\eqref{eq:adic-completion-pairing} with completion of the evaluation map defines
\[
 \Lambda_I\delta\Hom_R(P,Q)\otimes_{\delta R}^{\square}\Lambda_I\delta P
 \longrightarrow\Lambda_I\delta Q.
\]
The adjoint maps define an enriched functor $P\mapsto\Lambda_I\delta P$ from $\Lambda_I\delta\Perf(R)$ to solid $\delta R$-modules.
Under the Morita identification of right modules over $\delta\Perf(R)$ with solid $\delta R$-modules, enriched right Kan extension along $\delta\mathcal T\hookrightarrow\delta\Perf(R)$ carries $\delta r_{\mathcal T}P$ to $\Lambda_I\delta P$.
This right Kan extension is fully faithful by Proposition~\ref{thm:support-completion}, so it defines the comparison functor
\begin{equation}
  \Phi:\Lambda_I\delta\Perf(R)
        \longrightarrow\widehat{\Perf(R)}_{\mathcal T}^{\square},
  \qquad P\longmapsto P.
  \label{eq:adic-comparison-functor}
\end{equation}

Put
\begin{equation}
  A=\lim_n\delta(R/I^n),
  \label{eq:adic-solid-ring}
\end{equation}
with the multiplication induced by the quotient rings, and write $\widehat R_I=\lim_n R/I^n$.
For every profinite $S$, the transitions $C(S,R/I^{n+1})\to C(S,R/I^n)$ are surjective: a locally constant function has finitely many values, which can be lifted individually.
Thus $A$ is concentrated in degree zero and has sections $A(S)=C(S,\widehat R_I)$, where $\widehat R_I$ carries its inverse-limit topology and $C(S,\widehat R_I)$ denotes continuous functions.

When $I$ is weakly proregular, the lax monoidal equivalence $q_{\delta R}$ identifies $A$ with $\Lambda_I\delta R$ as a commutative $\delta R$-algebra.
For every solid $\delta R$-module $M$, we then define
\begin{equation}
  \alpha_M:A\otimes_{\delta R}^{\square}M
     \xrightarrow{\,1\otimes\eta_M\,}
     A\otimes_{\delta R}^{\square}\Lambda_I M
     \longrightarrow\Lambda_I M,
  \label{eq:adic-scalar-map}
\end{equation}
where the last arrow is the action supplied by~\eqref{eq:adic-completion-pairing}.
The lax monoidal structure makes $\alpha$ an $A$-linear natural transformation compatible with the tensor pairings.
Since $q$ respects those pairings, $q_M\circ\alpha_M$ is the map to the quotient limit induced by the projections $A\to\delta(R/I^n)$.
Applying $\alpha$ to mapping objects therefore defines an enriched functor
\[
  A\otimes_{\delta R}^{\square}\delta\Perf(R)
       \xrightarrow{\alpha}\Lambda_I\delta\Perf(R)
       \xrightarrow{\Phi}\widehat{\Perf(R)}_{\mathcal T}^{\square},
\]
whose source has mapping objects $A\otimes_{\delta R}^{\square}\delta\Hom_R(P,Q)$ and scalar-extended composition.

Let $\mathcal K_I$ be the idempotent-complete stable envelope of $\widehat{\Perf(R)}_{\mathcal T}$.
For weakly proregular $I$, the composite above induces $A\to\widehat{\Perf(R)}_{\mathcal T}^{\square}(R,R)$.
Point evaluation, together with Proposition~\ref{thm:solid-enrichment}, induces the corresponding map $\widehat R_I\to\Hom_{\mathcal K_I}(R,R)$.
Restriction along these algebra maps defines
\begin{equation}
\begin{aligned}
  \operatorname{ev}_R:\mathcal K_I
       &\longrightarrow\mathcal D(\widehat R_I),
       &U&\longmapsto\Hom_{\mathcal K_I}(R,U),\\
  \operatorname{ev}_R^{\square}:
  \RMod_{\widehat{\Perf(R)}_{\mathcal T}^{\square}}(\DS)
       &\longrightarrow\Mod_A(\DS),
       &N&\longmapsto N(R).
\end{aligned}
  \label{eq:adic-evaluation-functors}
\end{equation}
Precomposition defines the right actions, and commutativity of $A$ and $\widehat R_I$ identifies right and left modules.

The equivalence $\mathcal K_I\simeq\Perf(\widehat R_I)$ in part~\textup{(iii)} of Proposition~\ref{thm:adic-morphisms} below is proved in \cite[Propositions~B.2 and~B.4]{BGW}, based on \cite[Theorem~4.2]{PSYCompletion}.

\begin{proposition}
\label{thm:adic-morphisms}
Let $R$ be a commutative ring, let $I=(f_1,\ldots,f_r)$, and put $\mathcal T=\thick(K_1)\subseteq\Perf(R)$.
\begin{enumerate}
\item
The functor $\Phi$ in~\eqref{eq:adic-comparison-functor} is a $\DS$-enriched equivalence.
For $P,Q\in\Perf(R)$, $\Phi$ induces the equivalences
\begin{equation}
  \widehat{\Perf(R)}_{\mathcal T}^{\square}(P,Q)
  \xleftarrow[\ \sim\ ]{\Phi_{P,Q}}\Lambda_I\delta\Hom_R(P,Q),
  \label{eq:adic-morphism-koszul}
\end{equation}
compatibly with units and composition.
\item
If $I$ is weakly proregular, the enriched functor $\Phi\circ\alpha$ is also an equivalence, with maps
\begin{equation}
  \widehat{\Perf(R)}_{\mathcal T}^{\square}(P,Q)
  \xleftarrow[\ \sim\ ]{\Phi_{P,Q}\circ\alpha_{\delta\Hom_R(P,Q)}}
  A\otimes_{\delta R}^{\square}\delta\Hom_R(P,Q).
  \label{eq:adic-morphism-scalar}
\end{equation}
\item
If $I$ is weakly proregular, $\operatorname{ev}_R$ identifies $\mathcal K_I$ with $\Perf(\widehat R_I)$, and $\operatorname{ev}_R^{\square}$ is an equivalence of presentable $\DS$-module categories.
\end{enumerate}
\end{proposition}

\begin{proof}
Under the Morita identification used to define $\Phi$, the functor $\Lambda_{\mathcal T}$ of Proposition~\ref{thm:support-completion} is $\Lambda_I$.
Since $P$ is perfect, $\sHom_{\delta R}(\delta P,-)$ is tensor product over $\delta R$ with $\delta P^\vee$, where $P^\vee=\Hom_R(P,R)$, and therefore commutes with the Koszul limit.
Using the right Kan-extension realization and completion adjunction, we compute
\[
\begin{aligned}
 \widehat{\Perf(R)}_{\mathcal T}^{\square}(P,Q)
 &\simeq\sHom_{\delta R}(\Lambda_I\delta P,\Lambda_I\delta Q)\\
 &\simeq\sHom_{\delta R}(\delta P,\Lambda_I\delta Q)
 \simeq\Lambda_I\delta\Hom_R(P,Q).
\end{aligned}
\]
By the construction of $\Phi$ in~\eqref{eq:adic-comparison-functor}, $\Phi_{P,Q}$ is the inverse of this composite, using precomposition with the completion unit as in~\eqref{eq:completed-unit-precomposition}.
This proves~\eqref{eq:adic-morphism-koszul} for the specified functor.

Assume that $I$ is weakly proregular.
The complex $\Hom_R(P,Q)$ is perfect over $R$, so tensoring with its solid image commutes with the quotient limit.
Proposition~\ref{prop:quotient-towers} and the quotient-projection description of $\alpha$ identify $\alpha_{\delta\Hom_R(P,Q)}$ with this equivalence, proving~\eqref{eq:adic-morphism-scalar}.

Taking $P=Q=R$ in~\eqref{eq:adic-morphism-scalar} identifies the algebra map defining $\operatorname{ev}_R^\square$ as an equivalence; point evaluation carries it to the corresponding algebra equivalence for $\operatorname{ev}_R$.
Since restricted Yoneda and discrete realization are exact and $\Perf(R)=\thick(R)$, the image of $R$ thickly generates $\mathcal K_I$, and exactness of internal Hom implies that the enriched representables are thickly generated by the representable at $R$.
Every enriched module is a weighted colimit of representables, so the representable at $R$ generates the right-module category under colimits and the $\DS$-action.
By Morita theory over $\DZ$ and $\DS$, both evaluation functors are equivalences.
\end{proof}

\begin{remark}
\label{prop:completed-base-change}
Let $R\to R'$ be any homomorphism of commutative rings, and put $I'=IR'$ for a finitely generated ideal $I=(f_1,\ldots,f_r)\subseteq R$.
By~\eqref{eq:adic-completeness-test} and $R'\otimes_R R[1/f_i]=R'[1/f_i]$, restriction sends derived $I'$-complete solid $\delta R'$-modules to derived $I$-complete solid $\delta R$-modules.
Its left adjoint between the complete subcategories is $\Lambda_{I'}(\delta R'\otimes_{\delta R}^{\square}-)$.
The same adjunction shows that the map induced by $\eta_M$ is an equivalence
\begin{equation}
  \Lambda_{I'}(\delta R'\otimes_{\delta R}^{\square}M)
  \xrightarrow{\ \sim\ }
  \Lambda_{I'}(\delta R'\otimes_{\delta R}^{\square}\Lambda_I M)
  \label{eq:completed-base-change}
\end{equation}
for every solid $\delta R$-module $M$.
Successive completed extensions compose coherently by the same adjunction.
If $R'$ is perfect over $R$, tensoring with $\delta R'$ commutes with the Koszul limit, so $\delta R'\otimes_{\delta R}^{\square}\Lambda_I M$ is already $I'$-complete.
The same assertions hold in the derived module categories $\mathcal D(R)$ and $\mathcal D(R')$.
\end{remark}

\section{Solid scalar extension}
\label{sec:scalar-extension}

We compare solid scalar extension with derived completion using the maps $\alpha_M$ in~\eqref{eq:adic-scalar-map}, and apply this comparison to categorical completion.
For a commutative ring $R$ essentially of finite type over $\Z$ and an ideal $I\subseteq R$, we first prove that the solid tensor product over $\delta R$ of two bounded-above derived $I$-complete modules is derived $I$-complete.
Over any commutative noetherian ring $R$, the map $\alpha_{\delta V}$ in~\eqref{eq:adic-scalar-map} is an equivalence for every $V\in\mathcal D(R)$ whose cohomology is finitely generated over $R$ in every degree.
The categorical comparison uses enrichment over $\DS$ after forgetting the coefficient-algebra actions on the mapping objects.

\subsection{Tensor products of complete modules}

Cohomological bounds on solid modules refer to their underlying objects of $\DS$.

\begin{proposition}
\label{thm:solid-idempotence}
Let $R$ be a commutative ring essentially of finite type over $\Z$, let $I\subseteq R$ be an ideal, and put $A=\lim_n\delta(R/I^n)$.
Then:
\begin{enumerate}
\item If $M,N\in\Mod_{\delta R}(\DS)$ are bounded above and derived $I$-complete, then $M\otimes_{\delta R}^{\square}N$ is derived $I$-complete.
\item The multiplication map is an equivalence
\begin{equation}
  A\otimes_{\delta R}^{\square}A\xrightarrow{\ \sim\ }A.
  \label{eq:solid-completion-idempotence}
\end{equation}
Consequently restriction $\Mod_A(\DS)\to\Mod_{\delta R}(\DS)$ is fully faithful.
\end{enumerate}
\end{proposition}

\begin{proof}
Suppose first that $R$ is of finite type over $\Z$.
The Eilenberg--Mac Lane ring $\mathrm H R$ is noetherian and has $\pi_0(\mathrm H R)=R$.
Taking $\mathrm H\mathbb Z$-modules in solid spectra identifies solid $\mathrm H R$-modules and their relative tensor product with $\Mod_{\delta R}(\DS)$ and $\otimes_{\delta R}^{\square}$; the induced analytic structure agrees with the solid structure by \cite[Proposition~12.22 and Remark~12.23]{ScholzeAnalytic}.
This identification preserves the completeness test~\eqref{eq:adic-completeness-test}.
We may shift $M$ and $N$ independently to make them connective, since completeness is preserved by shifts.
Assertion~\textup{(i)} now follows from \cite[Lemma~1.26(3)]{CordovaFedeli}.

Write $R=S^{-1}R_0$, with $R_0$ of finite type over $\Z$, and choose a finitely generated ideal $I_0\subseteq R_0$ with $I=I_0R$.
The equivalence $R\otimes_{R_0}R\simeq R$ makes restriction along $\delta R_0\to\delta R$ fully faithful and identifies
\[
 M\otimes_{\delta R_0}^{\square}N
 \simeq M\otimes_{\delta R}^{\square}N.
\]
For each generator $f$ of $I_0$, the identity $R\otimes_{R_0}R_0[1/f]\simeq R[1/f]$ and extension--restriction adjunction identify the completeness tests in~\eqref{eq:adic-completeness-test}.
Applying the finite-type case over $R_0$ proves~\textup{(i)}.

The description following~\eqref{eq:adic-solid-ring} shows that $A$ is concentrated in degree zero, and Proposition~\ref{prop:quotient-towers} identifies $A$ with $\Lambda_I\delta R$.
Both $A$ and $A\otimes_{\delta R}^{\square}A$ are derived $I$-complete, the latter by~\textup{(i)}.
By Lemma~\ref{lem:adic-koszul-detection}, the unit $\eta_{\delta R}\otimes_{\delta R}^{\square}1_A:A\to A\otimes_{\delta R}^{\square}A$ becomes an equivalence after tensoring with $\delta K_1$.
Multiplication is its left inverse, so the same lemma proves~\eqref{eq:solid-completion-idempotence}.
For an $A$-module $N$, the extension--restriction counit is multiplication tensored over $A$ with $N$; hence it is an equivalence and restriction is fully faithful.
For finite-type $R$, the same idempotence also follows directly from \cite[Lemma~1.26(4)]{CordovaFedeli}.
\end{proof}

\begin{remark}
The conclusion of Proposition~\ref{thm:solid-idempotence}\textup{(i)} can fail when one factor is unbounded above.
Take $R=\Z$, $I=(p)$, and $A=\underline{\Z_p}$ with its $p$-adic topology.
The sums $\bigoplus_{n\geq0}A[n]$ and $\bigoplus_{n\geq0}A[-n]$ agree with their corresponding products: sums and products in solid abelian groups are exact, and at most one summand contributes in each cohomological degree.
Both sums are therefore complete, since products of complete modules are complete.
Using~\eqref{eq:solid-completion-idempotence}, we obtain
\[
 \left(\bigoplus_{n\geq0}A[n]\right)
 \otimes_{\Z_\square}
 \left(\bigoplus_{m\geq0}A[-m]\right)
 \simeq\bigoplus_{n,m\geq0}A[n-m].
\]
The terms with $n=m$ form the degree-zero direct summand $\bigoplus_{n\geq0}A$.
By Proposition~\ref{prop:quotient-towers}, point evaluation of its completion unit is
\[
 \begin{gathered}
 \bigoplus_{n\geq0}\Z_p\longrightarrow
       \lim_j\bigoplus_{n\geq0}\Z/p^j,\\
 \lim_j\bigoplus_{n\geq0}\Z/p^j
 =\left\{(a_n)\in\prod_{n\geq0}\Z_p:
       a_n\longrightarrow0\text{ in the }p\text{-adic topology}\right\}.
 \end{gathered}
\]
The sequence $(p^n)_{n\geq0}$ belongs to the target and not to the source.
Thus this summand, and hence the tensor product, is not complete.
\end{remark}

\subsection{Comparison with derived completion}

Let $R$ be a commutative noetherian ring, let $I\subseteq R$ be an ideal, and put $A=\lim_n\delta(R/I^n)$.
For a profinite set $S$, write $A(S)=\lim_n C(S,R/I^n)$.
For $M=A$,~\eqref{eq:adic-section-projection} specializes to
\begin{equation}
  A(S)\otimes_R V
  \xrightarrow{\sim}
  \bigl(A\otimes_{\delta R}^{\square}\delta V\bigr)(S),
  \qquad V\in\mathcal D(R).
  \label{eq:flat-profinite-coefficients}
\end{equation}

\begin{lemma}
\label{lem:flat-profinite-coefficients}
For noetherian $R$ and $A$ as above, $A(S)$ is flat over $R$ for every profinite set $S$.
Consequently, the functor $V\mapsto A\otimes_{\delta R}^{\square}\delta V$ from $\mathcal D(R)$ to $\Mod_{\delta R}(\DS)$ is $t$-exact.
\end{lemma}

\begin{proof}
Each $C(S,R/I^n)$ is a filtered colimit of finite free $R/I^n$-modules, indexed by finite clopen partitions of $S$.
The transition maps are surjective by the finite-value lifting argument following~\eqref{eq:adic-solid-ring}.
Thus $A(S)$ is flat over $R$ by \cite[Tag~0912]{Stacks}.
Since the free solid objects on profinite sets are compact projective generators, their section functors detect the $t$-structure.
Flatness and~\eqref{eq:flat-profinite-coefficients} therefore imply $t$-exactness.
\end{proof}

The following proposition proves a solid refinement of \cite[Proposition~4.12]{Shaul} for the map $\alpha$ defined in~\eqref{eq:adic-scalar-map}.

\begin{proposition}
\label{thm:finite-cohomology}
Let $R$ be a commutative noetherian ring, let $I\subseteq R$ be an ideal, and put $A=\lim_n\delta(R/I^n)$.
If $V\in\mathcal D(R)$ has finitely generated cohomology in every degree, then $\alpha_{\delta V}$ is an equivalence:
\begin{equation}
  A\otimes_{\delta R}^{\square}\delta V
  \xrightarrow[\ \sim\ ]{\alpha_{\delta V}}\Lambda_I\delta V.
  \label{eq:finite-cohomology-comparison}
\end{equation}
The maps $\alpha_{\delta V}$, defined for every $V\in\mathcal D(R)$, respect the tensor pairings induced by multiplication on $A$ and~\eqref{eq:adic-completion-pairing}.
\end{proposition}

\begin{proof}
Choose $r$ generators of $I$.
Each $K_n$ is a finite free complex in degrees $[-r,0]$.
The countable limit in~\eqref{eq:koszul-completion} is the fiber of a map between products.
Products are exact in solid abelian groups \cite[Theorems~1.10 and~5.8(i)]{Scholze}, and taking this fiber raises the upper cohomological bound by at most one.
Thus, for any solid $\delta R$-module $M$,
\begin{equation}
\begin{aligned}
  \Lambda_I M&\in\DS^{\geq a-r}
    &&\text{if }M\in\DS^{\geq a},\\
  \Lambda_I M&\in\DS^{\leq b+1}
    &&\text{if }M\in\DS^{\leq b}.
\end{aligned}
  \label{eq:completion-amplitude}
\end{equation}
Scalar extension is $t$-exact by Lemma~\ref{lem:flat-profinite-coefficients}.
Suppose that $V$ is bounded above.
By \cite[Tags~064Q and~066E]{Stacks}, it is pseudo-coherent, so there are perfect complexes $P_n\to V$ with cofibers in $\mathcal D(R)^{\leq-n}$.
Proposition~\ref{thm:adic-morphisms}\textup{(i)--(ii)}, applied to $R$ and $P_n$, shows that $\alpha_{\delta P_n}$ is an equivalence.
Scalar extension sends the cofiber of $P_n\to V$ to degrees at most $-n$, and completion sends it to degrees at most $-n+1$, by~\eqref{eq:completion-amplitude}.
Hence, in each fixed degree and for sufficiently large $n$, naturality identifies the comparison for $V$ with the comparison for $P_n$.

For arbitrary $V$, fix a degree $q$ and choose $b>q+r$.
The cofiber of $\tau_{\leq b}V\to V$ lies in degrees at least $b+1$.
By~\eqref{eq:completion-amplitude}, its completion lies in degrees at least $b+1-r>q$, while scalar extension preserves the original bound.
Thus this map induces isomorphisms in degree $q$ after both functors in~\eqref{eq:finite-cohomology-comparison}.
Since $\tau_{\leq b}V$ is bounded above with finitely generated cohomology in every degree, naturality and the preceding case prove the assertion in degree $q$.
Varying $q$ and using nondegeneracy proves the equivalence.
Compatibility with the tensor pairings follows from the construction of $\alpha$ after~\eqref{eq:adic-scalar-map}.
\end{proof}

For a small $\Perf(R)$-linear category $\mathcal C$ and $I=(f_1,\ldots,f_r)$, put
\[
  \mathcal T_I=\thick\{K_1\otimes_R X:X\in\mathcal C\},
  \qquad h_X=\delta\mathcal C(-,X).
\]
We apply $\Gamma_I$ and $\Lambda_I$ objectwise to enriched right $\delta\mathcal C$-modules, using the central $R$-action.
The category $A\otimes_{\delta R}^{\square}\delta\mathcal C$ has the objects of $\mathcal C$ and mapping objects $A\otimes_{\delta R}^{\square}\delta\mathcal C(X,Y)$, with units and composition induced by the unit and multiplication of $A$.
The maps $\alpha_{\delta\mathcal C(X,Y)}$ in~\eqref{eq:adic-scalar-map} define a functor $\alpha$ to $\Lambda_I\delta\mathcal C$ of \S~\ref{subsec:adic-morphisms}, acting as the identity on object labels.

Enriched Yoneda for $\Lambda_I\delta\mathcal C$, followed by restriction along the completion functor $\delta\mathcal T_I\to\Lambda_I\delta\mathcal C$, takes $X$ to $(\Lambda_Ih_X)|_{\mathcal T_I}$, with action induced by composition in $\Lambda_I\delta\mathcal C$.
The completion unit identifies this module with $\delta r_{\mathcal T_I}X$: at $K_1\otimes_R Y$, its component is $\delta K_1^\vee\otimes_{\delta R}^{\square}\eta_{h_X(Y)}$, an equivalence by $K_1^\vee\simeq K_1[-r]$ and Lemma~\ref{lem:adic-koszul-detection}; exactness and retracts extend this to $\mathcal T_I$.
The resulting functor lands in the full enriched category on the modules $\delta r_{\mathcal T_I}X$, defining $\Phi_{\mathcal C}$ in
\begin{equation}
  A\otimes_{\delta R}^{\square}\delta\mathcal C
     \xrightarrow{\alpha}\Lambda_I\delta\mathcal C
     \xrightarrow{\Phi_{\mathcal C}}
       \widehat{\mathcal C}_{\mathcal T_I}^{\square}.
  \label{eq:central-comparison-functors}
\end{equation}

\begin{proposition}
\label{thm:categorical-scalar-extension}
Let $R$ be a commutative noetherian ring, let $I=(f_1,\ldots,f_r)$, and let $\mathcal C$ be a small $\Perf(R)$-linear category.
Assume that $H^n\mathcal C(X,Y)$ is finitely generated over $R$ for every $X,Y$ and every $n\in\Z$.
For $\mathcal T_I$ and $A$ as above, the functors $\alpha$ and $\Phi_{\mathcal C}$ in~\eqref{eq:central-comparison-functors} are $\DS$-enriched equivalences.
On the original objects, their maps are
\begin{equation}
  \widehat{\mathcal C}_{\mathcal T_I}^{\square}(X,Y)
  \xleftarrow[\ \sim\ ]{(\Phi_{\mathcal C})_{X,Y}}\Lambda_I\delta\mathcal C(X,Y)
  \xleftarrow[\ \sim\ ]{\alpha_{\delta\mathcal C(X,Y)}}
  A\otimes_{\delta R}^{\square}\delta\mathcal C(X,Y).
  \label{eq:central-categorical-completion}
\end{equation}
These maps identify the units and composition maps of the three categories.
The exact extensions of $\alpha$ and $\Phi_{\mathcal C}$ to idempotent-complete stable envelopes, and restriction along $\alpha$ and $\Phi_{\mathcal C}$ on $\DS$-enriched right-module categories, are equivalences.
If $R$ is essentially of finite type over $\Z$, the map $\delta R\to A$ also satisfies~\eqref{eq:solid-completion-idempotence}.
\end{proposition}

\begin{proof}
The modules $h_{K_1\otimes_R X}\simeq\delta K_1\otimes_{\delta R}^{\square}h_X$ are fixed by $\Gamma_I$.
Every enriched right module is a weighted colimit of representables.
Applying $\Gamma_I$ and using that $K_\infty$ lies in the localizing subcategory generated by $K_1$ shows that these modules generate its image under colimits and the $\DS$-action.
By the generation statement following Proposition~\ref{thm:support-completion}, the same modules generate the image of $\Gamma_{\mathcal T_I}$, since $\mathcal T_I$ is the thick closure of the objects $K_1\otimes_R X$ and Yoneda preserves finite cofibers and retracts.
Thus the two functors are colocalizations onto the same full subcategory.
Taking right adjoints identifies $\Lambda_{\mathcal T_I}$ with objectwise $\Lambda_I$, compatibly with the completion units.

By~\eqref{eq:completed-kan-morphisms} and enriched Yoneda,
\[
 \widehat{\mathcal C}_{\mathcal T_I}^{\square}(X,Y)
 \simeq\sHom_{\delta\mathcal C}(h_X,\Lambda_Ih_Y)
 \simeq\Lambda_I\delta\mathcal C(X,Y).
\]
The first equivalence is right Kan extension followed by precomposition with $\eta_{h_X}$.
The composite of this chain with $(\Phi_{\mathcal C})_{X,Y}$ is the identity, since composition in $\Lambda_I\delta\mathcal C$ restricts along the completion unit to Yoneda evaluation.
Thus $(\Phi_{\mathcal C})_{X,Y}$ is the inverse of the displayed equivalence.
Proposition~\ref{thm:finite-cohomology} makes $\alpha$ an equivalence on all mapping objects.
Both functors preserve the original object labels, so they are enriched equivalences.
Their exact extensions to the idempotent-complete stable envelopes are therefore equivalences.
On right-module categories, enriched left Kan extension is inverse to restriction by~\eqref{eq:enriched-presheaf-universal}.
The last assertion is Proposition~\ref{thm:solid-idempotence}.
\end{proof}

The same mapping formula holds on the original objects of a small stable category equipped with an action of $\Perf(R)$ exact in each variable, under the same finiteness hypothesis on its mapping cohomology: apply the proposition to its idempotent completion.
Idempotent completion preserves the original mapping complexes, and the cohomology groups of the added mapping complexes are retracts of finitely generated $R$-modules.
The thick subcategory of $\Kar(\mathcal C)$ generated by $K_1\otimes_R X$ for $X\in\mathcal C$ is $\Kar(\mathcal T_I)$, so Proposition~\ref{prop:testing-envelope} identifies the mapping objects of the two categorical completions on objects of $\mathcal C$.

\begin{example}
\label{ex:discrete-padic-source}
The assumption that $R$ is essentially of finite type over $\Z$ in Proposition~\ref{thm:solid-idempotence} cannot be replaced by noetherianity.
Take the discrete ring $R=\Z_p$ and $A=\lim_n\delta(R/p^n)=\underline{\Z_p}$, where $A$ carries the $p$-adic topology.
By Scholze's multiplication equivalence $A\otimes_{\Z_{\square}}A\simeq A$ \cite[Example~6.4]{Scholze} and associativity,
\[
\begin{aligned}
  A\otimes_{\delta R}^{\square}A
  &\simeq\bigl(A\otimes_{\Z_{\square}}A\bigr)
       \otimes_{\delta(R\otimes_{\Z}R)}^{\square}\delta R\\
  &\simeq A\otimes_{\delta(R\otimes_{\Z}R)}^{\square}\delta R.
\end{aligned}
\]
The ring $R\otimes_{\Z}R$ acts on the two factors of $A\otimes_{\Z_\square}A$ separately; multiplication identifies this action with the action on $A$ through $R\otimes_{\Z}R\to R$.
Since $R=\Z_p$ is flat over $\Z$, the derived tensor product $R\otimes_{\Z}R$ is a discrete ring.
The projection formula~\eqref{eq:point-projection} over $R\otimes_{\Z}R$, with $\Gamma A=R$, identifies point evaluation of multiplication with the Hochschild augmentation $\HH_\bullet(\Z_p/\Z)=R\otimes_{R\otimes_{\Z}R}R\to R$.
The source has degree $-1$ cohomology
\[
  H^{-1}\!\left(R\otimes_{R\otimes_\Z R}R\right)
  \simeq\operatorname{Tor}^{R\otimes_\Z R}_1(R,R)
  \simeq\Omega^1_{R/\Z},
\]
where the last identification is the conormal description for the kernel of multiplication.
Localizing the module of differentials, we obtain $\Omega^1_{\Z_p/\Z}\otimes_{\Z_p}\mathbb Q_p=\Omega^1_{\mathbb Q_p/\mathbb Q}$, which is nonzero because $\mathbb Q_p/\mathbb Q$ has positive transcendence degree.
Thus multiplication on $\underline{\Z_p}$ is an equivalence over $\Z_\square$ but fails over $\delta\Z_p$.
\end{example}

\section{Perfect complexes on formal completions and punctures}\label{sec:formal-geometry}

Let $k=\Z$ or $\F_p$.
We identify the category of perfect complexes on the formal completion of a smooth proper $k$-scheme along any closed subscheme with the enriched Cauchy completion of its solid categorical completion.
Over $\F_p$, after point evaluation, the enriched Cauchy completion of the categorical puncture is equivalent to the category of all perfect complexes on the geometric puncture when the support is an effective Cartier divisor or a smooth closed subscheme.

\subsection{Tensor products of inverse systems of perfect complexes}

We first establish the tensor comparisons and pro-object identities needed for formal realization.

Let $X$ be a smooth proper $k$-scheme, let $\mathcal I\subset\mathcal O_X$ be a coherent ideal, set $Z=V(\mathcal I)$, and let
\[
 i_n:Z_n=V(\mathcal I^n)\hookrightarrow X.
\]
When regarded as an $\mathcal O_X$-algebra, $\mathcal O_{Z_n}$ denotes $i_{n*}\mathcal O_{Z_n}$.
The tensor product of presentable categories below extends the coefficients of $\mathcal D_{\mathrm{qc}}(X)$ along~\eqref{eq:discrete-adjunction}.
Write $\delta:\mathcal D_{\mathrm{qc}}(X)\to\mathcal E_X$ for external product with the coefficient unit $\delta k$.
Define the category $\mathcal E_X$ and a commutative algebra $\mathcal A$ in it by
\begin{equation}\label{eq:solid-geometric-category}
 \mathcal E_X=\mathcal D((k,\Z)_\square)\otimes_{\mathcal D(k)}\mathcal D_{\mathrm{qc}}(X),
 \qquad
 \mathcal A=\lim_n\delta(i_{n*}\mathcal O_{Z_n})\in\operatorname{CAlg}(\mathcal E_X).
\end{equation}
The functor $\delta$ is symmetric monoidal.
The tensor product $\otimes_{\mathcal E_X}$ combines tensor over $\mathcal O_X$ with solid tensor over $\delta k$, and the limit defining $\mathcal A$ is taken globally in $\mathcal E_X$.
Write $\sHom_{\mathcal E_X}$ and $\sHom_{\mathcal A}$ for the $\mathcal D((k,\Z)_\square)$-valued enriched mapping objects, and $\otimes_{\mathcal A}^{\square}$ for relative tensor product in $\Mod_{\mathcal A}(\mathcal E_X)$.
Define
\begin{equation}\label{eq:formal-perfect-category}
 \Perf(\widehat X_Z)=\lim_n\Perf(Z_n),
\end{equation}
using derived pullback as transition functors.
The tensor product $\otimes_{\mathcal O_{\widehat X_Z}}$ in this category is defined stagewise by $\otimes_{\mathcal O_{Z_n}}$.
For $L=(L_n)$ in this category, write
\begin{equation}\label{eq:formal-realization}
 \widetilde L=\lim_n\delta(i_{n*}L_n)\in\mathcal E_X.
\end{equation}
The stagewise actions equip $\widetilde L$ with an $\mathcal A$-module structure.

We use Tor-amplitude in the sense of \cite[Definition~15.68.1, Tag~0651]{Stacks}.

For countable inverse systems $(P_i)$ and $(Q_j)$ of perfect $k$-complexes, the limit projections and the monoidal structure of $\delta$ define the comparison
\begin{equation}\label{eq:properfect-tensor-map}
 \begin{aligned}
 \tau_{P_\bullet,Q_\bullet}:{}&
 \left(\lim_i\delta P_i\right)\otimes_{\delta k}^{\square}
 \left(\lim_j\delta Q_j\right)\\
 &\longrightarrow\lim_{i,j}\delta(P_i\otimes_k Q_j).
 \end{aligned}
\end{equation}

\begin{lemma}\label{lem:properfect-tensor}
Suppose there are fixed intervals $[a,b]$ and $[c,d]$ containing the Tor-amplitudes of every $P_i$ and every $Q_j$, respectively.
Then $\tau_{P_\bullet,Q_\bullet}$ in~\eqref{eq:properfect-tensor-map} is an equivalence in $\mathcal D((k,\Z)_\square)$.
\end{lemma}

\begin{proof}
Put $U=\colim_i P_i^\vee$ and $V=\colim_j Q_j^\vee$, where the duals on the right are taken in $\Perf(k)$.
These are bounded complexes with countable cohomology; since $\operatorname{gldim}(k)\leq1$, they admit bounded resolutions by countable free $k$-modules.
By internal tensor--Hom adjunction and perfect duality,
\[
 \lim_i\delta P_i\simeq\sHom_{\delta k}(\delta U,\delta k),
 \qquad
 \lim_j\delta Q_j\simeq\sHom_{\delta k}(\delta V,\delta k).
\]
For free modules, the required compatibility of duals with tensor products is the solid product formula
\begin{equation}\label{eq:solid-product-tensor}
 \left(\prod_I\delta k\right)\otimes_{\delta k}^{\square}
 \left(\prod_J\delta k\right)
 \simeq\prod_{I\times J}\delta k.
\end{equation}
For $k=\Z$ this is \cite[Proposition~6.3]{Scholze}.
The formula for $\F_p$ follows by extension of scalars, since reduction modulo $p$ commutes with products.
Finite totalization of the bounded free resolutions extends the formula to $U$ and $V$.
The result is
\[
 \sHom_{\delta k}(\delta U,\delta k)\otimes_{\delta k}^{\square}
 \sHom_{\delta k}(\delta V,\delta k)
 \simeq\sHom_{\delta k}\bigl(\delta(U\otimes_k V),\delta k\bigr),
\]
and the last term is the indicated inverse limit.
The map $\tau_{P_\bullet,Q_\bullet}$ is the canonical tensor--dual morphism, so these equivalences respect its associativity and symmetry constraints.
\end{proof}

For countable inverse systems $(P_i)$ and $(Q_j)$ in $\Perf(X)$, the same construction defines
\begin{equation}\label{eq:geometric-tensor-map}
 \begin{aligned}
 \tau^X_{P_\bullet,Q_\bullet}:{}&
 \left(\lim_i\delta P_i\right)\otimes_{\mathcal E_X}
 \left(\lim_j\delta Q_j\right)\\
 &\longrightarrow\lim_{i,j}\delta(P_i\otimes_{\mathcal O_X}Q_j).
 \end{aligned}
\end{equation}

\begin{lemma}\label{lem:geometric-tensor}
Assume that, for every $G\in\Perf(X)$, the complexes
\[
 \Hom_X(G,P_i),\qquad\Hom_X(G,Q_j)
\]
are perfect over $k$ with Tor-amplitudes bounded independently of the indices in their respective systems.
Then $\tau^X_{P_\bullet,Q_\bullet}$ in~\eqref{eq:geometric-tensor-map} is an equivalence in $\mathcal E_X$.
\end{lemma}

\begin{proof}
Choose a perfect compact generator $G$ of $\mathcal D_{\mathrm{qc}}(X)$, using \cite[Tag~09IS]{Stacks}.
After scalar extension, the $\mathcal D((k,\Z)_\square)$-valued functor $\sHom_{\mathcal E_X}(\delta G,-)$ detects equivalences in $\mathcal E_X$.
For the diagonal $\Delta:X\to X\times_k X$, write $P\boxtimes_k Q$ for the external tensor product.
For $P,Q\in\Perf(X)$, the projection formula identifies the global derived mapping complex as
\[
 \begin{aligned}
 &\Hom_X(G,P\otimes_{\mathcal O_X}Q)\\
 &\qquad\simeq\Hom_{X\times_k X}
 \bigl(\mathcal O_{X\times_k X},
       \Delta_*G^\vee\otimes_{\mathcal O_{X\times_k X}}(P\boxtimes_k Q)\bigr).
 \end{aligned}
\]
The kernel $\Delta_*G^\vee$ is perfect by smoothness.
Since $X$ is flat, quasi-compact and separated over $k$, external tensor products of perfect complexes generate $\Perf(X\times_k X)$ under finite colimits and retracts \cite[Lemma~8.11]{Toen}.

For $K_1,K_2\in\Perf(X)$, the K\"unneth equivalence extends after
solid scalar extension by separate colimit preservation and compatibility
with the coefficient action.
Moreover,
\[
 \sHom_{\mathcal E_X}\!\left(\delta K_1^\vee,\lim_i\delta P_i\right)
 \simeq\lim_i\delta\Hom_X(K_1^\vee,P_i),
\]
and similarly for $K_2$ and $(Q_j)$.
Thus K\"unneth identifies the comparison induced by the limit projections
for the kernel $K_1\boxtimes_kK_2$ with the map of
Lemma~\ref{lem:properfect-tensor} for
$\bigl(\Hom_X(K_1^\vee,P_i)\bigr)_i$ and
$\bigl(\Hom_X(K_2^\vee,Q_j)\bigr)_j$.
The comparison is therefore an equivalence for these kernels, hence for
$\Delta_*G^\vee$ by finite colimits and retracts.
Applying $\sHom_{\mathcal E_X}(\delta G,-)$ proves the assertion.
\end{proof}

An inverse system in the derived category of an abelian category, with cohomology in a fixed finite interval, is pro-zero whenever each of its cohomology systems is pro-zero \cite[Tag~0G3B]{Stacks}.
The same criterion holds in $\mathcal D_{\mathrm{qc}}(X)$ by applying it to complexes of $\mathcal O_X$-modules.

Let $L=(L_n)$ and $M=(M_n)$ belong to $\Perf(\widehat X_Z)$.
After passage to the cofinal diagonal in the two-index system, the stagewise balanced tensor products define a morphism of pro-objects
\begin{equation}\label{eq:formal-pro-product}
 \begin{aligned}
 \mu_{L,M}:{}&\{i_{n*}L_n\}_n\otimes_{\mathcal O_X}\{i_{m*}M_m\}_m\\
 &\longrightarrow\{i_{n*}(L_n\otimes_{\mathcal O_{Z_n}}M_n)\}_n.
 \end{aligned}
\end{equation}
For every fixed $r$, the quotient maps to $\mathcal O_{Z_r}$ and the compatible derived reductions of $L$ define
\begin{equation}\label{eq:formal-fixed-reduction}
 \operatorname{red}_{L,r}:
 \{i_{n*}L_n\otimes_{\mathcal O_X}\mathcal O_{Z_r}\}_{n\ge r}
 \longrightarrow\operatorname{const}(i_{r*}L_r).
\end{equation}

\begin{lemma}\label{lem:formal-pro-multiplication}
For each $L=(L_n)\in\Perf(\widehat X_Z)$, the complexes $i_{n*}L_n$ are perfect on $X$, with Tor-amplitudes contained in a finite interval independent of $n$.
The morphisms $\mu_{L,M}$ in~\eqref{eq:formal-pro-product} and $\operatorname{red}_{L,r}$ in~\eqref{eq:formal-fixed-reduction} are equivalences of pro-objects in $\mathcal D_{\mathrm{qc}}(X)$.
\end{lemma}

\begin{proof}
First take an affine open $\operatorname{Spec}R\subset X$ and choose generators $f_1,\ldots,f_s$ for the restricted ideal $I$.
The ideals $I_n=(f_1^n,\ldots,f_s^n)$ are cofinal with the powers of $I$.
Let $K_n$ be the Koszul complex on $f_1^n,\ldots,f_s^n$, in degrees $[-s,0]$.
By part~\ref{it:quotient-pro-objects} of Proposition~\ref{prop:quotient-towers}, the augmentation $\{K_n\}_n\to\{R/I_n\}_n$ is an equivalence of derived pro-objects.
Let $F_n$ be the fiber of $K_n\to R/I_n$.
For each $n$, choose $m\geq n$ such that $F_m\to F_n$ is zero in the derived category.
The transition
\[
 F_m\otimes_RR/I_m\longrightarrow F_n\otimes_RR/I_m
 \longrightarrow F_n\otimes_RR/I_n
\]
is then zero, so tensoring with the varying factor $R/I_n$ preserves the pro-zero fiber.

The differential in $K_n\otimes_R R/I_n$ vanishes.
Under the transition from $m$ to $n$, every summand indexed by a nonempty wedge of Koszul generators is multiplied by a nonempty product of factors $f_a^{m-n}$.
For $m\ge 2n$, this multiplier vanishes modulo $I_n$.
The augmentation therefore induces an equivalence of pro-objects in $\mathcal D(R)$
\[
 \{R/I_n\otimes_R R/I_n\}_n\xrightarrow{\ \sim\ }\{R/I_n\}_n.
\]
The same calculation with the second factor fixed at $R/I_r$ proves
\[
 \{R/I_n\otimes_R R/I_r\}_{n\ge r}
 \xrightarrow{\ \sim\ }\operatorname{const}(R/I_r).
\]
The diagonal is cofinal in the two-index system, and passage between $(I_n)$ and $(I^n)$ preserves these conclusions.

Apply \cite[Tag~09AW]{Stacks} to the $I$-adically complete ring $\widehat R_I$, using $\widehat R_I/I^n\widehat R_I\simeq R/I^n$.
Derived reduction induces the equivalence $\Perf(\widehat R_I)\simeq\lim_n\Perf(R/I^n)$; see also \cite[Proposition~A.4]{EfimovInfinity}.
Both comparisons are exact and preserve retracts in each formal perfect argument.
Since $\widehat R_I$ thickly generates, the two unit calculations prove that $\mu_{L,M}$ and $\operatorname{red}_{L,r}$ are equivalences on the affine open.

On $X$, compatibility with derived reduction implies
\[
 L_n\otimes_{\mathcal O_{Z_n}}\kappa(x)
 \simeq L_1\otimes_{\mathcal O_{Z_1}}\kappa(x)
 \qquad(x\in Z_1).
\]
Residue-field detection of Tor-amplitude \cite[Tag~068V]{Stacks} and quasi-compactness therefore imply that the Tor-amplitudes of all $L_n$ on $Z_n$ lie in a common interval $[a,b]$.
Since $X$ is regular of finite dimension, the $i_{n*}L_n$ are perfect on $X$ with Tor-amplitudes in $[a-\dim X,b]$.
The same argument establishes uniform bounds for $M$ and for the comparison fibers.
On a finite affine cover the fibers have pro-zero cohomology.
Choosing common later indices shows that the cohomology sheaves form pro-zero systems on $X$.
By \cite[Tag~0G3B]{Stacks}, these uniformly bounded fiber systems are pro-zero in $\mathcal D_{\mathrm{qc}}(X)$.
\end{proof}

\subsection{Perfect complexes on the formal completion}
\label{subsec:formal-recovery}

Let $\mathcal H$ be the stable idempotent-complete envelope of $\widehat{\Perf(X)}^{\square}_{\Perf_Z(X)}$ inside $\RMod_{\delta\Perf_Z(X)}\bigl(\mathcal D((k,\Z)_\square)\bigr)$, as in Definition~\ref{def:punctured-completion}.
For $P\in\Perf(X)$, we also write $P$ for its image under the completion functor.
We first describe the constructions on these objects, then extend them to $\mathcal H$ by exactness and preservation of retracts.
The category in~\eqref{eq:formal-perfect-category} has the solid mapping objects
\begin{equation}\label{eq:formal-solid-morphisms}
 \sHom_{\widehat X_Z}(L,M)
 =\lim_n\delta\Hom_{Z_n}(L_n,M_n),
\end{equation}
with composition induced by the stagewise compositions.
Evaluation at the point recovers its usual mapping complexes.
Derived pullback defines the restriction functor
\[
 \delta\Perf(X)\longrightarrow\Perf(\widehat X_Z),
 \qquad P\longmapsto(i_n^*P)_n.
\]
Its enriched maps are induced by the compatible coefficient extensions of derived pullback.

Define
\begin{equation}\label{eq:formal-realization-functor}
 \rho:\Perf(\widehat X_Z)\longrightarrow\Mod_{\mathcal A}(\mathcal E_X),
 \qquad L\longmapsto\widetilde L.
\end{equation}
Its enriched maps are adjoint to the limits of the stagewise evaluation maps
\[
 \begin{aligned}
 &\sHom_{\widehat X_Z}(L,M)\otimes_{\delta k}^{\square}\widetilde L\\
 &\quad\longrightarrow
   \delta\Hom_{Z_n}(L_n,M_n)\otimes_{\delta k}^{\square}\delta(i_{n*}L_n)
   \longrightarrow\delta(i_{n*}M_n).
 \end{aligned}
\]
These are $\mathcal A$-linear through $\mathcal A\to\delta\mathcal O_{Z_n}$, and their compatibility with transitions, composition and units follows from stagewise evaluation.
Likewise, the stagewise balanced tensor products define a lax symmetric monoidal structure
\[
 \widetilde L\otimes_{\mathcal A}^{\square}\widetilde M
 \longrightarrow\widetilde{L\otimes_{\mathcal O_{\widehat X_Z}}M},
 \qquad \rho(\mathcal O_{\widehat X_Z})=\mathcal A.
\]

To relate this realization to $\mathcal H$, first consider restricted enriched Yoneda:
\[
 \mathcal E_X\longrightarrow
 \RMod_{\delta\Perf(X)}\bigl(\mathcal D((k,\Z)_\square)\bigr),
 \qquad N\longmapsto\bigl(P\longmapsto\sHom_{\mathcal E_X}(\delta P,N)\bigr).
\]
Scalar extension of the $\mathcal D(k)$-linear module presentation of $\mathcal D_{\mathrm{qc}}(X)$ shows that this functor is an equivalence.
For $u:\Perf_Z(X)\hookrightarrow\Perf(X)$, use this equivalence and the Kan extensions of \S~\ref{subsec:full-subcategory} to define
\[
 \Gamma_Z=u_!u^*,\qquad \Lambda_Z=u_*u^*.
\]
Write $\eta_N:N\to\Lambda_ZN$ for the completion unit.

The equivalence between the Koszul and quotient pro-objects used in the proof of Lemma~\ref{lem:formal-pro-multiplication}, dualized on perfect complexes, identifies local cohomology with tensor product by $\colim_n\mathcal O_{Z_n}^{\vee}$.
Its image is generated by the tensor ideal $\Perf_Z(X)$ \cite[Tag~0A9A]{Stacks}, so scalar extension identifies it with $\Gamma_Z$.
Thus
\begin{equation}\label{eq:global-torsion-completion}
 \Gamma_Z(\delta\mathcal O_X)\simeq\colim_n\delta(\mathcal O_{Z_n}^{\vee}),
 \qquad
 \Lambda_Z N\simeq\lim_n(N\otimes_{\mathcal E_X}\delta\mathcal O_{Z_n}).
\end{equation}
The quotients are perfect on the regular finite-dimensional scheme $X$, and $(-)^\vee$ denotes their perfect duals.
The second formula is adjoint to $\Gamma_ZN\simeq\Gamma_Z(\delta\mathcal O_X)\otimes_{\mathcal E_X}N$; in particular, $\Lambda_Z(\delta\mathcal O_X)\simeq\mathcal A$.

For $N\in\mathcal E_X$, the projections $\mathcal A\to\delta\mathcal O_{Z_n}$ define a natural comparison
\begin{equation}\label{eq:global-completion-comparison}
 \alpha_N:\mathcal A\otimes_{\mathcal E_X}N
 \longrightarrow\Lambda_ZN
 \simeq\lim_n(N\otimes_{\mathcal E_X}\delta\mathcal O_{Z_n}).
\end{equation}
Dualizing the algebra structures on the perfect $\mathcal O_{Z_n}$ induces a counital coalgebra structure on $\Gamma_Z(\delta\mathcal O_X)$.
The counit identities and the idempotence of $\Gamma_Z$ identify it with the coalgebra of the smashing colocalization.
Its internal dual in $\mathcal E_X$, with multiplication induced by the coalgebra structure, is the inverse-limit algebra $\mathcal A$ of~\eqref{eq:solid-geometric-category}.
Thus $\alpha$ is the evaluation-induced morphism of monads $\mathcal A\otimes_{\mathcal E_X}-\to\Lambda_Z$, compatible with their units and multiplications.
If $P\in\Perf(X)$, tensoring with the dualizable object $\delta P$ commutes with limits, so this particular comparison is an equivalence:
\begin{equation}\label{eq:perfect-global-completion}
 \alpha_{\delta P}:\mathcal A\otimes_{\mathcal E_X}\delta P
 \xrightarrow{\ \sim\ }\Lambda_Z\delta P.
\end{equation}

For a derived $Z$-complete object $M\in\mathcal E_X$, the composite $\eta_M^{-1}\alpha_M$ is its canonical $\mathcal A$-action.
Enriched right Kan extension $u_*$, followed by restriction of monad
algebras along $\alpha$ and extension under finite cofibers and retracts,
defines
\begin{equation}\label{eq:geometric-generator-functor}
 J:\mathcal H\longrightarrow\Mod_{\mathcal A}(\mathcal E_X).
\end{equation}
For $P\in\Perf(X)$, \eqref{eq:perfect-global-completion} identifies
$J(P)\simeq\Lambda_Z\delta P\simeq
\mathcal A\otimes_{\mathcal E_X}\delta P$
as $\mathcal A$-modules, since $\alpha$ is a morphism of monads.

For a right $\mathcal H$-module $N$ and an $\mathcal A$-module $M$ in $\mathcal E_X$, define $J_!$ and $J^*$ by the coend and evaluation formulas
\begin{equation}\label{eq:geometric-coend-adjunction}
 \begin{aligned}
 J_!(N)&=\int^{P\in\mathcal H}N(P)\otimes_{\delta k}^{\square}J(P),\\
 J^*(M)(P)&=\sHom_{\mathcal A}(J(P),M).
 \end{aligned}
\end{equation}
The tensor is the coefficient action, and coend evaluation exhibits $J_!$ as the enriched left adjoint of $J^*$.
Define
\begin{equation}\label{eq:formal-cauchy-functor}
 \Phi=J^*\circ\rho:
 \Perf(\widehat X_Z)\longrightarrow
 \RMod_{\mathcal H}\bigl(\mathcal D((k,\Z)_\square)\bigr).
\end{equation}
For $P\in\Perf(X)$, the value of $\Phi(L)$ at $P$ is
\[
 \begin{aligned}
 \Phi(L)(P)
 &=\sHom_{\mathcal A}(\mathcal A\otimes_{\mathcal E_X}\delta P,\widetilde L)\\
 &\simeq\lim_n\delta\Hom_{Z_n}(i_n^*P,L_n),
 \end{aligned}
\]
by free-module adjunction, the limit formula for $\widetilde L$, and the adjunction for $i_n^*\dashv i_{n*}$.
The right $\mathcal H$-action is precomposition through $J$.

\begin{theorem}\label{thm:geometric-cauchy}
Let $X$ be smooth and proper over $k=\Z$ or $\F_p$, and let $Z\subset X$ be a closed subscheme.
For $\mathcal E_X$ and $\mathcal A$ defined in~\eqref{eq:solid-geometric-category}, multiplication is an equivalence
\[
 \mathcal A\otimes_{\mathcal E_X}\mathcal A\xrightarrow{\ \sim\ }\mathcal A.
\]
The coend functor in~\eqref{eq:geometric-coend-adjunction} is an equivalence of presentable $\mathcal D((k,\Z)_\square)$-linear categories
\begin{equation}\label{eq:geometric-free-envelope}
 J_!:\RMod_{\mathcal H}\bigl(\mathcal D((k,\Z)_\square)\bigr)
 \xrightarrow{\ \sim\ }\Mod_{\mathcal A}(\mathcal E_X).
\end{equation}
The functor $\Phi$ of~\eqref{eq:formal-cauchy-functor} takes values in Cauchy modules and induces a $\mathcal D((k,\Z)_\square)$-enriched equivalence with the solid mapping objects defined in~\eqref{eq:formal-solid-morphisms}
\begin{equation}\label{eq:geometric-cauchy-equivalence}
 \Phi:\Perf(\widehat X_Z)\xrightarrow{\ \sim\ }\Cauchy(\mathcal H).
\end{equation}
Under~\eqref{eq:geometric-free-envelope}, the Cauchy modules correspond exactly to the dualizable $\mathcal A$-modules.
For every $L=(L_n)\in\Perf(\widehat X_Z)$ and $r\geq1$, the canonical reduction map is an equivalence of $\delta\mathcal O_{Z_r}$-modules:
\begin{equation}\label{eq:formal-realization-reduction}
 \rho(L)\otimes_{\mathcal A}^{\square}\delta\mathcal O_{Z_r}
 \xrightarrow{\ \sim\ }\delta(i_{r*}L_r).
\end{equation}
The restriction of $P\in\Perf(X)$ corresponds to $\mathcal A\otimes_{\mathcal E_X}\delta P$.
Consequently, $\mathcal H$ corresponds to the subcategory of algebraizable perfect complexes in the sense of \cite[Introduction and Remark~0.3]{EfimovInfinity}:
\[
 \Perf_{\mathrm{alg}}(\widehat X_Z)
 =\thick\{P|_{\widehat X_Z}:P\in\Perf(X)\}.
\]
\end{theorem}

\begin{proof}
For $L=(L_n)$ and $M=(M_n)$ in $\Perf(\widehat X_Z)$, Lemma~\ref{lem:formal-pro-multiplication} implies that $i_{n*}L_n$ and $i_{n*}M_n$ have uniformly bounded Tor-amplitude.
For fixed $G\in\Perf(X)$, properness and smoothness make $\Hom_X(G,i_{n*}L_n)$ perfect over $k$ \cite[Tag~0DJT]{Stacks}.
The finite cohomological dimension of $X$ bounds their cohomology uniformly, and $\operatorname{gldim}(k)\leq1$ then bounds their Tor-amplitudes uniformly; the same argument applies to $M$.
The tensor--limit map $\tau^X$ in~\eqref{eq:geometric-tensor-map} for the systems $(i_{n*}L_n)_n$ and $(i_{n*}M_n)_n$, followed by the inverse-limit realization of $\mu_{L,M}$ in~\eqref{eq:formal-pro-product}, is therefore an equivalence by Lemmas~\ref{lem:geometric-tensor} and~\ref{lem:formal-pro-multiplication}.
Setting $L$ or $M$ equal to the formal structure sheaf, we obtain
\begin{equation}\label{eq:formal-realization-monoidal}
 \mathcal A\otimes_{\mathcal E_X}\mathcal A\simeq\mathcal A,
 \qquad
 \widetilde L\otimes_{\mathcal E_X}\widetilde M
 \simeq\widetilde{L\otimes_{\mathcal O_{\widehat X_Z}} M},
 \qquad
 \mathcal A\otimes_{\mathcal E_X}\widetilde L\simeq\widetilde L.
\end{equation}
Idempotence of $\mathcal A$ identifies $\otimes_{\mathcal A}^{\square}$ with $\otimes_{\mathcal E_X}$ on $\mathcal A$-modules.
The displayed maps come from stagewise multiplication and reduction, so they are the structure maps of $\rho$.
Thus $\rho$ is unital symmetric monoidal, and $\widetilde L$ is dualizable with dual $\widetilde{L^\vee}$.

We next prove that $J_!$ is an equivalence.
For $P,Q\in\Perf(X)$, right Kan extension and free-module adjunction identify the map on mapping objects induced by $J$ with
\[
 \begin{aligned}
 \mathcal H(P,Q)
 &\simeq\sHom_{\mathcal E_X}(\Lambda_Z\delta P,\Lambda_Z\delta Q)\\
 &\xrightarrow{\eta_{\delta P}^{*}}
   \sHom_{\mathcal E_X}(\delta P,\Lambda_Z\delta Q)\\
 &\simeq\sHom_{\mathcal A}(J(P),J(Q)).
 \end{aligned}
\]
The last identification uses~\eqref{eq:perfect-global-completion} for $Q$.
Compatibility of $\alpha$ with the monad units identifies the map induced by $J$ with the displayed precomposition map.
This map is an equivalence because $\Lambda_Z\delta Q$ is complete.
Exactness and preservation of retracts extend the comparison to $\mathcal H$, so $J$ is fully faithful.
For $P\in\Perf(X)$, free-module adjunction identifies
$\sHom_{\mathcal A}(J(P),-)$ with evaluation at $P$ after forgetting
the $\mathcal A$-action, so it preserves weighted colimits.
Finite cofibers and retracts extend this assertion to every object of
$\mathcal H$, while free-module bar resolutions show that the $J(P)$
generate under weighted colimits.
Since $J$ is fully faithful, the presheaf universal property~\eqref{eq:enriched-presheaf-universal} makes the unit of $J_!\dashv J^*$ an equivalence, and generation implies essential surjectivity of $J_!$, proving~\eqref{eq:geometric-free-envelope}.

The free module $\mathcal A=J(\mathcal O_X)$ is Cauchy.
For a dualizable module, mapping out of it amounts to mapping out of $\mathcal A$ after tensoring with its dual.
Hence all dualizable modules are Cauchy, including $\widetilde L$, and $J^*$ carries it to the Cauchy module $\Phi(L)$.
For $L,M\in\Perf(\widehat X_Z)$, the mapping objects satisfy
\begin{align}
 \sHom_{\mathcal A}(\widetilde L,\widetilde M)
 &\simeq\sHom_{\mathcal A}
       (\mathcal A,\widetilde{L^\vee\otimes_{\mathcal O_{\widehat X_Z}} M})\notag\\
 &\simeq\lim_n\delta\Hom_{Z_n}(L_n,M_n),
       \label{eq:formal-realization-mappings}
\end{align}
by free-module adjunction, perfect duality on each $Z_n$, and limit preservation by $\sHom_{\mathcal E_X}(\delta\mathcal O_X,-)$.
Under these identifications, the map on mapping objects induced by $\rho$ is an equivalence.
Indeed, when $L=\mathcal O_{\widehat X_Z}$, the construction by stagewise evaluation identifies that map with the identity of
$\lim_n\delta\Hom_{Z_n}(\mathcal O_{Z_n},M_n)$.
The strong monoidal structure of $\rho$ and duality reduce a general pair $(L,M)$ to the pair $(\mathcal O_{\widehat X_Z},L^\vee\otimes_{\mathcal O_{\widehat X_Z}}M)$.
Thus $\rho$, and hence $\Phi$, is fully faithful as an enriched functor.

For essential surjectivity, let $P$ be Cauchy in $\Mod_{\mathcal A}(\mathcal E_X)$.
For every $n$, extension along
$\mathcal A\to\delta(i_{n*}\mathcal O_{Z_n})$
preserves Cauchy objects.
Since the closed immersion $i_n$ is affine, there is a $\mathcal D(k)$-linear equivalence
\[
 \mathcal D_{\mathrm{qc}}(Z_n)
 \simeq\Mod_{i_{n*}\mathcal O_{Z_n}}\bigl(\mathcal D_{\mathrm{qc}}(X)\bigr).
\]
After extending coefficients, this equivalence becomes
\[
 \mathcal D((k,\Z)_\square)\otimes_{\mathcal D(k)}\mathcal D_{\mathrm{qc}}(Z_n)
 \xrightarrow{\ \sim\ }
 \Mod_{\delta(i_{n*}\mathcal O_{Z_n})}(\mathcal E_X),
\]
which sends the discrete realization of $L$ to $\delta(i_{n*}L)$ with its induced module structure.
Lemma~\ref{lem:discrete-cauchy} therefore supplies $L_n\in\Perf(Z_n)$ with
\begin{equation}\label{eq:cauchy-reductions}
 P\otimes_{\mathcal A}^{\square}\delta(i_{n*}\mathcal O_{Z_n})
 \simeq\delta(i_{n*}L_n).
\end{equation}
Because $\delta$ is symmetric monoidal, the discrete functors of Lemma~\ref{lem:discrete-cauchy} commute with derived extension of scalars.
By associativity of extension of scalars, the reductions of the left side of~\eqref{eq:cauchy-reductions} are compatible.
Full faithfulness of those discrete functors lifts these compatibilities coherently to the $L_n$, so $(L_n)$ belongs to $\Perf(\widehat X_Z)$.

Let $U:\Mod_{\mathcal A}(\mathcal E_X)\to\mathcal E_X$ be the forgetful functor.
To reconstruct $P$ from the reductions in~\eqref{eq:cauchy-reductions}, we prove that $U(P)$ is derived $Z$-complete, that is, the completion unit $U(P)\to\Lambda_ZU(P)$ is an equivalence.
By~\eqref{eq:absolute-yoneda} and~\eqref{eq:geometric-free-envelope}, $P$ is an absolute weighted colimit of the modules $J(Q)$, $Q\in\mathcal H$.
Their underlying objects are complete: this follows from~\eqref{eq:perfect-global-completion} for the objects coming from $\Perf(X)$ and then from closure of the complete subcategory under finite cofibers and retracts.
The functor $\Lambda_Z$ is enriched as the right adjoint of the coefficient-linear functor $\Gamma_Z$.
The unit is an enriched natural transformation, the mate of $\Gamma_Z\to\id$.
The enriched functors $U$ and $\Lambda_ZU$ preserve this absolute weighted colimit.
The completion unit is an equivalence on every $U(J(Q))$, so it is also an equivalence on $U(P)$.
Since $\mathcal A\otimes_{\mathcal E_X}\mathcal A\simeq\mathcal A$, associativity identifies $U(P)\otimes_{\mathcal E_X}\delta\mathcal O_{Z_n}$ with $U(P\otimes_{\mathcal A}^{\square}\delta\mathcal O_{Z_n})$.
The extension-of-scalars units now define a canonical $\mathcal A$-linear map
\[
 P\longrightarrow
 \lim_n\bigl(P\otimes_{\mathcal A}^{\square}\delta\mathcal O_{Z_n}\bigr)
 \simeq\widetilde L,
\]
where the limit is taken in $\Mod_{\mathcal A}(\mathcal E_X)$ and the last equivalence uses~\eqref{eq:cauchy-reductions} with its compatible module structures.
Its underlying map is the completion unit by~\eqref{eq:global-torsion-completion}, hence an equivalence; by conservativity of $U$, we conclude that $P\simeq\rho(L)$ as $\mathcal A$-modules.

Finally, we prove~\eqref{eq:formal-realization-reduction}.
The functor $-\otimes_{\mathcal E_X}\delta\mathcal O_{Z_r}$ commutes with limits, since $\mathcal O_{Z_r}$ is perfect on $X$.
The equivalence $\operatorname{red}_{L,r}$ in~\eqref{eq:formal-fixed-reduction} therefore shows that~\eqref{eq:formal-realization-reduction} is an equivalence of $\delta\mathcal O_{Z_r}$-modules.
Thus compatible discrete reduction is inverse to $\rho$, and $J^*\circ\rho$ is the asserted equivalence $\Phi$.
Every Cauchy object is therefore dualizable; the converse was proved above.
For $P\in\Perf(X)$, the projection formula identifies the natural map from $\mathcal A\otimes_{\mathcal E_X}\delta P$ to $\rho((i_n^*P)_n)$ with~\eqref{eq:perfect-global-completion}.
This identifies the restriction of $P$ with $J(P)$, and hence with the representable $h_P$ under $J^*$, proving the assertion about the thick subcategory they generate.
\end{proof}

Over $\F_p$, point evaluation of formal restriction recovers Efimov's equivalence with pseudo-perfect $\Perf_Z(X)$-modules \cite[Theorem~2.14]{EfimovInfinity} via
\[
 L\longmapsto\left(T\longmapsto\lim_n\Hom_{Z_n}(i_n^*T,L_n)\right).
\]

For a noetherian ring, Balmer--Sanders identify the dualizable derived-complete complexes for the completed tensor product with perfect complexes over the completion \cite[Theorem~1.4 and Corollary~5.2]{BalmerSanders}.
Efimov studies nuclear modules on affine formal schemes \cite{EfimovLimits}, and Matsukawa reconstructs formal schemes from their symmetric monoidal categories of nuclear modules \cite[abstract]{Matsukawa}.

\begin{remark}\label{rem:nonalgebraizable-formal-line}
Let $Z$ be the line at infinity in $X=\mathbb P^2_{\F_p}$.
In the coordinates of Efimov's Picard group calculation \cite[Proposition~8.1]{EfimovInfinity}, choose $h\in y^{-1}\F_p[[y^{-1}]]$ transcendental over $\F_p(y)$.
Such a series exists because the algebraic series form a countable set, whereas $y^{-1}\F_p[[y^{-1}]]$ is uncountable.
The class $1-x^{-1}h$ defines a line bundle on $\widehat X_Z$ whose punctured restriction is not algebraizable as a perfect complex by \cite[Theorem~8.6]{EfimovInfinity}.
Punctured restriction preserves finite cofibers and retracts and sends restrictions from $\Perf(X)$ into the algebraizable subcategory.
The bundle therefore lies outside the thick closure $\Perf_{\mathrm{alg}}(\widehat X_Z)$ of restrictions from $X$, and Theorem~\ref{thm:geometric-cauchy} places it in the enriched Cauchy completion.
\end{remark}

\subsection{The categorical puncture}
\label{subsec:geometric-puncture}

Theorem~\ref{thm:geometric-cauchy} identifies $\mathcal K=\Cauchy(\mathcal H)$ of Definition~\ref{def:punctured-completion} with the formal category enriched by~\eqref{eq:formal-solid-morphisms}.
Thus the canonical fully faithful enriched functor $\iota_{\mathrm{form}}:\mathcal Q_{\mathrm{form}}\to\mathcal Q$ of Proposition~\ref{prop:cauchy-quotient} is
\begin{equation}\label{eq:geometric-quotient-inclusion}
 \begin{aligned}
 \iota_{\mathrm{form}}:{}&
 \Kar\!\left(\Perf(\widehat X_Z)/\delta\Perf_Z(X)\right)\\
 &\longrightarrow\mathcal Q=\Cauchy\!\left(\mathcal H/\delta\Perf_Z(X)\right).
 \end{aligned}
\end{equation}

Suppose in addition that $X$ has pure relative dimension $d$ over $k=\Z$ or $\F_p$.
By relative Serre duality, $P\mapsto P\otimes_{\mathcal O_X}\omega_{X/k}[d]$ is the coherent Serre autoequivalence of $\Perf(X)$~\cite[Tag~0BRT]{Stacks}.
Tensoring by $\omega_{X/k}$ preserves support, so it induces autoequivalences of $\mathcal H$, of $\mathcal Q_0=\mathcal H/\delta\Perf_Z(X)$, and of $\mathcal Q=\Cauchy(\mathcal Q_0)$.
Under Theorem~\ref{thm:geometric-cauchy}, the induced action on $\Cauchy(\mathcal H)\simeq\Perf(\widehat X_Z)$ is tensor product with $\omega_{X/k}|_{\widehat X_Z}$; this action also descends to $\mathcal Q_{\mathrm{form}}$.
Applying Theorems~\ref{thm:punctured-solid-serre} and~\ref{thm:intro-punctured-serre}, we therefore obtain
\begin{equation}\label{eq:geometric-serre}
 \sqdual\mathcal Q(P,R)
 \simeq
 \mathcal Q(R,P\otimes_{\mathcal O_X}\omega_{X/k})[d-1].
\end{equation}
This is an internal solid duality on the categorical puncture for arbitrary closed support, over either coefficient ring.
It restricts along $\iota_{\mathrm{form}}$ in~\eqref{eq:geometric-quotient-inclusion}.
If the relative dimensions vary, the same formula applies on each open-and-closed pure-dimensional component.

We now identify the categorical puncture with perfect complexes on the punctured formal neighborhood in the following cases.
The notation $\widehat X_Z\setminus Z$ refers to the punctured formal neighborhood defined by completed affine charts and descent in~\cite{HPV}.

\begin{theorem}\label{thm:geometric-puncture}
Let $X$ be a smooth proper $\F_p$-scheme.
Suppose that $Z\subset X$ is either an effective Cartier divisor or a smooth closed subscheme.
Equip $\Perf(\widehat X_Z)$ with the solid mapping objects in~\eqref{eq:formal-solid-morphisms}, and form the quotients and Cauchy completion in~\eqref{eq:geometric-quotient-inclusion} over $\mathcal D((\F_p,\Z)_\square)$.
Then $\iota_{\mathrm{form}}$ is an enriched equivalence.
After point evaluation, we obtain an equivalence of $\F_p$-linear stable categories
\[
 \Gamma\mathcal Q\simeq\Perf(\widehat X_Z\setminus Z),
\]
including all perfect complexes on that punctured formal neighborhood.
The same identification holds for $\Gamma\mathcal Q_{\mathrm{form}}$ via $\Gamma\iota_{\mathrm{form}}$.
When $X$ has pure dimension $d$, the solid mapping objects on this puncture satisfy~\eqref{eq:geometric-serre}.
\end{theorem}

We use the module description in \S~\ref{subsec:punctured-module-model} and the quotient comparison in \S~\ref{subsec:classical-puncture-localization} to reduce the theorem to essential surjectivity of $\iota_{\mathrm{form}}$.
In \S~\ref{subsec:cartier-puncture}, the \v Cech descent argument proves this for an effective Cartier divisor.
In \S~\ref{subsec:smooth-center-puncture}, we blow up a smooth center and apply the equivalence established in \S~\ref{subsec:cartier-puncture} to the exceptional divisor.

\subsection{A module description of the categorical puncture}
\label{subsec:punctured-module-model}

Let $L_Z$ be the smashing localization away from $Z$ in $\mathcal E_X$, and put
\begin{equation}\label{eq:punctured-structure-algebra}
 \mathcal B
 =\mathcal A\otimes_{\mathcal E_X}L_Z(\delta\mathcal O_X),
 \qquad
 \mathcal M=\Mod_{\mathcal B}(\mathcal E_X).
\end{equation}
Both factors in $\mathcal B$ are idempotent commutative algebras, so $\mathcal B$ is idempotent as well.
Let $j:\delta\Perf_Z(X)\to\mathcal H$ be the fully faithful restriction of~\eqref{eq:support-completion-functor}, and let $U$ forget the $\mathcal A$-module structure.
The equivalence $J_!$ in~\eqref{eq:geometric-free-envelope} transports the localization of Proposition~\ref{prop:free-to-complete} to $\Lambda_ZU$:
\[
 \Lambda_ZUJ_!\simeq u_*j^*
 \quad\text{on }\RMod_{\mathcal H}\bigl(\mathcal D((k,\Z)_\square)\bigr).
\]
For $P\in\Perf(X)$, we use the enriched construction of~\eqref{eq:geometric-generator-functor} and the counit of $u^*\dashv u_*$ to obtain
$u^*UJ_!(h_P)\simeq u^*\Lambda_Z\delta P\simeq u^*\delta P\simeq j^*h_P$.
Finite cofibers, retracts, and the weighted-colimit universal property~\eqref{eq:enriched-presheaf-universal} extend this to $u^*UJ_!\simeq j^*$ on all modules; applying $u_*$ proves the displayed identity for the action on complete modules in~\eqref{eq:complete-solid-action}.

For $T\in\Perf_Z(X)$, the completion unit identifies $J(T)=\mathcal A\otimes_{\mathcal E_X}\delta T$ with $\delta T$.
Thus, for $N\in\Mod_{\mathcal A}(\mathcal E_X)$,
\[
 (j^*J^*N)(T)\simeq\sHom_{\mathcal A}(\delta T,N).
\]
In $\Mod_{\mathcal A}(\mathcal E_X)$, the modules $\delta T$, $T\in\Perf_Z(X)$, generate the kernel of $-\otimes_{\mathcal E_X}L_Z(\delta\mathcal O_X)$ under colimits and the $\mathcal D((k,\Z)_\square)$-action, by the support colocalization~\eqref{eq:global-torsion-completion} and free-module adjunction.
Their right orthogonal in $\Mod_{\mathcal A}(\mathcal E_X)$ is the essential image of this localization, namely $\Mod_{\mathcal B}(\mathcal E_X)$.
We therefore deduce from~\eqref{eq:punctured-module-kernel} an equivalence of presentable $\mathcal D((k,\Z)_\square)$-linear categories
\begin{equation}\label{eq:geometric-boundary-modules}
 \RMod_{\mathcal Q_0}\bigl(\mathcal D((k,\Z)_\square)\bigr)
 \simeq\Mod_{\mathcal B}(\mathcal E_X).
\end{equation}

The closed monoidal category $\mathcal M$ has two Hom constructions.
We write $[P,R]_{\mathcal B}$ for its internal Hom, which is a $\mathcal B$-module in $\mathcal E_X$.
The notation $\sHom_{\mathcal B}(P,R)$ continues to denote the enriched mapping object in $\mathcal D((k,\Z)_\square)$.

\begin{lemma}\label{lem:punctured-cauchy-dualizable}
Under~\eqref{eq:geometric-boundary-modules}, the objects of $\mathcal Q$ correspond exactly to the tensor-dualizable objects of $\mathcal M$.
\end{lemma}

\begin{proof}
The free modules $g_F=\mathcal B\otimes_{\mathcal E_X}\delta F$, $F\in\Perf(X)$, are dualizable Cauchy generators and include the unit: free-module adjunction proves the Cauchy assertion, and bar resolutions establish generation under weighted colimits, as in the proof of Theorem~\ref{thm:geometric-cauchy}.
Let $P$ be Cauchy.
For every $F\in\Perf(X)$,
\[
 \begin{aligned}
 \sHom_{\mathcal B}(g_F,[P,R]_{\mathcal B})
 &\simeq
 \sHom_{\mathcal B}(P\otimes_{\mathcal B}^{\square}g_F,R)\\
 &\simeq
 \sHom_{\mathcal B}(P,g_F^\vee\otimes_{\mathcal B}^{\square}R).
 \end{aligned}
\]
The right side preserves colimits and solid coefficient tensors.
Since the Cauchy generators $g_F$ jointly detect equivalences, $[P,-]_{\mathcal B}$ has the same preservation properties.
Put $P^\vee=[P,\mathcal B]_{\mathcal B}$.
The evaluation map
\[
 P^\vee\otimes_{\mathcal B}^{\square}R\longrightarrow[P,R]_{\mathcal B}
\]
is an equivalence for $R=g_F$ by duality of $g_F$, hence for every $R$ by generation.
Thus $P$ is dualizable.
Conversely, mapping from a dualizable $P$ is mapping from the Cauchy unit after tensoring with $P^\vee$, so $P$ is Cauchy.
\end{proof}

On a principal-adic chart, the following lemma identifies dualizable solid modules with perfect complexes over the ring of point sections.
For a topological ring $B$, write $\underline B$ for its associated condensed ring; the point evaluation $\underline B(*)$ is the underlying ring.

\begin{lemma}\label{lem:principal-adic-perfect}
Let $R$ be a noetherian $\F_p$-algebra, let $f\in R$ be a nonzerodivisor, and let $A=\widehat R_{(f)}$.
Equip $A$ with its $f$-adic topology and $B=A[1/f]$ with its complete Tate topology.
The associated condensed ring $\underline B$ is the solid algebra
\[
 \underline B=\left(\varprojlim_n\delta(R/f^n)\right)[1/f]
 \quad\text{in }\mathcal D((\F_p,\Z)_\square).
\]
Its dualizable solid modules form $\thick(\underline B)$.
Point evaluation identifies their $\F_p$-linear stable category with the perfect complexes over the ring $\underline B(*)$.
\end{lemma}

\begin{proof}
If $A=0$, the assertion is immediate.
The ring $A$ is a noetherian ring of definition, and $f$ is a topologically nilpotent unit of $B$.
Choose as $B^+$ the integral closure of $A$ in $B$.
Huber's sheaf theorem applies because the ring of definition is noetherian~\cite[Theorem~2.2]{HuberAdic}.
Thus $(B,B^+)$ is a complete sheafy Tate Huber pair, as required for the analytic results of~\cite{Andreychev}.
The surjective quotient tower identifies $\underline A\simeq\lim_n\delta(R/f^n)$: on extremally disconnected profinite sets, locally constant functions have finite image and lift along the quotient maps, so the derived limit is concentrated in degree zero.
A continuous map from a profinite set into $B$ has compact image and hence lands in some $f^{-m}A$.
Thus $\underline B\simeq\colim_{m\geq0}f^{-m}\underline A$, which proves the displayed solid description.

For a dualizable solid $\underline B$-module $P$, the internal Hom in $P\simeq[P^\vee,\underline B]_{\underline B}$ is also the condensed internal Hom, since the solid localization is closed monoidal and the target is solid.
Thus~\cite[Lemma~5.27]{Andreychev} places $P$ in the analytic category.
Analytic completion preserves this duality because it is symmetric monoidal and compatible with solid localization.
By~\cite[Corollary~5.51.1]{Andreychev}, $P$ therefore comes from a perfect $\underline B(*)$-complex, whose solid realization is built from $\underline B$ by finite cofibers and retracts; hence $P\in\thick(\underline B)$.
Conversely, every object of $\thick(\underline B)$ is dualizable in the solid module category.
For $P,Q\in\thick(\underline B)$, the object $\sHom_{\underline B}(P,Q)$ is obtained from $\underline B$ by finite stable operations and retracts in the two variables.
Point evaluation therefore identifies the mapping complexes:
\[
 \Gamma\sHom_{\underline B}(P,Q)
 \simeq\Hom_{\underline B(*)}(\Gamma P,\Gamma Q).
\]
This proves the assertion and its compatibility with scalar extension between such charts.
\end{proof}

The inverse limit topology in Lemma~\ref{lem:principal-adic-perfect} is part of the hypothesis.
On the completion at an ideal with several generators, the inverse limit topology need not be the adic topology defined by a chosen generator; the lemma therefore need not apply after that generator is inverted.
The smooth-center argument below instead reduces the full puncture to a Cartier divisor by a blowup.

\subsection{A quotient description of the puncture}
\label{subsec:classical-puncture-localization}

\begin{lemma}\label{lem:supported-formal-algebraization}
Let $X$ be a regular noetherian scheme of finite dimension, and let $Z\subset X$ be a closed subscheme.
Formal restriction is fully faithful on $\Perf_Z(X)$.
For each affine open $\operatorname{Spec}R\subset X$, write $I\subset R$ for the ideal defining $Z$.
Every bounded coherent complex on $\widehat X_Z$ whose cohomology on each completed affine chart $\operatorname{Spec}\widehat R_I$ is supported on $V(I\widehat R_I)$ is the formal restriction of an object of $\Perf_Z(X)$.
\end{lemma}

\begin{proof}
On an affine chart, flat completion preserves the bounded coherent mapping complex between supported perfect complexes because a power of the ideal annihilates its cohomology.
These comparisons satisfy descent, proving full faithfulness.
For a bounded coherent complex with the stated support, coherence implies that each cohomology module on $\operatorname{Spec}\widehat R_I$ is annihilated by a power of $I\widehat R_I$.
Boundedness and quasi-compactness ensure that a common power of the ideal annihilates all these cohomology sheaves, so they are coherent on a finite thickening of $Z$.
These sheaves algebraize, and regularity makes them perfect on $X$.
By full faithfulness on supported perfect complexes, we lift the finitely many extension maps in the Postnikov filtration and thereby algebraize the entire complex.
\end{proof}

For $T\in\Perf_Z(X)$ and $P\in\Perf(X)$, the restricted Yoneda formula~\eqref{eq:support-completion-functor} reads
\[
 \mathcal H(jT,P)\simeq\delta\Hom_X(T,P).
\]
Discreteness extends under finite cofibers and retracts to all objects $P\in\mathcal H$.
Theorem~\ref{thm:geometric-cauchy} identifies $h_{jT}\in\Cauchy(\mathcal H)$ with the formal restriction of $T$.
Applying Proposition~\ref{prop:cauchy-quotient}, we therefore obtain the $\mathcal D(k)$-enriched equivalence
\begin{equation}\label{eq:geometric-ordinary-formal-quotient}
 \Gamma\mathcal Q_{\mathrm{form}}
 \simeq\Kar\!\left((\Gamma\Perf(\widehat X_Z))/\Perf_Z(X)\right).
\end{equation}
The quotient on the right is formed over $\mathcal D(k)$, for either $k=\Z$ or $k=\F_p$.

\begin{lemma}\label{lem:classical-perfect-puncture}
Let $k$ be any field, let $X$ be a smooth finite-dimensional scheme of finite type over $k$, and let $Z\subset X$ be closed.
All categories and the quotient in this statement are $k$-linear stable $\infty$-categories enriched over $\mathcal D(k)$.
There is a canonical equivalence
\begin{equation}\label{eq:classical-perfect-puncture}
 \Kar\!\left(\Perf(\widehat X_Z)/\Perf_Z(X)\right)
 \xrightarrow{\ \sim\ }
 \Perf(\widehat X_Z\setminus Z).
\end{equation}
\end{lemma}

\begin{proof}
Write $\operatorname{Coh}^{-}$ for the category of almost perfect complexes, and let $\mathcal K_Z$ be the kernel of restriction from $\widehat X_Z$ to $\widehat X_Z\setminus Z$.
Since $X$ is quasi-compact and quasi-separated, Hennion--Porta--Vezzosi \cite[Theorem~7.3 and proof of Corollary~7.5]{HPV} identify this restriction as a $t$-exact localization:
\[
 \operatorname{Coh}^{-}(\widehat X_Z)/\mathcal K_Z
 \simeq\operatorname{Coh}^{-}(\widehat X_Z\setminus Z).
\]
The smooth affine charts are regular and excellent, so their completions are regular noetherian rings \cite[Tags~0AH2 and~0H7S]{Stacks}.
At each maximal ideal the completed local ring agrees with that of $X$ at the corresponding point of $Z$, so the dimensions are bounded by $\dim X$.
Thus the bounded coherent complexes on the completed charts and their punctures are precisely the perfect complexes.

A perfect complex on the puncture lifts to an almost perfect formal complex.
Truncating this lift below a lower bound for its restriction makes it bounded coherent and, by $t$-exactness, preserves its restriction.
For a quotient roof between bounded coherent complexes, choose a common lower cohomological bound $a$ for the endpoints and replace
\[
 A\longleftarrow C\longrightarrow B
 \qquad\text{by}\qquad
 A\longleftarrow\tau^{\geq a}C\longrightarrow B.
\]
Truncation adjunction factors both maps through $\tau^{\geq a}C$, which is bounded coherent.
By $t$-exactness, the new left map still becomes an equivalence after restriction.
The same procedure applies to refinements and every shift of the endpoints, proving full faithfulness on mapping spectra.
Lemma~\ref{lem:supported-formal-algebraization} identifies the bounded kernel with $\Perf_Z(X)$.
Passing to retracts proves~\eqref{eq:classical-perfect-puncture}.
\end{proof}

For smooth proper $X/\F_p$, the point evaluation of $\mathcal Q_{\mathrm{form}}$ is equivalent to $\Perf(\widehat X_Z\setminus Z)$ by~\eqref{eq:geometric-ordinary-formal-quotient} and~\eqref{eq:classical-perfect-puncture}.
To finish Theorem~\ref{thm:geometric-puncture}, we must show that $\iota_{\mathrm{form}}$ in~\eqref{eq:geometric-quotient-inclusion} is essentially surjective over $\mathcal D((\F_p,\Z)_\square)$ when $Z$ is an effective Cartier divisor or a smooth closed subscheme.

\subsection{\v Cech descent for a Cartier divisor}
\label{subsec:cartier-puncture}

Let $k=\F_p$ and suppose that $Z$ is an effective Cartier divisor.
Choose a finite affine open cover $(U_i)_{1\leq i\leq r}$ of $X$ on which $Z$ is principal.
For every nonempty increasing index set $J$, the intersection $U_J$ is affine because $X$ is separated.
Write $j_J:U_J\hookrightarrow X$, choose an equation $f_J$ for $Z|_{U_J}$, and put
\[
 R_J=\mathcal O(U_J),\qquad
 A_J=\varprojlim_n\delta(R_J/f_J^n),\qquad
 B_J=A_J[1/f_J].
\]
Identify $B_J$-modules in the scalar extension of $\mathcal D(R_J)$
with solid $B_J$-modules over $\mathcal D((\F_p,\Z)_\square)$,
with $\delta R_J$-action induced by $\delta R_J\to B_J$.
For $J\subset K$, restriction induces $A_J\to A_K$ and $B_J\to B_K$,
since the restricted defining ideals agree and their generators differ
by a unit; the transition functors are extension of scalars.

For each thickening, apply the finite alternating \v Cech resolution associated to this cover.
Scalar extension preserves that finite resolution.
The functors $j_{J*}$ in the scalar-extended categories are right adjoints and therefore preserve inverse limits.
Passing to the inverse limit over the thickenings, we obtain
\[
 \mathcal A\simeq\operatorname{Tot}_{J}j_{J*}A_J.
\]
Here and below $\operatorname{Tot}_{J}$ is the finite alternating totalization, with one term for each nonempty increasing index set.
We tensor this resolution with $L_Z(\delta\mathcal O_X)$ and apply the projection formula to obtain a finite resolution of $\mathcal B$ as a $\mathcal B$-module
\begin{equation}\label{eq:finite-cech-boundary}
 \mathcal B\simeq\operatorname{Tot}_{J}j_{J*}B_J.
\end{equation}
This construction uses the natural maps $j_J^*\mathcal A\to A_J$.
It does not require open pullback to commute with the inverse limit defining $\mathcal A$.

For $N\in\mathcal M$, define $F_JN$ by pullback to $U_J$ followed by extension of scalars along $j_J^*\mathcal B\to B_J$:
\begin{equation}\label{eq:completed-local-restriction}
 F_JN=B_J\otimes_{j_J^*\mathcal B}^{\square}j_J^*N.
\end{equation}
Under the affine identifications with solid modules, we use associativity of extension of scalars to obtain coherent equivalences
\[
 B_K\otimes_{B_J}^{\square}F_JN\simeq F_KN
 \qquad(J\subset K).
\]
The projection formula and finiteness of~\eqref{eq:finite-cech-boundary} imply
\[
 N\simeq\operatorname{Tot}_{J}j_{J*}F_JN.
\]
By extension--restriction adjunction, we then deduce an equivalence of solid mapping objects
\begin{equation}\label{eq:completed-cech-homs}
 \sHom_{\mathcal B}(M,N)
 \simeq
 \operatorname{Tot}_{J}
 \sHom_{B_J}(F_JM,F_JN).
\end{equation}
Thus the functor $N\mapsto(F_JN)_J$ into the \v Cech limit of the module categories is fully faithful, including on arbitrary modules.

We check its compatibility with formal realization before using it for essential surjectivity.
Let $L=(L_n)$ be any formal perfect complex, and set
\[
 P_J=A_J\otimes_{j_J^*\mathcal A}^{\square}j_J^*\widetilde L.
\]
Pullback along $j_J$ and extension of scalars to $A_J$ are symmetric monoidal, so $P_J$ is dualizable.
Applying~\eqref{eq:formal-realization-reduction} after pullback to $U_J$ and extension of scalars to $A_J$ induces compatible maps
\[
 P_J\longrightarrow
 P_J\otimes_{A_J}^{\square}\delta(R_J/f_J^n)
 \xrightarrow{\ \sim\ }\delta(L_n|_{U_J\cap Z_n}).
\]
Since $A_J=\lim_n\delta(R_J/f_J^n)$ as an $A_J$-module and tensoring
with $P_J$ preserves limits, these maps induce an equivalence
\begin{equation}\label{eq:completed-restriction-formal}
 P_J\simeq\varprojlim_n\delta(L_n|_{U_J\cap Z_n}).
\end{equation}
This applies to every formal perfect complex, including nonalgebraizable
ones; after inverting $f_J$, it identifies $F_J$ with classical punctured
restriction on $\mathcal Q_{\mathrm{form}}$.

Now let $P\in\mathcal M$ be dualizable.
Each $F_JP$ is dualizable, so Lemma~\ref{lem:principal-adic-perfect} places it in $\operatorname{thick}(B_J)$ and identifies it with a perfect $B_J(*)$-complex.
For $J\subset K$, after point evaluation, we have
\[
 \Gamma\!\left(B_K\otimes_{B_J}^{\square}F_JP\right)
 \simeq B_K(*)\otimes_{B_J(*)}\Gamma(F_JP).
\]
Indeed, the natural base-change map is an equivalence on $B_J$, and both functors preserve finite colimits and retracts.
The coherent transition maps therefore define an object of $\Perf(\widehat X_Z\setminus Z)$.
Lemma~\ref{lem:classical-perfect-puncture} realizes that object in the idempotent completion of the quotient of all formal perfect complexes.
Via~\eqref{eq:geometric-ordinary-formal-quotient}, choose the corresponding object $R$ of $\mathcal Q_{\mathrm{form}}$.
The compatibility~\eqref{eq:completed-restriction-formal} identifies, for every $J$, the image under $F_J$ of the $\mathcal B$-module corresponding to $\iota_{\mathrm{form}}(R)$ in~\eqref{eq:geometric-boundary-modules} with $F_JP$, compatibly with their descent data.
Full faithfulness~\eqref{eq:completed-cech-homs} lifts this identification to an equivalence with $P$.
By Lemma~\ref{lem:punctured-cauchy-dualizable}, this proves essential surjectivity in Theorem~\ref{thm:geometric-puncture} for Cartier divisors.

\subsection{Reduction of a smooth center to a Cartier divisor}
\label{subsec:smooth-center-puncture}

Let $k=\F_p$ and suppose that $Z$ is smooth.
Components on which $Z=X$ or $Z$ is empty have zero punctured category and can be treated separately.
On the remaining components, let
\[
 \pi:Y=\operatorname{Bl}_Z(X)\longrightarrow X,
 \qquad E=\pi^{-1}Z.
\]
The schemes $Y$ and $E$ are smooth, and $E$ is an effective Cartier divisor~\cite[Tag~0FUT]{Stacks}.
We first verify the direct-image calculation for every exceptional ideal power.
The natural map $\pi^*\mathcal I_Z^n\to\mathcal O_Y(-nE)$ has adjoint
\begin{equation}\label{eq:blowup-powers}
 b_n:\mathcal I_Z^n\longrightarrow \pi_*\mathcal O_Y(-nE).
\end{equation}
For every $n\geq0$, the map $b_n$ is an equivalence, compatible with the maps to $\mathcal O_X$.

Since $X$ and $Z$ are smooth over $k$, the closed immersion $Z\hookrightarrow X$ is regular~\cite[Tag~067U]{Stacks}.
\'Etale-locally on $X$, the pair $(X,Z)$ is pulled back from the coordinate pair $(\mathbb A_k^{m+c},\mathbb A_k^m)$~\cite[Tag~0FUE]{Stacks}.
Flat base change therefore reduces the calculation to the blowup of the zero section of $\mathbb A^c_{\mathbb A_k^m}$.
In the normal directions that blowup is the total space of $\mathcal O_{\mathbb P^{c-1}_{\mathbb A_k^m}}(-1)$.
Its affine projection to $\mathbb P^{c-1}_{\mathbb A_k^m}$ sends $\mathcal O_Y(-nE)$ to
\[
 \bigoplus_{\ell\geq0}\mathcal O_{\mathbb P^{c-1}_{\mathbb A_k^m}}(n+\ell).
\]
For $n\geq0$, all higher direct images of these summands to $\mathbb A_k^m$ vanish.
The direct image of the displayed sum is the ideal of polynomials of normal degree at least $n$, under the canonical map to the polynomial algebra in the normal coordinates.
This identifies $b_n$ in~\eqref{eq:blowup-powers} as an equivalence compatible with the ideal inclusions for every $n$, including $n=0$.
Flat base change transports this calculation from the local model to $X$.

For $n\geq1$, the induced morphism $\pi_n:nE\to Z_n$ defines a map of commutative algebras
\begin{equation}\label{eq:blowup-thickenings}
 \mathcal O_X/\mathcal I_Z^n\longrightarrow\pi_*\mathcal O_{nE}.
\end{equation}
Its underlying map is the map on cofibers induced by $b_n$ and $b_0$, so it is an equivalence by~\eqref{eq:blowup-powers}.
These algebra equivalences are compatible with the quotient transitions by construction.
Let
\[
 \mathcal E_Y=\mathcal D((k,\Z)_\square)
       \otimes_{\mathcal D(k)}\mathcal D_{\qc}(Y),\qquad
 \mathcal A_Y=\varprojlim_n\delta\mathcal O_{nE},\qquad
 \mathcal B_Y=\mathcal A_Y\otimes_{\mathcal E_Y}L_E(\delta\mathcal O_Y).
\]
The scalar-extended pushforward $\pi_*$ is a colimit-preserving right adjoint, agrees with $\delta\pi_*$ on discrete objects, and preserves limits.
Moreover $\pi^*L_Z(\delta\mathcal O_X)=L_E(\delta\mathcal O_Y)$.
From~\eqref{eq:blowup-thickenings} and the projection formula, we therefore deduce
\begin{equation}\label{eq:blowup-punctured-pushforward}
 \pi_*\mathcal A_Y\simeq\mathcal A,
 \qquad
 \pi_*\mathcal B_Y\simeq\mathcal B.
\end{equation}

Consider pullback along $\pi$ followed by extension of scalars to $\mathcal B_Y$:
\begin{equation}\label{eq:completed-blowup-pullback}
 F_\pi:\Mod_{\mathcal B}(\mathcal E_X)
 \longrightarrow\Mod_{\mathcal B_Y}(\mathcal E_Y),
 \qquad
 P\longmapsto
 \mathcal B_Y\otimes_{\pi^*\mathcal B}^{\square}\pi^*P.
\end{equation}
The right adjoint of $F_\pi$ is $\pi_*$, with the $\mathcal B$-action induced by $\mathcal B\to\pi_*\mathcal B_Y$.
The relative projection formula identifies the unit $P\to\pi_*F_\pi P$ with
\[
 P\longrightarrow P\otimes_{\mathcal B}^{\square}\pi_*\mathcal B_Y.
\]
The formula applies to arbitrary modules: it holds for the dualizable perfect generators and their solid coefficient tensors, and extends by colimits and free-module resolutions.
By~\eqref{eq:blowup-punctured-pushforward}, the unit is an equivalence.
Thus $F_\pi$ is fully faithful.
It is also symmetric monoidal, and hence preserves dualizable modules.

We again check compatibility with the entire formal category.
For $L=(L_n)\in\Perf(\widehat X_Z)$, pull back its realization to $Y$ and extend scalars to $\mathcal A_Y$:
\[
 P_L=\mathcal A_Y\otimes_{\pi^*\mathcal A}^{\square}\pi^*\widetilde L.
\]
Pullback and extension of scalars are symmetric monoidal, so $P_L$ is dualizable over $\mathcal A_Y$.
Writing $e_n:nE\hookrightarrow Y$, we apply~\eqref{eq:formal-realization-reduction} to obtain
\begin{equation}\label{eq:blowup-formal-reductions}
 P_L\otimes_{\mathcal A_Y}^{\square}\delta\mathcal O_{nE}
 \simeq
 \delta\mathcal O_{nE}
 \otimes_{\pi^*\delta\mathcal O_{Z_n}}^{\square}\pi^*\delta(i_{n*}L_n)
 \simeq\delta(e_{n*}\pi_n^*L_n).
\end{equation}
As in the argument for~\eqref{eq:completed-restriction-formal}, dualizability and the compatible reductions~\eqref{eq:blowup-formal-reductions} imply
\[
 P_L\simeq\lim_n\delta(e_{n*}\pi_n^*L_n).
\]
Thus realization identifies formal pullback of every formal perfect complex with pullback followed by extension of scalars to $\mathcal A_Y$.
Localization, passage to the quotient, and retracts preserve this
compatibility.

Let $P$ be a dualizable $\mathcal B$-module, and write $\mathcal Q_Y$ for the categorical puncture associated to $(Y,E)$.
Using~\eqref{eq:geometric-boundary-modules} and Lemma~\ref{lem:punctured-cauchy-dualizable} for both $(X,Z)$ and $(Y,E)$, we regard $F_\pi$ as an enriched functor $\mathcal Q\to\mathcal Q_Y$.
The equivalence established in \S~\ref{subsec:cartier-puncture} for the effective Cartier divisor $E$ identifies $\Gamma\mathcal Q_Y$ with $\Perf(\widehat Y_E\setminus E)$.
Using this identification together with~\eqref{eq:geometric-ordinary-formal-quotient} and Lemma~\ref{lem:classical-perfect-puncture}, we express the compatibility with formal pullback proved above by the following commutative square
\[
 \begin{CD}
 \Gamma\mathcal Q_{\mathrm{form}}
     @>{\Gamma(F_\pi\circ\iota_{\mathrm{form}})}>>
       \Gamma\mathcal Q_Y\\
 @V{\simeq}VV @VV{\simeq}V\\
 \Perf(\widehat X_Z\setminus Z)
     @>{\pi^*}>{\simeq}>
       \Perf(\widehat Y_E\setminus E).
 \end{CD}
\]
The bottom arrow is an equivalence by~\cite[Corollary~2.16(1)]{EfimovInfinity}: both $X$ and $Y$ are smooth and proper over $\F_p$, and $\pi$ restricts to an isomorphism $Y\setminus E\to X\setminus Z$.
Let $L_Y\in\Perf(\widehat Y_E\setminus E)$ correspond to $F_\pi P$ under the right vertical equivalence.
Choose a perfect complex $L_X\in\Perf(\widehat X_Z\setminus Z)$ and an equivalence $\pi^*L_X\simeq L_Y$.
The left vertical equivalence supplies $R\in\mathcal Q_{\mathrm{form}}$ whose image is $L_X$.
Commutativity of the square then implies $F_\pi\iota_{\mathrm{form}}(R)\simeq F_\pi P$.
Full faithfulness of $F_\pi$ implies $\iota_{\mathrm{form}}(R)\simeq P$.
This proves essential surjectivity for smooth closed centers and completes the proof of Theorem~\ref{thm:geometric-puncture}.

\section{Hochschild cochain comparisons}\label{sec:hochschild}

For a full subcategory of modules containing the representables, restriction identifies its Hochschild cochains with those of the source category as $E_2$-algebras (Proposition~\ref{prop:restriction-centers}).
For adic algebras, scalar extension respects cup products and compares solid Hochschild cochains with the solid derived completion of the cochains of the original ring (Proposition~\ref{thm:hochschild-comparison}).

\subsection{Restriction to representables}

Let $\mathcal V$ be a closed presentably symmetric monoidal stable category whose tensor preserves colimits in each variable.
For a small $\mathcal V$-enriched category $\mathcal A$ put $\mathcal A^e=\mathcal A^{\op}\otimes_{\mathcal V}\mathcal A$ and $\Delta_{\mathcal A}(a,b)=\mathcal A(a,b)$.
The diagonal is the kernel of the identity functor of $\RMod_{\mathcal A}(\mathcal V)$, and the presheaf universal property determines a colimit-preserving $\mathcal V$-linear endofunctor and its transformations from its kernel.
Thus the Hochschild cochain object $\HH^\bullet_{\mathcal V}(\mathcal A)=\sHom_{\mathcal A^e}(\Delta_{\mathcal A},\Delta_{\mathcal A})$, with internal Hom taken between left $\mathcal A^e$-modules, is the internal endomorphism object of this identity functor.
This description equips the Hochschild cochain object with its usual center $E_2$-structure \cite[Proposition~2.1.5]{BR}.

\begin{proposition}\label{prop:restriction-centers}
Let $\mathcal V$ be closed, presentably symmetric monoidal and stable.
Let $\mathcal A$ be a small $\mathcal V$-enriched category and $\mathcal H$ a small full $\mathcal V$-enriched subcategory of $\RMod_{\mathcal A}(\mathcal V)$ containing the representables.
Restriction along $j:\mathcal A\to\mathcal H$, $a\mapsto h_a=\mathcal A(-,a)$, is an equivalence of $E_2$-algebras in $\mathcal V$,
\[
 j^*:\HH^\bullet_{\mathcal V}(\mathcal H)\xrightarrow{\sim}\HH^\bullet_{\mathcal V}(\mathcal A).
\]
\end{proposition}

\begin{proof}
Restrict the diagonal first along $j^{\op}\otimes\id$.
The result is the evaluation bimodule
\[
 U:\mathcal A^{\op}\otimes_{\mathcal V}\mathcal H\longrightarrow\mathcal V,
 \qquad U(a,M)=\mathcal H(h_a,M)\simeq M(a).
\]
Restricting once more along $\id\otimes j$ recovers $\Delta_{\mathcal A}$.
These restrictions induce maps in $\mathcal V$
\begin{equation}\label{eq:center-evaluation-restrictions}
 \HH^\bullet_{\mathcal V}(\mathcal H)
 \xrightarrow{(j^{\op}\otimes\id)^*}\operatorname{End}(U)
 \xrightarrow{(\id\otimes j)^*}\HH^\bullet_{\mathcal V}(\mathcal A)
\end{equation}
whose composite is $j^*$.
Here $\operatorname{End}(U)$ is the $\mathcal V$-valued endomorphism object in $\Fun_{\mathcal V}(\mathcal A^{\op}\otimes_{\mathcal V}\mathcal H,\mathcal V)$.
Multiplication in all three terms is composition of internal natural transformations.

Currying in the variable $M$ identifies $M\mapsto U(-,M)$ with the inclusion of $\mathcal H$ into right $\mathcal A$-modules.
In this description, $(j^{\op}\otimes\id)^*$ in~\eqref{eq:center-evaluation-restrictions} applies restriction to maps between the right representables of $\mathcal H$.
At $M,N$ the resulting comparison is
\[
 \mathcal H(M,N)\longrightarrow\sHom_{\mathcal A}(M,N).
\]
It is an equivalence because $\mathcal H$ is full on the $\mathcal V$-valued mapping objects.
Taking the end for natural transformations proves that $(j^{\op}\otimes\id)^*$ is an equivalence.

In the other variable $a\mapsto U(a,-)$ takes values in covariant representables $\mathcal H(h_a,-)$ and is contravariant in $a$.
Restriction to $\mathcal A$ sends them to $\mathcal A(a,-)$.
For $a,b\in\mathcal A$, covariant enriched Yoneda identifies the map between their internal Homs with
\[
 \mathcal H(h_b,h_a)\xrightarrow{\sim}\mathcal A(b,a),
\]
the inverse of the Yoneda embedding on this Hom.
The order $b,a$ comes from the variance of covariant representables.
Taking the end in the $\mathcal A^{\op}$ variable proves that $(\id\otimes j)^*$ in~\eqref{eq:center-evaluation-restrictions} is an equivalence.
Both arrows are restrictions of transformations and preserve composition.

For the $E_2$-structure, consider the $\mathcal V$-linear colimit-preserving extension $j_!:\RMod_{\mathcal A}(\mathcal V)\to\RMod_{\mathcal H}(\mathcal V)$ as an arrow between presentable $\mathcal V$-linear categories.
Its relative center $\HH^\bullet_{\mathcal V}(j_!)$ is the internal endomorphism object of the unit in the monoidal endomorphism category of this arrow.
It has maps of $E_2$-algebras to the two centers, with underlying pullback
\begin{equation}\label{eq:relative-center-pullback}
 \HH^\bullet_{\mathcal V}(j_!)\simeq
 \HH^\bullet_{\mathcal V}(\mathcal A)\mathbin{\times}_{\operatorname{End}(j_!)}\HH^\bullet_{\mathcal V}(\mathcal H)
\end{equation}
in $\mathcal V$ \cite[Proposition~2.2.4]{BR}.
Here $\operatorname{End}(j_!)$ is the internal object of natural endomorphisms of $j_!$ in the category of colimit-preserving $\mathcal V$-linear functors.
The maps to this object are postcomposition $\eta\mapsto j_!\eta$ from the $\mathcal A$ center and precomposition $\theta\mapsto\theta j_!$ from the $\mathcal H$ center.
The unit of $j_!\dashv j^*$ is invertible on representables by fullness of $j$, and on modules by weighted colimits.
Thus postcomposition with $j_!$ identifies the source center with $\operatorname{End}(j_!)$.
The projection from the pullback to the $\mathcal H$ center is an equivalence.
Its inverse followed by projection to the $\mathcal A$ center takes a transformation to its action on $j_!h_a$.
Since $j_!h_a$ is represented by $ja$, evaluation on representables identifies this map with restriction along $j^{\op}\otimes j$, the composite already proved invertible.
\end{proof}

For dg categories over a commutative ring, Lowen--Van den Bergh obtain restriction to the representables for full subcategories of cofibrant modules, with the brace operations \cite[\S~2.4 and Theorems~4.1.2 and~4.4.1]{LowenVanDenBergh}.
For a proper dg category over a field and full intermediate subcategories of cohomologically finite modules (with finite-dimensional total cohomology at each object), see Goodbody \cite[Theorem~5.4 and Remark~5.5]{Goodbody}.

Apply Proposition~\ref{prop:restriction-centers} to $\mathcal A=\delta\mathcal T$ and the full subcategory $\mathcal H$ of \S~\ref{subsec:punctured-completion}, which contains the representables.
Write $\HH^\bullet_\square$ for the cochain object over $\mathcal D((k,\Z)_\square)$.
The resulting equivalence is
\begin{equation}\label{eq:completion-center-restriction}
 j^*:\HH^\bullet_\square(\mathcal H)\xrightarrow{\sim}\HH^\bullet_\square(\delta\mathcal T).
\end{equation}
The target of $j^*$ in~\eqref{eq:completion-center-restriction} is computed by the internal end for $\delta\mathcal T$ in $\mathcal D((k,\Z)_\square)$.
Comparison with the discrete realization of the $\mathcal D(k)$-valued Hochschild cochains requires compatibility with the end defining those cochains.
For a commutative ring $R$ essentially of finite type over $\Z$ and an ideal $I\subseteq R$, the next subsection compares $\HH^\bullet_\square(A/\Z)$, where $A=\lim_n\delta(R/I^n)$, with $\Lambda_I\delta\HH^\bullet(R/\Z)$ using scalar extension and finite approximations to the diagonal.

\subsection{Adic Hochschild cochains}\label{subsec:adic-hochschild}

In this subsection, we compare complexes in $\DZ$ with solid complexes in $\DS$, and all Hochschild objects are relative to $\Z$.
For flat noetherian algebras, Shaul compares derived completions of Hochschild cochains under weak proregularity of the induced ideal of the enveloping algebra \cite[Theorem~4.1]{Shaul}.
With noetherian tensor square and finitely generated bimodule coefficients he also identifies the completed Ext groups \cite[Theorem~4.13]{Shaul}.
Proposition~\ref{thm:finite-cohomology} establishes a solid refinement of Shaul's comparison for unbounded complexes with finitely generated cohomology in every degree \cite[Proposition~4.12]{Shaul}.
The solid comparison uses the maps $\alpha_M$ in~\eqref{eq:adic-scalar-map} and an internal Hom from the derived diagonal.

For a commutative ring $R$ viewed in degree zero and a commutative algebra $A$ in $\mathcal D(\Z_\square)$, write $R^e=R\otimes_\Z R$ and $A^e=A\otimes_{\Z_\square}A$ for their derived enveloping algebras.
Composition of endomorphisms of the diagonals defines the cup products on the Hochschild cochains
\[
 \HH^\bullet(R/\Z)=\Hom_{R^e}(R,R),\qquad
 \HH^\bullet_\square(A/\Z)=\sHom_{A^e}(A,A).
\]

Over a connective algebra $D$ in $\mathcal D(\Z)$, finite-cell modules are built from zero by finitely many attachments of shifts of $D$.
A bounded-above module admitting finite-cell approximations with cofibers in arbitrarily negative degrees is \emph{pseudo-coherent}, or \emph{almost perfect} \cite[Definition~7.2.4.10 and Proposition~7.2.4.11(5)]{LurieHA}.

\begin{lemma}\label{lem:finite-diagonal-approximations}
Let $R$ be a commutative ring essentially of finite type over $\Z$.
For $D=R\otimes_\Z R$ and each $n\geq1$, there is a finite-cell $D$-module $E_n\to R$ with cofiber in $\mathcal D(D)^{\leq-n}$.
\end{lemma}

\begin{proof}
Write $R=S^{-1}R_0$ with $R_0$ of finite type over $\Z$.
Since $R_0$ is discrete and of finite type over $\Z$, it is noetherian and almost finitely presented over $\Z$ \cite[Remark~7.2.4.28 and Proposition~7.2.4.31]{LurieHA}.
Base change makes $D_0=R_0\otimes_\Z R_0$ almost finitely presented over $R_0$, so $D_0$ is connective and noetherian by the same proposition.
The discrete diagonal $R_0$ is cyclic over $H^0(D_0)$, hence almost perfect \cite[Proposition~7.2.4.17]{LurieHA}.
The proof of \cite[Proposition~7.2.4.11(5)]{LurieHA} supplies the finite-cell approximations.
Localizing both factors at $S$ preserves finite cells and is $t$-exact, so the cofiber bounds persist for the localized diagonal
\[
 D\otimes_{D_0}R_0\simeq R\otimes_{R_0}R\simeq R.
\]
\end{proof}

To apply these approximations to cochains, we compare scalar extension of the derived mapping complex with internal Hom into the completed target.
Let $D$ be a connective commutative algebra object of $\DZ$ with a morphism $D\to R$ to a noetherian ring viewed in degree zero, and let $A=\lim_m\delta(R/I^m)$ for an ideal $I\subseteq R$.
For $P\in\mathcal D(D)$, the $R$-action on $\Hom_D(P,R)$ is the action through its target.
Applying $\delta$ to evaluation and extending the target coefficients to $A$ induces a map
\[
  \bigl(A\otimes_{\delta R}^{\square}\delta\Hom_D(P,R)\bigr)
       \otimes_{\Z_{\square}}\delta P\longrightarrow A,
\]
which is $\delta D$-linear in the second factor and $A$-linear in the first.
Its tensor--Hom adjoint is the comparison
\begin{equation}
  \beta_P:A\otimes_{\delta R}^{\square}\delta\Hom_D(P,R)
       \longrightarrow\sHom_{\delta D}(\delta P,A).
  \label{eq:pseudocoherent-hom-base-change}
\end{equation}

\begin{lemma}
\label{lem:pseudocoherent-hom-base-change}
Let $D$ be a connective commutative algebra object of $\DZ$, equipped with a morphism of commutative algebra objects $D\to R$ to a noetherian ring $R$, viewed in degree zero.
Suppose that for each $n\geq1$ there is a finite-cell $D$-module $P_n\to R$ with cofiber in $\mathcal D(D)^{\leq-n}$.
For every ideal $I\subseteq R$ and $A=\lim_m\delta(R/I^m)$, the map $\beta_R$ of~\eqref{eq:pseudocoherent-hom-base-change} is an equivalence.
\end{lemma}

\begin{proof}
The map $\beta_D$ is the identity of $A$.
Viewed as contravariant functors of $P$, both sides of~\eqref{eq:pseudocoherent-hom-base-change} send finite colimits to finite limits and preserve retracts.
Hence $\beta_P$ is an equivalence for every finite-cell $P$, and in particular for $P=P_n$.
Let $C_n=\operatorname{cofib}(P_n\to R)$.
Truncation orthogonality implies
\[
  \Hom_D(C_n,R)\in\mathcal D(R)^{\geq n}.
\]
Its scalar extension remains in degrees at least $n$, by Lemma~\ref{lem:flat-profinite-coefficients}.
For $\sHom_{\delta D}(\delta C_n,A)$, tensor--Hom adjunction takes the following form after evaluation on a profinite set $S$:
\[
  \Hom_{\DS}\bigl(\Z[S]^{\square},
          \sHom_{\delta D}(\delta C_n,A)\bigr)
  \simeq
  \Hom_{\Mod_{\delta D}(\DS)}
       \bigl(\Z[S]^{\square}\otimes_{\Z_{\square}}\delta C_n,A\bigr).
\]
Because $\delta D$ is connective, the standard $t$-structure on $\Mod_{\delta D}(\DS)$ is detected on underlying solid complexes.
The source $\Z[S]^{\square}\otimes_{\Z_{\square}}\delta C_n$ lies in degrees at most $-n$ by right $t$-exactness of tensor.
The target $A$ is concentrated in degree zero, as noted after~\eqref{eq:adic-solid-ring}, so truncation orthogonality in this module category puts the mapping complex on the right in $\mathcal D(\Z)^{\geq n}$.
Evaluation on profinite sets detects the $t$-structure, so $\sHom_{\delta D}(\delta C_n,A)$ belongs to $\DS^{\geq n}$.
For a fixed integer $q$ and $n>q+1$, these bounds and the fiber sequences obtained from $P_n\to R\to C_n$ identify degree-$q$ cohomology of both contravariant Hom constructions at $R$ with that at $P_n$.
Naturality reduces $H^q(\beta_R)$ to $H^q(\beta_{P_n})$, which is invertible.
Nondegeneracy of the standard $t$-structure on $\DS$ proves the result.
\end{proof}

Assume now that $R$ is essentially of finite type over $\Z$, with $I\subseteq R$ and $A=\lim_n\delta(R/I^n)$ as above.
We construct the cochain map as a morphism of $E_1$ $A$-algebras.
Regard $R^e$ as an $R$-algebra through its first factor, and put
\[
  B=A\otimes_{\Z_{\square}}\delta R.
\]
The strong symmetric monoidal functor $V\mapsto A\otimes_{\delta R}^{\square}\delta V$ from $\mathcal D(R)$ to $\Mod_A(\DS)$ sends $R^e$ to $B$ and its diagonal $R$ to the $B$-module $A$.
This extension changes the coefficients of the diagonal calculation, with Hochschild base $\Z$ throughout.
Applying this functor to module evaluation and composition defines a map of endomorphism $E_1$ $A$-algebras
\[
  A\otimes_{\delta R}^{\square}\delta\Hom_{R^e}(R,R)
       \longrightarrow\sHom_B(A,A).
\]
A map $\delta R\to A$ of commutative solid algebras is a homotopy epimorphism if multiplication $A\otimes_{\delta R}^{\square}A\to A$ is invertible \cite[Definition~2.14 and Lemma~2.15]{BenBassatMukherjee}.
Proposition~\ref{thm:solid-idempotence} proves that this multiplication map is an equivalence for the adic algebra $A$.
The map $B\to A^e$ is the base change of $\delta R\to A$ along the second-factor map $\delta R\to B$, $r\mapsto1\otimes r$, since
\[
  A^e\simeq B\otimes_{\delta R}^{\square}A.
\]
Base change preserves the multiplication equivalence, so restriction along $B\to A^e$ is fully faithful and therefore induces an equivalence
\[
  \sHom_{A^e}(A,A)\longrightarrow\sHom_B(A,A)
\]
of endomorphism $E_1$ $A$-algebras.
Compose the preceding coefficient map with the inverse of this restriction equivalence to define
\begin{equation}
  \gamma:A\otimes_{\delta R}^{\square}\delta\HH^{\bullet}(R/\Z)
       \longrightarrow\HH_{\square}^{\bullet}(A/\Z).
  \label{eq:hochschild-cochain-map}
\end{equation}
The second arrow of~\eqref{eq:completed-hochschild-cochains} is $\alpha_{\delta\HH^{\bullet}(R/\Z)}$ from~\eqref{eq:adic-scalar-map}; its target is the solid derived completion $\Lambda_I\delta\HH^{\bullet}(R/\Z)$.

\begin{proposition}
\label{thm:hochschild-comparison}
\label{cor:intro-hochschild}
Let $R$ be a commutative ring essentially of finite type over $\Z$, let $I\subseteq R$ be an ideal, and put $A=\lim_n\delta(R/I^n)$.
The maps $\gamma$ of~\eqref{eq:hochschild-cochain-map} and $\alpha_{\delta\HH^{\bullet}(R/\Z)}$ of~\eqref{eq:adic-scalar-map} are equivalences in $\mathcal D(\Z_\square)$:
\begin{equation}
  \HH_{\square}^{\bullet}(A/\Z)
  \xleftarrow[\ \sim\ ]{\gamma}
       A\otimes_{\delta R}^{\square}\delta\HH^{\bullet}(R/\Z)
  \xrightarrow[\ \sim\ ]{\alpha_{\delta\HH^{\bullet}(R/\Z)}}
       \Lambda_I\delta\HH^{\bullet}(R/\Z).
  \label{eq:completed-hochschild-cochains}
\end{equation}
The cochain scalar-extension equivalence $\gamma$ respects the $E_1$ $A$-algebra structure given by the cup product.
In particular, for every integer $q$ there is an isomorphism of solid $A$-modules concentrated in degree zero:
\begin{equation}
  H^q\HH_{\square}^{\bullet}(A/\Z)
  \simeq A\otimes_{\delta R}^{\square}
            \delta\operatorname{Ext}_{R^e}^q(R,R).
  \label{eq:completed-hochschild-cohomology}
\end{equation}
\end{proposition}

\begin{proof}
Using associativity of relative tensor products and the multiplication equivalence of Proposition~\ref{thm:solid-idempotence}, we compute
\[
 A^e\otimes_{\delta R^e}^{\square}\delta R\simeq A\otimes_{\delta R}^{\square}A\xrightarrow{\sim}A.
\]
The extension--restriction adjunction identifies
\[
  \sHom_{A^e}(A,A)\simeq\sHom_{\delta R^e}(\delta R,A).
\]
Under this identification, $\gamma$ is the evaluation comparison $\beta_R$ for $D=R^e$.
Indeed, the analogous adjunction identifies $\sHom_B(A,A)$ with $\sHom_{\delta R^e}(\delta R,A)$, and the coefficient map defining $\gamma$ is adjoint to the evaluation used to define $\beta_R$.
Apply Lemma~\ref{lem:pseudocoherent-hom-base-change} to the connective commutative algebra object $D=R^e=R\otimes_{\Z}R$, with its multiplication morphism $R^e\to R$ and the finite-cell approximations of Lemma~\ref{lem:finite-diagonal-approximations}.
It proves that $\gamma$ is an equivalence, and its construction preserves the $E_1$ $A$-algebra structure.

For a fixed integer $q$, choose $n\geq1$ with $n>q+1$; the cofiber bound implies $\operatorname{Ext}_{R^e}^q(R,R)\simeq H^q\Hom_{R^e}(E_n,R)$.
Since $E_n$ is finite-cell, $\Hom_{R^e}(E_n,R)$ is a perfect $R$-complex; its cohomology is finitely generated because $R$ is noetherian.
Proposition~\ref{thm:finite-cohomology} therefore makes $\alpha_{\delta\HH^{\bullet}(R/\Z)}$ an equivalence, including the entire unbounded cochain complex.
Equation~\eqref{eq:completed-hochschild-cohomology} follows from the $t$-exactness in Lemma~\ref{lem:flat-profinite-coefficients}.
\end{proof}

Proposition~\ref{prop:quotient-towers} identifies the completed cochain complex in~\eqref{eq:completed-hochschild-cochains} with the inverse limit of the derived reductions of $\delta\HH^\bullet(R/\Z)$:
\[
  \Lambda_I\delta\HH^\bullet(R/\Z)
  \simeq\lim_n\bigl(\delta(R/I^n)\otimes_{\delta R}^{\square}\delta\HH^\bullet(R/\Z)\bigr).
\]

\section{Examples}
\label{sec:examples}

We compute endomorphism objects in three punctured categories over the relative solid base $\mathcal D((k,\Z)_\square)$.
A product modulo the direct sum exhibits the failure of shifted $k$-linear duality after point evaluation.
Our example based on a proper exterior dg algebra has boundary cohomology unbounded in both directions.
Over $\Z$, completion at $(p)$ recovers the $p$-adic integers, and the puncture exhibits the internal dual of the $p$-adic field in $\mathcal D(\Z_\square)$.

\subsection{A product modulo its direct sum}

Let $k=\F_p$, and let $\mathcal T$ be the stable idempotent-complete category generated by countably many orthogonal objects $t_i$, with $\mathcal T(t_i,t_i)=k$.
Every object of $\mathcal T$ lies in $\thick\{t_i:i\in F\}$ for some finite set $F$.
Since the mapping complexes between the generators are perfect over $k$, $\mathcal T$ is proper.
The sum of the traces on those finite sets defines a coherent $0$-Calabi--Yau structure.
Let $P$ be the pointwise perfect right module with $P(t_i)=k$ for every $i$, and take $\mathcal H$ to be the full stable idempotent-complete subcategory generated by $\delta P$ and the representables.

\begin{example}\label{ex:product-quotient}
Put
\begin{equation}\label{eq:product-quotient-algebras}
 A=\prod_{i\geq1}\delta k,
 \qquad J=\bigoplus_{i\geq1}\delta k,
 \qquad B=\cofib(J\longrightarrow A).
\end{equation}
Orthogonality identifies $\mathcal H(P,P)$ with $A$ and the coend $I(P,P)$ in~\eqref{eq:punctured-factorization-coend} with $J$.
Composition through $h_{t_i}$ sends $1\in\delta k$ to the $i$th coordinate idempotent $e_i\in A$.
Thus $\alpha_{P,P}$ is the coordinate inclusion $J\to A$, so~\eqref{eq:punctured-quotient-triangle} identifies $B$ with $\mathcal Q_0(qP,qP)$.

The idempotent $e_i$ exhibits $h_{t_i}$ as a retract of $P$, so $\mathcal H=\thick(P)$.
For a finite set $F$ of indices, put $e_F=\sum_{i\in F}e_i$.
The quotient map $A\to\mathcal Q_0(qP,qP)$ sends every $e_i$ to zero and hence inverts every $1-e_F$.
It therefore induces an algebra map from the central localization
\[
 \colim_F(1-e_F)A\longrightarrow\mathcal Q_0(qP,qP)\simeq B.
\]
The transition for $F\subseteq G$ in this colimit is $a\mapsto(1-e_G)a$.
Taking the filtered colimit of $e_FA\to A\to(1-e_F)A$ and comparing with~\eqref{eq:punctured-quotient-triangle} shows that this map is an equivalence on underlying $A$-modules, hence an equivalence of algebras.
This identifies the quotient composition on $B$ with the multiplication on the commutative solid algebra $\colim_F(1-e_F)A$.
Since $qP$ generates the quotient, enriched Morita theory identifies its module category with that of $B$.

The product $A$ and direct sum $J$ realize bounded countable pro-perfect and ind-perfect objects, respectively, so their cofiber $B$ realizes a bounded countable derived Tate object \cite[Definition~4.2 and Theorem~4.3]{MPV}.
Internal duality identifies $\sqdual J$ with $A$, and Lemma~\ref{lem:unbounded-solid-biduality} then implies $\sqdual A\simeq J$.
The coordinate evaluation pairings identify the transpose of $J\to A$ with the same inclusion.
Dualizing~\eqref{eq:product-quotient-algebras}, we therefore obtain
\[
 \sqdual B\simeq\fib(J\longrightarrow A)\simeq B[-1].
\]

Point evaluation preserves products, direct sums and cofibers by \S~\ref{subsec:solid-coefficients}.
It sends $J\to A$ to the inclusion of the direct sum into the product, so the cofiber is concentrated in degree zero with
\[
 B(*)=\left(\prod_i k\right)\bigg/\left(\bigoplus_i k\right).
\]
The constant sequence $1$ has nonzero class in this quotient.
Since $k$ is a field, the derived $k$-linear dual $\Hom_k(B(*),k)$ is concentrated in degree zero, while $B(*)[-1]$ has nonzero cohomology in degree one.
Thus point evaluation does not identify the internal solid dual with the algebraic dual, and the corresponding shifted $k$-linear self-duality fails, as in \cite[Example~5.7]{GaoGao}.
\end{example}

\subsection{An unbounded boundary from one proper generator}

\begin{example}\label{ex:unbounded-frobenius}
Let $k=\Z$ or $\F_p$, and consider the dg algebra
\[
 A_0=\Lambda_k(\epsilon,\eta),
 \qquad |\epsilon|=1,\quad |\eta|=3,\quad d_{A_0}=0.
\]
The coefficient of $\epsilon\eta$ defines a strict cyclic trace $A_0\to k[-4]$ with perfect graded pairing.
This trace equips $\mathcal T=\Perf(A_0)$ with a coherent proper $4$-Calabi--Yau structure.
Let $M=k$ be the augmentation right module, and let $\mathcal H$ be the full stable subcategory generated by $\delta M$ and the representables.
Its objects are pointwise perfect, and Corollary~\ref{cor:punctured-cy} applies.

A semifree right-module resolution of $M$ has underlying graded module $\bigoplus_{a,b\geq0}e_{a,b}A_0$, with
\begin{equation}\label{eq:exterior-resolution}
 |e_{a,b}|=2b,
 \qquad
 d(e_{a,b})=e_{a-1,b}\epsilon+e_{a,b-1}\eta,
\end{equation}
where we interpret any term with a negative index as zero.
The augmentation sends $e_{0,0}$ to $1$ and the other generators to zero.
The exterior relations on $\epsilon$ and $\eta$ imply $d^2=0$, and the filtration by $a+b$ is semifree.
In one variable of degree $e=1$ or $3$, the analogous resolution has generators in degree $a(e-1)$; its differential pairs each positive-index generator with the preceding generator multiplied by the exterior variable.
The augmentation identifies its cohomology with $k$ in degree zero.
Tensoring the two one-variable resolutions over $k$ recovers the complex in~\eqref{eq:exterior-resolution} and proves that it resolves $M$.

The strong symmetric monoidal, colimit-preserving functor $\delta$ carries this resolution to a resolution of $\delta M$ by free $\delta A_0$-module cells.
After applying internal Hom into $\delta k$, the differential vanishes because $\epsilon$ and $\eta$ act by zero, so
\begin{equation}\label{eq:unbounded-completed-hom}
 \mathcal H(M,M)\simeq\prod_{a,b\geq0}\delta k[2b].
\end{equation}
Morita invariance of the coend~\eqref{eq:punctured-factorization-coend} identifies $I(M,M)$ with $\mathcal H(M,h_{A_0})\otimes_{\delta A_0}^{\square}\delta M$.
Equation~\eqref{eq:serre-reverse-mixed} identifies the first factor with $\delta k[-4]$, with its left $\delta A_0$-action induced by the augmentation.
Tensoring the resolution with $k[-4]$ and using that $\delta$ preserves relative tensor products, we obtain
\begin{equation}\label{eq:unbounded-factorization}
 I(M,M)\simeq\delta\bigl(k[-4]\otimes_{A_0}M\bigr)
 \simeq\bigoplus_{a,b\geq0}\delta k[-2b-4],
\end{equation}
again with zero differential.

The source of the composition map $\alpha_{M,M}:I(M,M)\to\mathcal H(M,M)$ is discrete, so the coefficient adjunction computes its mapping class over $k$ after point evaluation of the target.
For each free source summand, we deduce from~\eqref{eq:unbounded-completed-hom} that
\[
 H^0\Hom_k\bigl(k[-2b-4],\mathcal H(M,M)(*)\bigr)
 \cong H^{2b+4}\bigl(\mathcal H(M,M)(*)\bigr)=0.
\]
Maps out of a direct sum are products of mapping spaces, so $\alpha_{M,M}$ is nullhomotopic.
From the quotient triangle~\eqref{eq:punctured-quotient-triangle}, we therefore deduce a noncanonical equivalence of underlying solid complexes $B=\mathcal Q_0(qM,qM)\simeq\mathcal H(M,M)\oplus I(M,M)[1]$.
Its cohomology objects are
\begin{equation}\label{eq:unbounded-example-cohomology}
 H^n(B)=
 \begin{cases}
 \displaystyle\prod_{a\geq0}\delta k,&n=-2b,\quad b\geq0,\\[3pt]
 \displaystyle\bigoplus_{a\geq0}\delta k,&n=2b+3,\quad b\geq0,\\[3pt]
 0,&\text{otherwise}.
 \end{cases}
\end{equation}
Corollary~\ref{cor:punctured-cy}, applied with $d=4$, supplies the bimodule duality $\sqdual B\simeq B[3]$.
The endomorphism complex $B$ is unbounded in both directions, and every nonzero $H^n(B(*))$ is infinitely generated over $k$.
Its unbounded cohomology places $B$ outside the image of the bounded countable derived Tate objects in~\cite[Theorem~4.3]{MPV}.
For $k=\F_p$, evaluating~\eqref{eq:unbounded-completed-hom} at a point shows that $H^0\Hom_{A_0}(k,k)\cong\prod_{a\geq0}k$, so the finite-dimensionality hypothesis in~\cite[Theorem~5.6]{GaoGao} fails.

Over $k=\F_7$, the endomorphism object $B$ above also occurs in the solid punctured completion of $\Perf(X)$ for a smooth proper scheme $X$, as follows.
Let $E/k$ be an elliptic curve, let
\[
 Y=\{z_0^5+z_1^5+z_2^5+z_3^5+z_4^5=0\}
 \subset\mathbb P^4_{\F_7},
 \qquad X=E\times Y,
\]
and choose the rational point $x=(0_E,(1:-1:0:0:0))$.
The partial derivatives $5z_i^4$ have no common projective zero, so $Y$ is smooth.
By adjunction, $\omega_Y\simeq\mathcal O_Y(5-5)\simeq\mathcal O_Y$ \cite[Tag~0E9Z]{Stacks}.
The hypersurface sequence $0\to\mathcal O_{\mathbb P^4_k}(-5)\to\mathcal O_{\mathbb P^4_k}\to\mathcal O_Y\to0$ and~\cite[Tag~01XT]{Stacks} imply $H^*(Y,\mathcal O_Y)=k\oplus k[-3]$.
For the elliptic curve, an invariant differential trivializes $\omega_E$, and the genus-one condition implies $H^*(E,\mathcal O_E)=k\oplus k[-1]$ \cite[Tag~072J]{Stacks}.
Both global-section algebras have no negative cohomology and have $H^0=k$, so we choose strictly unital minimal $A_\infty$ models as in~\cite[\S\S~3.3 and~3.5]{KellerAinf}.
For a nonunit generator in degree $e=1$ or $3$, an operation of arity $n\geq2$ has output degree $ne+2-n$, outside the respective cohomology range.
All these operations vanish; operations of arity at least three involving the unit vanish by strict unitality.
Since $E$ and $Y$ are proper over the field $k$, pullback and cup product induce the K\"unneth equivalence \cite[Tag~0FLQ]{Stacks}; using the factorwise models, we then obtain
\[
 \begin{gathered}
 \Hom_E(\mathcal O_E,\mathcal O_E)\otimes_k
 \Hom_Y(\mathcal O_Y,\mathcal O_Y)\\
 \xrightarrow{\ \sim\ }\Hom_X(\mathcal O_X,\mathcal O_X)\simeq A_0
 \end{gathered}
\]
as $E_1$-algebras.

Take $\mathcal C=\Perf(X)$ and $\mathcal T=\thick(\mathcal O_X)$.
The algebra equivalence above determines the module calculations below, while the chosen trivializations of $\omega_E$ and $\omega_Y$ make $[4]$ a coherent Serre functor on $\mathcal C$.
Since $X$ is smooth and proper of dimension four, $\mathcal C$ is proper, and $\mathcal T$ is Serre-stable.
The object $\mathcal O_x$ is perfect by smoothness.
Under the Morita equivalence $\mathcal T\simeq\Perf(A_0)$, its restricted Yoneda module is represented by $\Hom_X(\mathcal O_X,\mathcal O_x)\simeq k$.
A strictly unital minimal model of this right module is concentrated in degree zero \cite[\S~4.3(b)]{KellerAinf}.
An operation with $r\geq1$ positive-degree algebra arguments has output degree at least $r+(1-r)=1$, so it vanishes.
Thus the restricted module is the augmentation module.
The full stable subcategory generated by this module and the representables embeds in the stable idempotent-complete envelope of the solid categorical completion $\widehat{\mathcal C}_{\mathcal T}^{\square}$.
Both categories use the same full module Homs, so the mixed Homs and composition map in~\eqref{eq:punctured-factorization-coend}--\eqref{eq:punctured-quotient-triangle} are unchanged at this module.
The quotient endomorphism object of the image of $\mathcal O_x$ is therefore $B$, compatibly with composition.

If $\mathcal O_x$ belonged to $\mathcal T$, its augmentation module $k$ would be perfect over $A_0$.
Properness of $A_0$ would then make $\Hom_{A_0}(k,k)$ perfect over $k$, contradicting~\eqref{eq:unbounded-completed-hom} after point evaluation.
Hence $\mathcal O_x\notin\mathcal T$ and $\mathcal T\subsetneq\Perf(X)$.
Since $\mathcal T$ contains the tensor unit $\mathcal O_X$, it cannot be a tensor ideal: a tensor ideal containing $\mathcal O_X$ equals $\Perf(X)$.
Every $\Perf_Z(X)$ is a tensor ideal, so $\mathcal T\ne\Perf_Z(X)$ for every closed subscheme $Z\subseteq X$.
Theorem~\ref{thm:intro-punctured-serre} thus applies to a categorical completion that does not arise from a closed support.
\end{example}

\subsection{Arithmetic completion and duality}

\begin{example}\label{ex:arithmetic}
Let $p$ be prime, $k=\Z$, $\mathcal C=\Perf(\Z)$ and $\mathcal T=\Perf_{(p)}(\Z)=\thick(\Z/p)$.
The usual trace on perfect $\Z$-complexes makes the identity a coherent Serre functor.
Since $(p)$ is weakly proregular, Proposition~\ref{thm:adic-morphisms}\textup{(iii)} makes $\operatorname{ev}_{\Z}^{\square}$ in~\eqref{eq:adic-evaluation-functors} an equivalence.
It sends the representable at $\Z$ to the free rank-one module over $\underline{\Z_p}=\lim_n\delta(\Z/p^n)$ and the generator $\Z/p$ to the cofiber of $p:\underline{\Z_p}\to\underline{\Z_p}$.
For a solid $\underline{\Z_p}$-module $N$, the internal Hom from this cofiber to $N$ is $\fib(p:N\to N)$.
Thus~\eqref{eq:punctured-module-kernel} identifies the local modules as those on which $p$ acts invertibly.
Localization is extension of scalars along $\underline{\Z_p}\to\underline{\Z_p}[1/p]$.
The image of $\Z$ in the punctured category therefore corresponds to the free rank-one module over $\underline{\Z_p}[1/p]$, with this algebra as its endomorphism algebra.
For the $p$-adic topologies, $\underline{\Z_p}[1/p]\simeq\underline{\Q_p}$: every continuous map from a profinite set to $\Q_p$ has compact image contained in some $p^{-n}\Z_p$.
Applying Theorem~\ref{thm:intro-punctured-serre} to this quotient object, we obtain
\begin{equation}\label{eq:qp-duality}
 \sHom_{\Z_\square}(\underline{\Q_p},\delta\Z)\simeq\underline{\Q_p}[-1].
\end{equation}
\end{example}

\appendix
\section{Coefficient changes for internal Homs}\label{app:centralizers}

We characterize the modules for which solid coefficient change preserves internal Homs to every target.
We then compare completion indexed by $\delta\mathcal M$ with completion indexed by the full subcategory on the realized modules $\delta M$, and show that pointwise perfection alone does not make the comparison an equivalence.

\subsection{Detection of perfect sources}

Let $\mathcal B$ be a small $\mathcal D(k)$-enriched category over a discrete commutative ring $k$.
For right $\mathcal B$-modules $M,N$ and a profinite set $S$, recall the comparisons of~\eqref{eq:discreteness-comparison} and~\eqref{eq:module-coefficient-comparison}:
\[
 \begin{aligned}
 \chi_{S;M,N}:{}&
 C(S,k)\otimes_k\Hom_{\mathcal B}(M,N)
 \longrightarrow\Hom_{\mathcal B}(M,C(S,k)\otimes_kN),\\
 \kappa_{M,N}:{}&
 \delta\Hom_{\mathcal B}(M,N)
 \longrightarrow\sHom_{\delta\mathcal B}(\delta M,\delta N).
 \end{aligned}
\]
The second target is formed over $\mathcal D((k,\Z)_\square)$.
By the right-module version of~\eqref{eq:compact-enriched-modules}, perfect modules are the compact objects of $\RMod_{\mathcal B}(\mathcal D(k))$.

\begin{proposition}\label{prop:perfection-detection}
For a right $\mathcal B$-module $M$, the following are equivalent:
\begin{enumerate}
\item $M$ is perfect.
\item The map $\chi_{S;M,N}$ is an equivalence for every right $\mathcal B$-module $N$ and every profinite $S$.
\item The map $\kappa_{M,N}$ is an equivalence for every right $\mathcal B$-module $N$.
\end{enumerate}
\end{proposition}

\begin{proof}
For perfect $M$, its exact mapping-complex functor preserves colimits.
Since $C(S,k)$ is a filtered colimit of finite free $k$-modules, condition (ii) follows.
The identity $\kappa_{M,N}(S)=\theta_{S;M,N}\circ\chi_{S;M,N}$ from Proposition~\ref{thm:family-formula} and detection by the $k[S]^\square$ identify (ii) and (iii).

Conversely, for an infinite set $J$ and its profinite one-point compactification $S=J\cup\{\infty\}$, the $k$-module $C(S,k)$ decomposes as follows:
\[
 C(S,k)\cong k\oplus\bigoplus_{j\in J}k,
 \qquad
 f\longmapsto\bigl(f(\infty),(f(j)-f(\infty))_{j\in J}\bigr).
\]
Continuity at $\infty$ makes the second component finitely supported.
Under this decomposition, $\chi_{S;M,N}$ is the direct sum of the identity on $\Hom_{\mathcal B}(M,N)$ and the canonical coproduct comparison
\[
 \bigoplus_{j\in J}\Hom_{\mathcal B}(M,N)
 \longrightarrow\Hom_{\mathcal B}\left(M,\bigoplus_{j\in J}N\right).
\]
Condition (ii) makes this comparison an equivalence.
For a family $(N_j)$, take $N=\bigoplus_jN_j$ and let $p_j:N\to N_j\to N$ be the projection idempotent.
Naturality in each summand intertwines $\bigoplus_j\Hom_{\mathcal B}(M,p_j)$ on the source with $\Hom_{\mathcal B}(M,\bigoplus_jp_j)$ on the target.
Passing to their retracts, we obtain
\[
 \bigoplus_j\Hom_{\mathcal B}(M,N_j)
 \xrightarrow{\ \sim\ }\Hom_{\mathcal B}\left(M,\bigoplus_jN_j\right).
\]
Since $\Hom_{\mathcal B}(M,-)$ is exact, it now preserves all colimits by~\cite[Proposition~1.4.4.1(2)]{LurieHA}.
Thus $M$ is compact and hence perfect by~\eqref{eq:compact-enriched-modules}.
\end{proof}

Proposition~\ref{prop:perfection-detection} requires comparison with every target module; discreteness of a particular completion requires the comparisons only between its evaluation modules.
If the full inclusion $\mathcal T\hookrightarrow\mathcal C$ of \S~\ref{subsec:full-subcategory} admits an enriched right adjoint $q$, then
\[
 r_{\mathcal T}X\simeq\mathcal T(-,qX),
 \qquad
 \widehat{\mathcal C}_{\mathcal T}^{\square}(X,Y)
 \simeq\delta\mathcal T(qX,qY)
\]
by adjunction and solid Yoneda.

\subsection{Comparison of the two completions}

Let $\mathcal M$ be a full small subcategory of left modules on a small $\mathcal D(k)$-enriched category $\mathcal C$.
As in Remark~\ref{rem:fixed-testing-representation}, let $\mathcal M_{\mathrm{full}}$ be the full subcategory of $\LMod_{\delta\mathcal C}(\mathcal D((k,\Z)_\square))$ on the solid realizations $\delta M$.
The maps $\kappa_{M,N}$ in~\eqref{eq:module-coefficient-comparison} define $q:\delta\mathcal M\to\mathcal M_{\mathrm{full}}$, sending $M$ to $\delta M$.
Evaluation at $X$ defines a left $\mathcal M_{\mathrm{full}}$-module $\varepsilon_X$ with values $\varepsilon_X(\delta M)=\delta(M(X))$ and action
\[
 \sHom_{\delta\mathcal C}(\delta M,\delta N)
 \otimes_{\delta k}^{\square}\delta(M(X))
 \longrightarrow\delta(N(X)).
\]
By the definition of $\kappa$, these actions identify $q^*\varepsilon_X\simeq\delta e_X$.
Let $\widehat{\mathcal C}_{\mathcal M}^{\square,\mathrm{full}}$ be the full enriched category on the modules $\varepsilon_X$, indexed by $X\in\mathcal C$.
Restriction induces
\[
 \sHom_{\mathcal M_{\mathrm{full}}}(\varepsilon_X,\varepsilon_Y)
 \longrightarrow\sHom_{\delta\mathcal M}(\delta e_X,\delta e_Y),
\]
and hence a functor
\begin{equation}\label{eq:full-testing-comparison}
 \rho_q:\widehat{\mathcal C}_{\mathcal M}^{\square,\mathrm{full}}
 \longrightarrow\widehat{\mathcal C}_{\mathcal M}^{\square}.
\end{equation}
Restriction identifies the evaluation maps, so this functor is under $\delta\mathcal C$.
If every $M\in\mathcal M$ is perfect as a left $\mathcal C$-module, Proposition~\ref{prop:perfection-detection} applied to $\mathcal C^{\op}$ shows that $q$ and $\rho_q$ are equivalences.
Pointwise perfection alone allows the following counterexample.

\begin{proposition}\label{prop:centralizer-counterexample}
Let $k=\F_p$, let $F=\bigoplus_{n\geq0}k$, and let $\mathcal C$ be the idempotent-complete stable envelope of the $\mathcal D(k)$-enriched category with two objects labelled $0,1$, diagonal endomorphism algebras $k$, and
\[
 \mathcal C(0,1)=F,
 \qquad \mathcal C(1,0)=0.
\]
Let $S_0,S_1$ be the extensions to $\mathcal C$ of the simple left modules $(k,0)$ and $(0,k)$ on this two-object category, and let $\mathcal M$ be the full subcategory on $S_0,S_1$.
Each $S_i$ takes every object of $\mathcal C$ to a perfect $k$-complex.
Write $F^*=\Hom_k(F,k)$ for the algebraic dual, concentrated in degree zero.
The solid completion formed from the coefficient extension $\delta\mathcal M$ has mapping object
\[
 \sHom_{\delta k}(\delta(F^*),\delta k)
\]
from $0$ to $1$, with point evaluation $F^{**}$ in degree zero.
The completion formed from $\mathcal M_{\mathrm{full}}$ has mapping object $\delta F$ from $0$ to $1$.
The map $\Gamma(\rho_q(0,1))$ is the nonsurjective bidual map $F\to F^{**}$.
\end{proposition}

\begin{proof}
Restriction from $\mathcal C$ to the original two-object category identifies the module categories by Proposition~\ref{prop:testing-envelope}.
Since $\delta$ preserves finite cofibers, shifts, and retracts, the representable-generation argument in the proof of Proposition~\ref{thm:solid-enrichment} establishes the same identification after coefficient change.
These equivalences identify the chosen modules, the full solid-enriched subcategories on their realizations, and evaluation at $0,1$.
Each $S_i$ is pointwise perfect on $\mathcal C$ because its original values are $k$ or $0$ and its extension preserves finite cofibers and retracts.
We therefore compute on the two-object category.

Over either $\mathcal V=\mathcal D(k)$ or $\mathcal V=\mathcal D((k,\Z)_\square)$, let $\mathcal B$ have two objects with unit endomorphism algebras and sole off-diagonal Hom $A=\mathcal B(0,1)$.
The left representables $P_0=(\mathbf 1_{\mathcal V},A)$ and $P_1=(0,\mathbf 1_{\mathcal V})=S_1$ fit into a cofiber sequence
\[
 A\otimes_{\mathcal V}P_1\longrightarrow P_0\longrightarrow S_0,
 \qquad S_0=(\mathbf 1_{\mathcal V},0),
\]
whose first map is zero at $0$ and $\id_A$ at $1$.
By enriched Yoneda, we have
\[
 \sHom_{\mathcal B}(P_0,S_1)=0,
 \qquad
 \sHom_{\mathcal B}(A\otimes_{\mathcal V}P_1,S_1)
 \simeq\sHom_{\mathcal V}(A,\mathbf 1_{\mathcal V}),
\]
so applying internal Hom to the cofiber sequence yields
\begin{equation}\label{eq:two-object-simple-hom}
 \sHom_{\mathcal B}(S_0,S_1)\simeq\sHom_{\mathcal V}(A,\mathbf 1_{\mathcal V})[-1].
\end{equation}
The diagonals remain the unit and the reverse Hom vanishes, so this arrow determines the full enriched category on the simples.
Evaluation at the original objects $0,1$ defines left modules over this full subcategory with values
\[
 \varepsilon_0=(\mathbf 1_{\mathcal V},0),
 \qquad
 \varepsilon_1=(0,\mathbf 1_{\mathcal V}).
\]
Their off-diagonal actions are zero, so they are its two simple modules.
Applying~\eqref{eq:two-object-simple-hom} a second time, we compute the mapping object as
\[
 \begin{aligned}
 &\sHom_{\mathcal V}
   \bigl(\sHom_{\mathcal V}(A,\mathbf 1_{\mathcal V})[-1],\mathbf 1_{\mathcal V}\bigr)[-1]\\
 &\qquad\simeq\sHom_{\mathcal V}
   \bigl(\sHom_{\mathcal V}(A,\mathbf 1_{\mathcal V}),\mathbf 1_{\mathcal V}\bigr).
 \end{aligned}
\]
For a map $u:A\to A'$, the corresponding functor between two-object categories is the identity on objects.
After restricting the cofiber sequence for $A'$ to the category with arrow $A$, the map of cofiber sequences is $u\otimes\id$ on the left, $(\id,u)$ in the middle, and $\id$ on $S_0$.
Consequently, under~\eqref{eq:two-object-simple-hom}, restriction induces precomposition with $u$:
\[
 \sHom_{\mathcal V}(A',\mathbf 1_{\mathcal V})[-1]
 \longrightarrow\sHom_{\mathcal V}(A,\mathbf 1_{\mathcal V})[-1].
\]

Over $\mathcal D(k)$,~\eqref{eq:two-object-simple-hom} with $A=F$ computes the mapping complex from $S_0$ to $S_1$ as $F^*[-1]$.
Thus $\delta\mathcal M$ has arrow $\delta(F^*)[-1]$, so the mapping object from $0$ to $1$ in $\widehat{\mathcal C}_{\mathcal M}^{\square}$ is $\sHom_{\delta k}(\delta(F^*),\delta k)$, with point evaluation $F^{**}$ by~\eqref{eq:module-point-comparison}.
For $\mathcal M_{\mathrm{full}}$, the arrow is
\[
 \sHom_{\delta k}(\delta F,\delta k)[-1]
 \simeq\left(\prod_{n\geq0}\delta k\right)[-1].
\]
Canonical evaluation is an equivalence by Lemma~\ref{lem:unbounded-solid-biduality}:
\[
 \iota:\delta F\xrightarrow{\ \sim\ }\sHom_{\delta k}\left(\prod_{n\geq0}\delta k,\delta k\right).
\]
The coordinate projections $F^*=\prod_{n\geq0}k\to k$ induce the unshifted map underlying $q$ on arrows,
\[
 \delta(F^*)\longrightarrow\prod_{n\geq0}\delta k,
\]
whose point evaluation is $\id_{F^*}$.
By the naturality of~\eqref{eq:two-object-simple-hom}, $\rho_q(0,1)$ is precomposition with this map:
\[
 \sHom_{\delta k}\left(\prod_{n\geq0}\delta k,\delta k\right)
 \longrightarrow\sHom_{\delta k}(\delta(F^*),\delta k).
\]
After $\iota$, its point evaluation sends $x\in F$ to the functional $\lambda\mapsto\lambda(x)$ in $F^{**}$.

To see nonsurjectivity, choose a nonprincipal ultrafilter on $\N$.
Since $k$ is finite, ultrafilter limits define a $k$-linear functional $\prod_{n\geq0}k\to k$ that vanishes on finite-support sequences and takes the constant sequence $1$ to $1$.
If evaluation against $x\in F$ vanished on finite-support sequences, evaluating on each coordinate vector would imply $x=0$, so it would also vanish on the constant sequence $1$.
Thus $\rho_q$ is not fully faithful, although both modules are pointwise perfect.
\end{proof}

\bibliographystyle{amsalpha}
\bibliography{solidcomplete}
\enlargethispage{\baselineskip}
\end{document}